\documentclass[11pt,letterpaper]{article}

\usepackage{amsthm}
\usepackage{amsmath}
\usepackage{amsfonts}

\usepackage[dvips]{graphicx}
\def\be{\begin{equation}}
\def\ee{\end{equation}}
\def\bea{\begin{eqnarray}}
\def\eea{\end{eqnarray}}
\def\ea{\end{array}}
\def\ds{\displaystyle}

\newcommand{\rline}  {{\RR}}
\newcommand{\rarrow} {{\,\rightarrow\,}}
\def\RR{{\rm I~\hspace{-1.15ex}R} }
\def\11{{\rm 1~\hspace{-1.5ex}1} }

\def\CC{\rm \hbox{C\kern-.56em\raise.4ex
         \hbox{$\scriptscriptstyle |$}\kern+0.5 em }}

\newcommand{\rfb}[1]{\mbox{\rm
   (\ref{#1})}\ifx\undefined\stillediting\else:\fbox{$#1$}\fi}

\newcommand{\m}     {{\hbox{\hskip 1pt}}}
\makeatletter
\def\section{\@startsection {section}{1}{\z@}{-3.5ex plus -1ex minus
    -.2ex}{2.3ex plus .2ex}{\large\bf}}
\makeatother

\newcommand{\FORALL} {{\hbox{$\hskip 11mm \forall \;$}}}

\newcommand{\Dscr} {{\cal D}}
\newcommand{\BE}{\begin{equation}}

\newcommand{\EEQ}{\end{equation}}

\newcommand{\half}   {{\frac{1}{2}}}

\font\eufm=eufm10\font\eufms=eufm10\font\eufmss=eufm10\newfam\eufam
\textfont\eufam=\eufm\scriptfont\eufam=\eufms\scriptscriptfont\eufam=\eufmss

\input{epsf}

\newtheorem{theorem}{Theorem}[section]

\newtheorem{proposition}{Proposition}[section]
\newtheorem{definition}{Definition}[section]

\begin{document}
\thispagestyle{empty}
\title{\bf The best decay rate of the damped plate equation in a square}
\author{Ka\"{\i}s Ammari  
\thanks{UR Analysis and Control of Pde, UR 13ES64, Department of Mathematics, Faculty of Sciences of
Monastir, University of Monastir, 5019 Monastir, Tunisia, \ email: kais.ammari@fsm.rnu.tn}    
 $\;\;$  Abdelkader Sa\"{\i}di
\thanks{ Institut de
Recherche Math\'ematique Avanc\'ee, University of Strasbourg, 7 rue Ren\'e Descartes,
F-67084 Strasbourg, France, \ email: saidi@math.unistra.fr}
}
\date{\today}
\maketitle
%
%
%
{\bf Abstract.~}
{\small  
In this paper we study the best decay rate of the solutions of a damped plate  equation in a square and with a homogeneous Dirichlet boundary conditions. We show that the fastest decay rate is given by the supremum of the real part of the spectrum of the infinitesimal generator of the underlying semigroup, if the damping coefficient is in $L^\infty(\Omega).$ Moreover, we give some numerical illustrations by spectral computation of the spectrum associated to the damped plate equation. The numerical results obtained for various cases of damping are in a good agreement with theoretical ones. Computation of the spectrum and energy of discrete solution of damped  plate show that the best decay rate is given by spectral abscissa of numerical solution.}
\vskip2mm
\begin{center}
\begin{minipage}[c]{3cm}
\ \hrule
\end{minipage}
\end{center}
\section{Introduction} \label{formulare}
Let $\Omega = (0,1) \times (0,1)\subset \RR^2$, $\partial
\Omega = \Gamma $. We consider the plate equation with interior dissipation, more precisely we have the following partial differential equations :
\be
\label{eq1} \frac{\partial^2 u}{\partial t^2}(x,t) +
\Delta^2 u(x,t) + a(x) \frac{\partial u}{\partial t}(x,t) = 0, \, (x,t) \in \Omega \times (0,+
\infty), \ee
with boundary conditions :
\be
\label{eq2} u(x,t) = 0, \, \Delta u(x,t) = 0, \,
(x,t) \in \Gamma \times (0,+ \infty),  \ee
and initial conditions :
\be
\label{eq3} u(x,0) = u_0(x), \,\,\,
\frac{\partial u}{\partial t}(x,0) = u_1(x) , \;\;\; x \in   \Omega,
\ee
where $a(x) \in L^{\infty}(\Omega)$ is a nonnegative damping coefficient. 

If $u$ is a solution of \rfb{eq1}-\rfb{eq3}
we define the energy of $u$ at instant $t$ 
by :
\be
\label{enercont}
 E(t) = \frac{1}{2} \int_{\Omega} \Bigr(  \left|\frac{\partial u}{\partial t} \right|^2 
+ \left| \Delta u \right|^2  \Bigr) dx.
\ee
Simple calculations show that a sufficiently smooth solution of 
\rfb{eq1}-\rfb{eq3} satisfies 
\be
\label{ESTEN}
E(t)-E(0) = - \int^t_{0} \int_{\omega} a(x) \left| \frac{\partial u}{\partial s} \right|^2  ds dx .\; \qquad \forall t \geq 0.
\ee
In particular \rfb{ESTEN} implies that 
\be
\label{dissip}
 E(t)\le E(0),\ \forall t\ge 0.
 \ee
Estimate above suggests that the natural wellposedness space 
for (\ref{eq1})-(\ref{eq3})  
is $ X = V \times L^2(\Omega),$ with inner product $\left\langle \left[f,g\right],\left[u,v\right]\right\rangle = 
\int_\Omega \left( \Delta f \, \Delta \bar{u} + g \bar{v} \right) \, dx,$ where $V = H^2(\Omega) \cap H^1_0(\Omega)$.

We have
the following wellposedness result
\begin{proposition}\label{exist}
Suppose that $(u^0,u^1)\in  V \times L^2(\Omega)$. Then the  
 problem \rfb{eq1}-\rfb{eq3}  
admits a unique solution  
$$u\in C(0,+\infty;V)\cap C^1(0,+\infty;L^2(\Omega))$$ 
Moreover $u$ satisfies the energy estimate \rfb{ESTEN}. 
\end{proposition}
If we denote by $U = \left[ u, \partial_t u \right]$, we can rewrite \rfb{eq1}-\rfb{eq3} in the form : 
\be
\label{1eq2}
\left\{
\begin{array}{ll}
\partial_t U - {\cal A}_a U = 0, \, Q = \Omega \times (0,+\infty),\\
U(x,0) = U^0(x), \, \Omega,
\end{array}
\right.
\ee
where $U^0 = (u^0,u^1), \, {\cal A}_a : {\cal D}({\cal A}_a) \subset X \rightarrow X$ is the operator defined by :
\begin{equation}
\label{Aa}
{\cal A}_a = \left(
\begin{array}{cc}
0 &  Id \\
 - \Delta^2 & - a(x) 
\end{array}
\right), 
\end{equation}
$$
{\cal D}({\cal A}_a) = \left\{ (u,v) \in [H^4 (\Omega) \cap V] \times V, \ds \Delta u_{|\Gamma} = 0 \right\}.
$$
In order to state the result on the optimal location of the actuator, 
we define the decay rate, depending on $a$, as 
$$ 
\omega(a)=\inf\{\omega\vert\hbox{ there exists }C=C(\omega)>0\hbox{ such that } 
E(t)\le C(\omega)e^{2\omega t}E(0), 
$$ 
\be 
\hbox{ for every solution of \rfb{eq1}-\rfb{eq3} with initial data 
in }V\times L^2(\Omega)\}, 
\label{DEFRATE}\ee 
and the spectral abscissa of ${\cal A}_a$,
\be
\label{as}
\mu(a) = \sup \left\{Re \lambda \, : \, \lambda \in \sigma({\cal A}_a) \right\},
\ee
where $E(t)$ is defined in \rfb{enercont} and $\sigma({\cal A}_a)$ denotes the spectrum of ${\cal A}_a$. It follows easily that 
\be
\label{BORNEINF}
\mu(a) \leq w(a).
\ee

According to \rfb{dissip} we have that $\omega(a)\le 0$
for all $a \in L^\infty(\Omega)$ is nonnegative. Moreover, if
$a \in L^\infty(\Omega)$ is nonnegative and satisfying the following condition:
\be
\label{hyp}
\exists \, c > 0 \; s.t., \; a(x) \geq c; \; a.e., \; \hbox{in an open subset} \; \omega \subset \Omega, \; meas (\omega) \neq 0.
\ee
We have, according to \cite{haraux,amm, takeo} (see Section \ref{backg}) for more details) that $\omega (a) < 0$.

The main result, on the optimal decay rate, is 

\begin{theorem} \label{locopt}
If $a \in L^\infty(\Omega)$ then
\be
\label{CONDUN} 
\mu(a) = w(a).
\ee
Moreover, if the assumption in damping coefficient \rfb{hyp} is holds, then all finite energy solutions of \rfb{eq1}-\rfb{eq3} are exponentially stable which implies that the fastest decay rate of the solutions of \rfb{eq1}-\rfb{eq3} satisfies \rfb{CONDUN}.
\end{theorem}

The problem of finding the optimal decay rate for
beams with  distributed interior damping is
difficult and has not a complete answer in the case of a variable (in space)
damping coefficient. We refer to 
\cite{note}, \cite{AHT}, \cite{AS}, \cite{aschlebeau}, \cite{lebeau}, \cite{Coxz}, \cite{Cox}, \cite{Freitas} and to references therein. 
Recently C. Castro and S. Cox in \cite{castro} made a decisive contribution
by showing that one can get an arbitrarily large decay rate by means of appropriate damping. By this way, they answer by the negative to an old
 conjecture according to which the best decay rate
 should be provided by the best constant damping.
The main  novelties brought in by this paper is that,
in the simpler case of a distributed interior damping, we can give the precise optimal decay rate and we illustrate this result numerically. One of the main ingredients of this study is a result showing that the eigenfunctions of the associated dissipative operator form a Riesz basis with parentheses in the energy space.

We remark that the corresponding optimal decay rate problem make sense since, in this case, the system is exponentially stable (see for instance 
\cite{haraux}, \cite{amm}).
 
The paper is organized as follows. Section \ref{backg} 
contains some background on optimal decay rate of dissipative systems. In section{theorique} we give the proof of the main result. 
Section \ref{discretization} is devoted to some illustrations of the main result an optimal decay rate by numerical spectral analysis strategy. We present numerical results of computation of the spectrum for different cases of damping. The discrete energy of solution is computed in each case and compared to the spectral abscissa. The latest is showed to be the best decay rate in each example presented.

\section{Some background on optimal decay rate for dissipative operators} \label{backg}
Let $H$ be a Hilbert space equipped with the norm $||.||_H$, and let
$A :\Dscr(A)  \subset H \rightarrow H$ be self-adjoint, positive and with compact invertible operator. Then, $A$ has a discrete spectrum, let $(\mu_k)_{k\geq 1}$ be its eigenvalues, each taken with its multiplicity. Denote a complete orthonormal system of eigenvectors of the operator $A$ that correspond to these eigenvalues by $(\varphi_k)_{k \geq 1}$. Moreover, we suppose that $(\mu_k)_{k\geq 1}$ satisfies the following generalized uniform gap:
\be
\label{gapcond}
\exists \, p \in \mathbb{N}^*, \, c > 0; \; \hbox{such that} \;  \mu_{k + p} - \mu_k \geq c, \, \forall \,k \in \mathbb{N}^*.
\ee
We introduce the scale of Hilbert
spaces $H_{\beta}$, $\beta\in\rline$, as follows\m: for every
$\beta \geq 0$, $H_{\beta}=\Dscr(A^{\beta})$, with the norm
$\|z \|_\beta=\|A^\beta z\|_H$. The space $H_{-\beta}$ is
defined by duality with respect to the pivot space $H$ as
follows\m: $H_{-\beta} =H_{\beta}^*$ for $\beta>0$. The
operator $A$ can be extended (or restricted) to each $H_\beta$,
such that it becomes a bounded operator
\be \label{A0extb}
A : H_\beta \rarrow H_{\beta-1} \FORALL \beta\in\rline \m.
\ee
Let a bounded
linear operator $B: U \rightarrow H$, where $U$ is another
Hilbert space which will be identified with its dual.

The systems we consider are described by
\be 
\label{dampedd1}
\ddot{x}(t) + A x(t) + \m BB^* \dot{x}(t) = 0 \m, 
\ee
\be 
\label{dampedd2}
x(0) \m=\m x_0, \, \dot{x}(0) \m=\m x_1.
\ee

We can rewrite the system \rfb{dampedd1}-\rfb{dampedd2} as a first order differential equation, by putting $z(t) = \left( \begin{array}{ll} x(t)  \\ \dot{x}(t) \end{array} \right)$ :
\be \label{dampedcc1}
\dot{z}(t) - {\cal A} z(t) + \m {\cal B}{\cal B}^* z(t)=0 \m, z(0) \m=\m z_0 = \left( \begin{array}{ll} x_0  \\ x_1 \end{array} \right), 
\ee
where 
$$
{\cal A} = \left(
\begin{array}{cc}
0 & I \\
- A & 0
\end{array}
\right) : {\cal D}({\cal A}) = H_1 \times H_\half \subset {\cal H} = H_\half \times H \rightarrow {\cal H}, \, 
{\cal B} = \left( \begin{array}{ll} 0  \\ B \end{array} \right) \in {\cal L}(U,{\cal H}).
$$
By the same way the system \rfb{dampedd1}-\rfb{dampedd2} can be also rewritten by : 
\be \label{dampedd}
\dot{z}(t) + {\cal A}_d z(t) = 0 \m, z(0) \m=\m z_0,\, 
\EEQ
where $$
{\cal A}_d = {\cal A} - {\cal B}{\cal B}^* : {\cal D}({\cal A}_d) =  H_1 \times H_\half \subset {\cal H} \rightarrow {\cal H}.$$

It is clear that the operator ${\cal A}$ is skew-adjoint on ${\cal H}$ and hence, 
it generates a strongly continous group of unitary operators on ${\cal H}$, 
denoted by $({\cal S}(t))_{t \in \rline}$.
 
Since ${\cal A}_d$ is dissipative and onto, it generates a contraction semigroup on ${\cal H}$, denoted by $({\cal S}_d (t))_{t \in \rline^+}.$ 

The system \rfb{dampedd1}-\rfb{dampedd2} is well-posed. More
precisely, the following classical result, holds.
\begin{proposition}\label{existb}
Suppose that $(x_0,x_1) \in H_\half \times H$. Then the
problem \rfb{dampedd1}-\rfb{dampedd2} admits a unique solution
$$(x,\dot{x}) \in C([0,\infty);H_\half \times H).$$
Moreover $w$ satisfies, for all $t\geq 0$, the energy estimate
\be
E(0) - E(t)  \m=\m  \int_0^t
\left\|B^* \dot{x}(s)\right\|_{U}^2, \, ds, \, 
\label{ESTENb}\ee
where $E(t) = \frac{1}{2} \, 
\|(x(t),\dot{x}(t))\|^2_{{\cal H}}$.
\end{proposition}
From \rfb{ESTENb} it follows that the mapping $t\mapsto
\|(x(t),\dot{x}(t))\|^2_{{\cal H}}$ is non
increasing. In many applications it is important to know if this
mapping decays exponentially when $t\to\infty$, i.e. if the
system \rfb{dampedd1}-\rfb{dampedd2} is exponentially stable. One
of the methods currently used for proving such exponential
stability results is based on an observability inequality for the
undamped system associated to the initial value problem
\be  
\ddot\phi(t) + A \phi(t) = 0, 
\label{eq3b}
\ee
\be 
\phi(0) = x_0, \, \dot{\phi}(0) = x_1. 
\label{eq4bb}
\ee
It is well known that \rfb{eq3}-\rfb{eq4bb} is well-posed in $H_1 \times H_\half$ and in ${\cal H}$. The result below, proved in \cite{haraux,amm}, shows that the exponential stability of \rfb{dampedd1}-\rfb{dampedd2} is equivalent to an observability inequality for \rfb{eq3b}-\rfb{eq4bb}. 
\begin{theorem}  \label{princx}
The system described by \rfb{dampedd1}-\rfb{dampedd2} is
exponentially stable in ${\cal H}$ if and only if
there exists $T,C_T > 0$ such that
\be
C_T \, \int_{0}^{T} ||{\cal B}^* {\cal S}(t) z_0||^2_{U} \, dt \geq 
||z_0||^2_{{\cal H}} \FORALL
z_0 \in {\cal H}_1.
\label{CONDUNb}
\ee
\end{theorem}

In order to state the result on the optimal decay rate, 
we define the decay rate, depending on  $B$, as 
$$ 
\omega(B)=\inf\{\omega\vert\hbox{ there exists }C=C(\omega)>0\hbox{ such that } 
E(t)\le C(\omega) \, e^{2\omega t}E(0), 
$$ 
\be 
\hbox{ for every solution of \rfb{dampedd1}-\rfb{dampedd2} with initial data  
in } H_\half \times H\}
\label{DEFRATEb}
\ee 
and the spectral abscissa as
\be
\label{asb}
\mu(B) = \sup \left\{Re \lambda \, : \, \lambda \in \sigma({\cal A}_d) \right\},
\ee
where $\sigma({\cal A}_d)$ denotes the spectrum of ${\cal A}_d$.

It follows easily that 
\begin{equation}\label{ineqinvab}
\mu(B)\leq \omega(B).
\end{equation}

We recall that a Riesz basis in a Hilbert space, is by definition, isomorphic to an orthonormal basis. 

\begin{definition} \label{RP} 
A system $(\phi_k)_{k\geq 1}$ of a space $H$ is called a basis with parentheses if the series $f = \ds \sum_{k \geq 1} c_k \phi_k$ converges in the norm of $H$ for any $f \in H$ after some arrangement of parentheses
that does not depend on $f$. If a system remains a basis after any permutation of the sets of its vectors corresponding to the terms of the series enclosed in parentheses, then such a system is called a Riesz basis with parentheses.
\end{definition}

The operator ${\cal A}_d$ is a perturbed self-adjoint operator, with the  self-adjoint part is with discrete spectrum which satisfies the gap condition \rfb{gapcond} and the perturbation, ${\cal B}{\cal B}^* \in {\cal L}({\cal H})$ is bounded. Then, according to \cite[Theorem 2]{shk}, we have that the eigenvectors of ${\cal A}_d$ forms a Riesz basis with parentheses and according to \cite{adz} we obtain on estimation of optimal decay rate. We have the following:

\begin{theorem}(\cite[Ammari-Dimassi-Zerzeri]{adz}) \label{princ}
If the observability inequality \rfb{CONDUNb} is holds then, 
\be
\label{estopt}
\omega(B) =  \mu (B) <0.
\ee
In other words if all finite energy solutions of \rfb{dampedd1}-\rfb{dampedd2} are 
exponentially stable then the fastest decay rate 
of the solutions of \rfb{dampedd1}-\rfb{dampedd2} satisfies \rfb{estopt}. 
\end{theorem}

\section{Proof of Theorem \ref{locopt}} \label{theorique}

The eigenvalue problem for the non self-adjoint, quadratic operator pencil generated by \rfb{eq1}-\rfb{eq3} is obtained by replacing $u$ in \rfb{eq1} by
$$
u(x,t) = e^{\lambda t} \phi(x). 
$$
We obtain from \rfb{eq1} the standard form
$$
({\cal A}_a - \lambda Id) \Phi = 0; \Phi = [\phi,\lambda\phi]=
\phi[1,\lambda].
$$
The condition for the existence of non trivial solutions is that $\lambda \in \sigma ({\cal A}_a)$ (the spectrum of
${\cal A}_a$). Since ${\cal D}({\cal A}_a)$ is compactly embedded in the energy space $V \times L^2(\Omega)$ then the spectrum
$\sigma({\cal A}_a)$ is discrete and the eigenvalues of ${\cal A}_a$ have a finite algebraic multiplicity. On the
other hand, since ${\cal A}_a$ is a bounded monotone perturbation of a skew-adjoint operator
(undamped ${\cal A}_0$), it follows from the Hille-Yosida theorem that ${\cal A}_a$ generates a $C_0$-semigroup
of contractions on the energy space $V \times L^2(\Omega)$.

According to \cite[Proposition A.1]{takeo} the spectra of ${\cal A}_a$ satisfies the gap assumption \rfb{gapcond}. So, by Theorem \ref{princ} we end the proof.


\hfill$\square$

%
\section{Discretization of the eigenvalue problem} \label{specanal}  \label{discretization}
The domain $\Omega$ is approximated by a net of equidistant discrete points $\Omega_h = M_{ij}$ with $i,j = 1,2,\dots,N$ and $N =(1/h) + 1$. The Laplacian operator is approximated in a standard way by a second order centered difference scheme : \\
\be
\Delta_h u = ( u_{i+1,j} + u_{i-1,j} + u_{i,j+1} + u_{i,j-1} -4 u_{i,j} )/h^2
\ee
where $u_{i,j} = u(x_i,y_j)$. The discrete approximation of the operator ${\cal A}_a$ defined by (\ref{Aa}) is then given by a nonsymmetric block matrix, $A_h$, of the form :
\be
A_h = \left[ \begin{array}{ll}
0 & -Id_N  \\
\Delta^2_h & \;\;\; a_h
\end{array}\right]
\ee
where $Id_N$ is the identity matrix of order $N$ and $\Delta^2_h$ is the discrete version of the bilaplacian.
$a_h$ is a diagonal matrix with $a_{kk} = a(x_i,y_j)$. The value of $k$ is determined by the numbering of the discrete
points $M_{ij}$. We use the implicitly restarted Arnodi-Lanszos method of Sorenson \cite{sorenson} to compute the
eigenvalues of matrix $A_h$. This method is a generalization of an inverse power method with subspace iteration 
\cite{aschlebeau},\cite{sorenson}.
\section{The case of square}
%
%
%
%
We consider a domain $\Omega=(0,1) \times (0,1)$ with a set of normalized eigenfunctions $\Phi_{n,m}(x,y) = \sqrt 2 sin(n \pi x) \sqrt 2 sin(m \pi y)$,
for all $(x,y) \in \Omega $ the solution $u(x,y,t)$
of the problem (\ref{eq1})-(\ref{eq2}) is given by :
\be \label{eq4v}
u(x,y,t)=\sum_{n,m} \alpha_{n,m}(t) \Phi_{n,m}(x,y)
\ee
replacing (\ref{eq4}) in equation (\ref{eq1}) and multiplying by a test eigenfunction  $\Phi_{k,l}(x,y)$
the solution for the low frequencies is solution of the equation :
\be
\frac{d^2 \alpha_{n,m}(t)}{dt^2} +  ( (n \pi )^2 + (m \pi )^2 )^2 \alpha_{n,m} (t) + \frac{d \alpha_{n,m}(t)}{dt} = 0 
\ee
In the first example we consider a damped plate with constant coefficient $a(x)=1$ and a 
damping region covering all the domain $\omega=\Omega$. We compute the spectrum of the damped plate (figure 1) and the energy of the first eigenmodes : n=m=3 and n=m=12 (figure 2 ). The energy of the solution (line 1) is compared to $E_0 e^{Re(\lambda_0) t}$  (line 2), where $E_0 = E(0)$ is the energy of initial conditions and $Re(\lambda_0) = inf \{Re(\lambda)/ \lambda \in \sigma(A) \}$, ($\sigma(A)$ : the spectrum of operator A). \\
In the other examples we consider a damping in different domains $\omega$ . The spectrum of the damped plate is
computed and energy for different eigenmodes is compared to  $E_0 e^{Re(\lambda_0) t}$ (line 2 : dashed line).
\vspace{2cm}
\begin{center}
\includegraphics[angle=0,width=6.5cm]{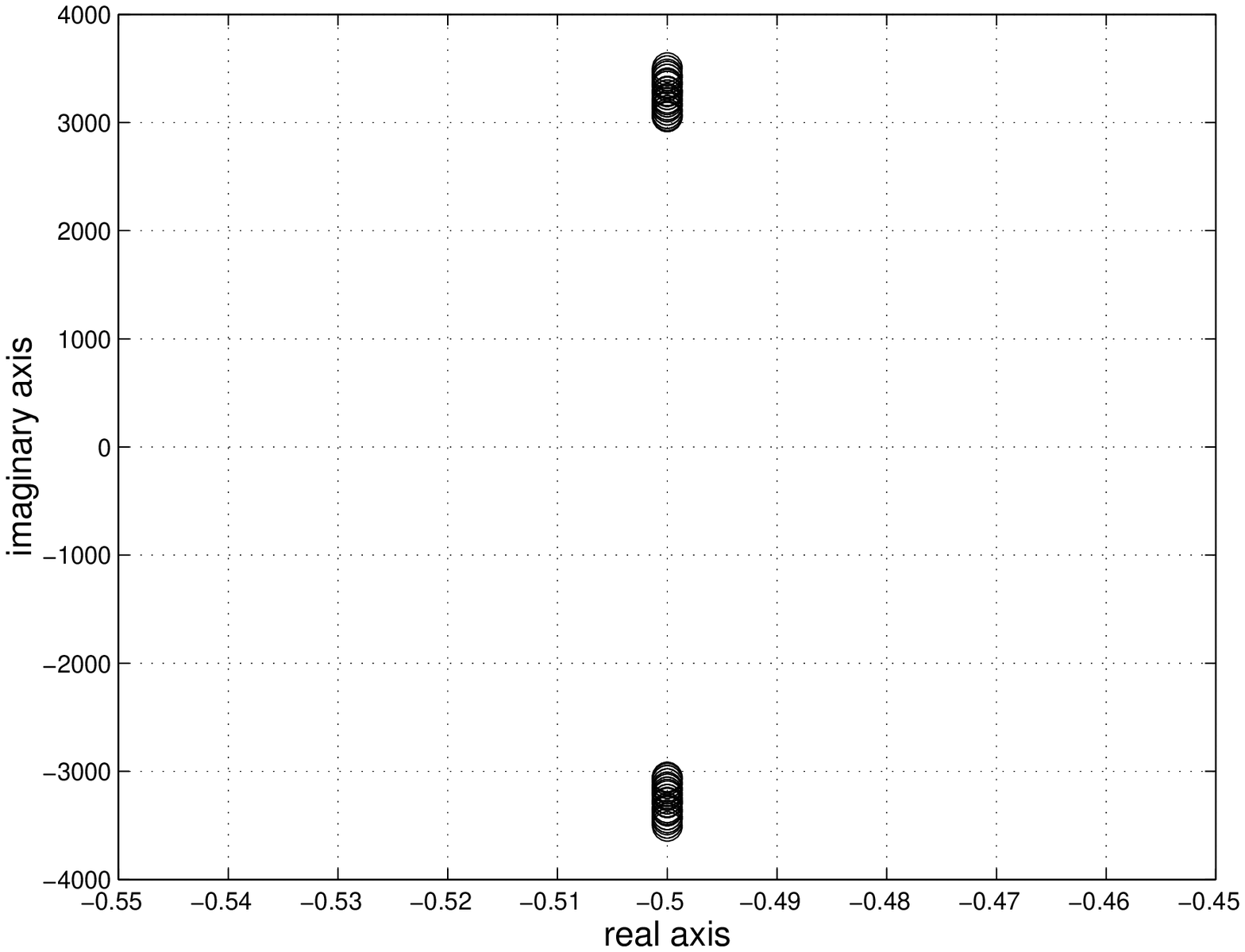}
\includegraphics[angle=0,width=6.5cm]{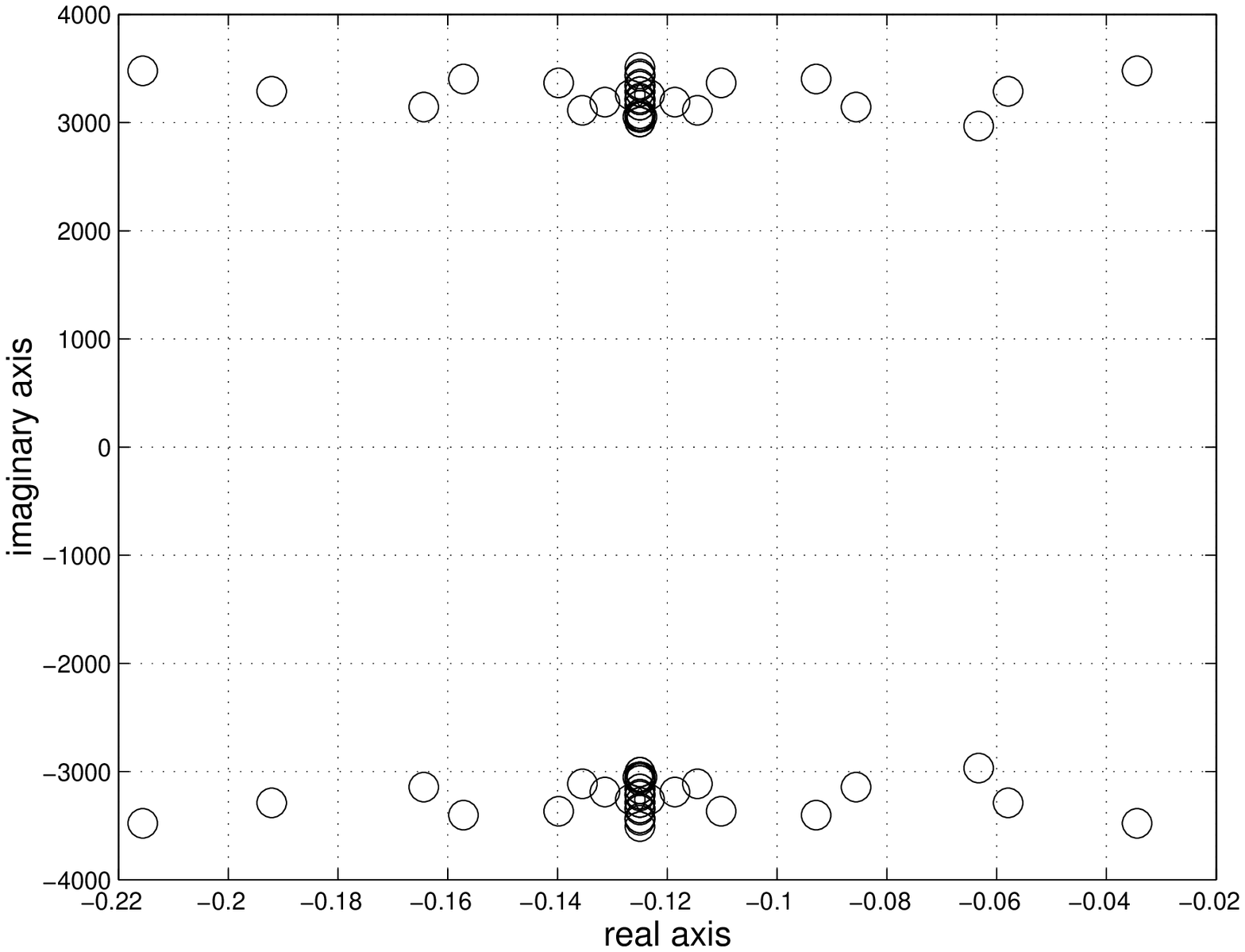}
\end{center}
\begin{center}
figure 1 : spectrum of the damped plate $a(x)=1$ , $ \omega = \Omega  $ (left)  , \\
$a(x)=1$ , $\omega = (0 , \frac{1}{2}) \times (0, \frac{1}{2}) $ (right).
\end{center}

\begin{center}
\includegraphics[angle=0,width=6.5cm]{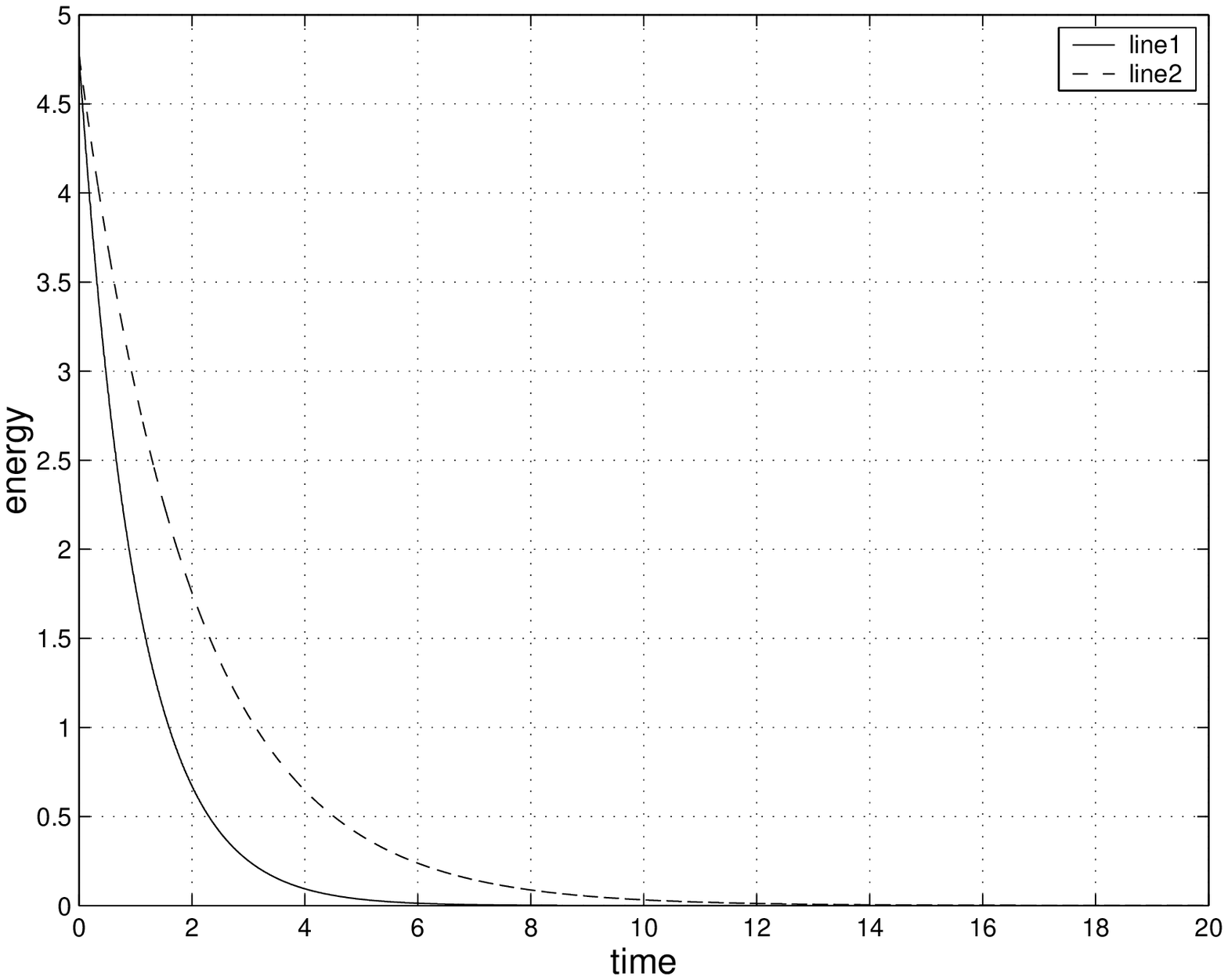}
\includegraphics[angle=0,width=6.5cm]{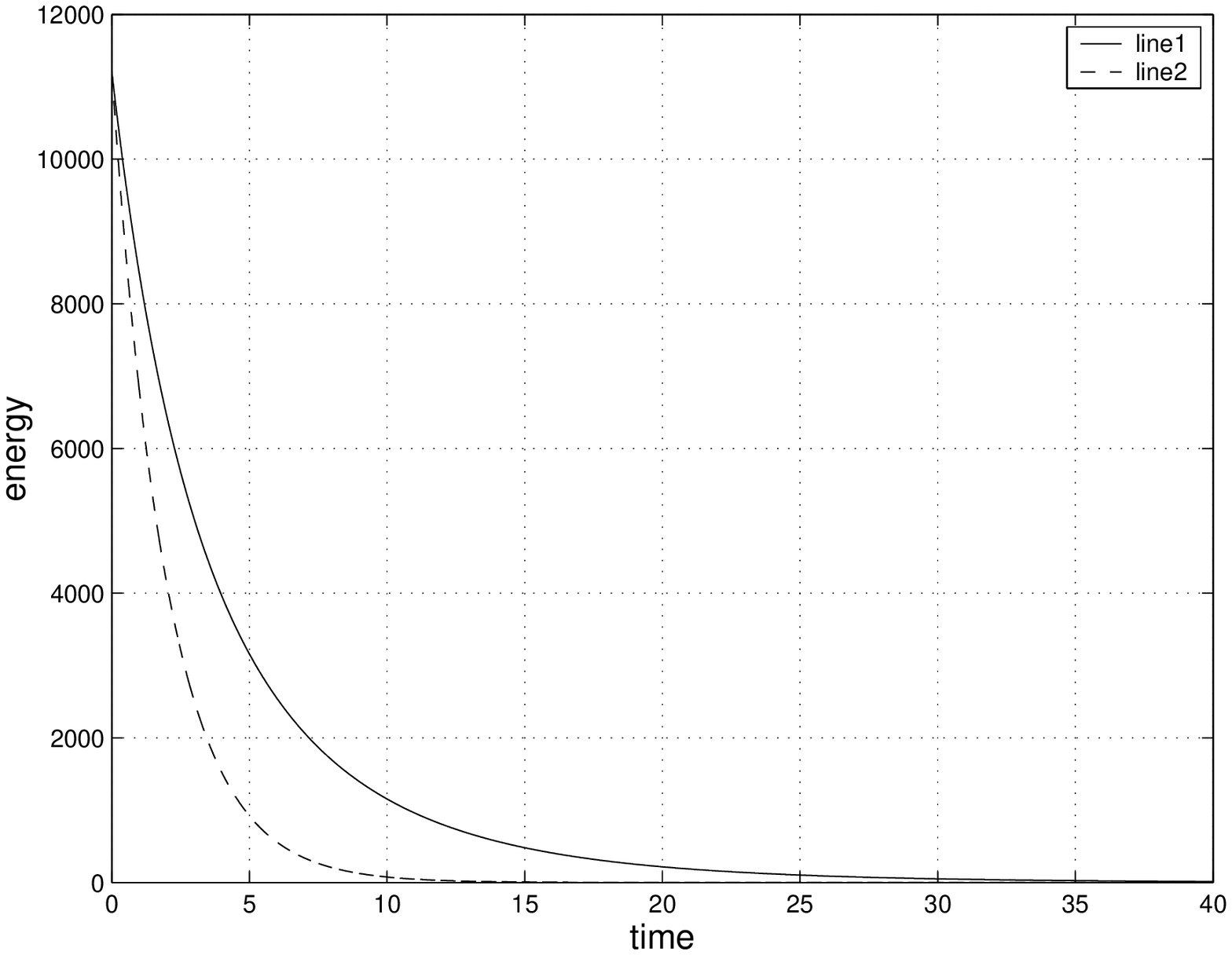}
\end{center}
\begin{center}
figure 2 : energy of the damped plate,  $a(x)=1$ , $\omega = \Omega  $, $n=m=3$ (left),  and $n=m=12$ (right).
\end{center}
\begin{center}
\includegraphics[angle=0,width=6.5cm]{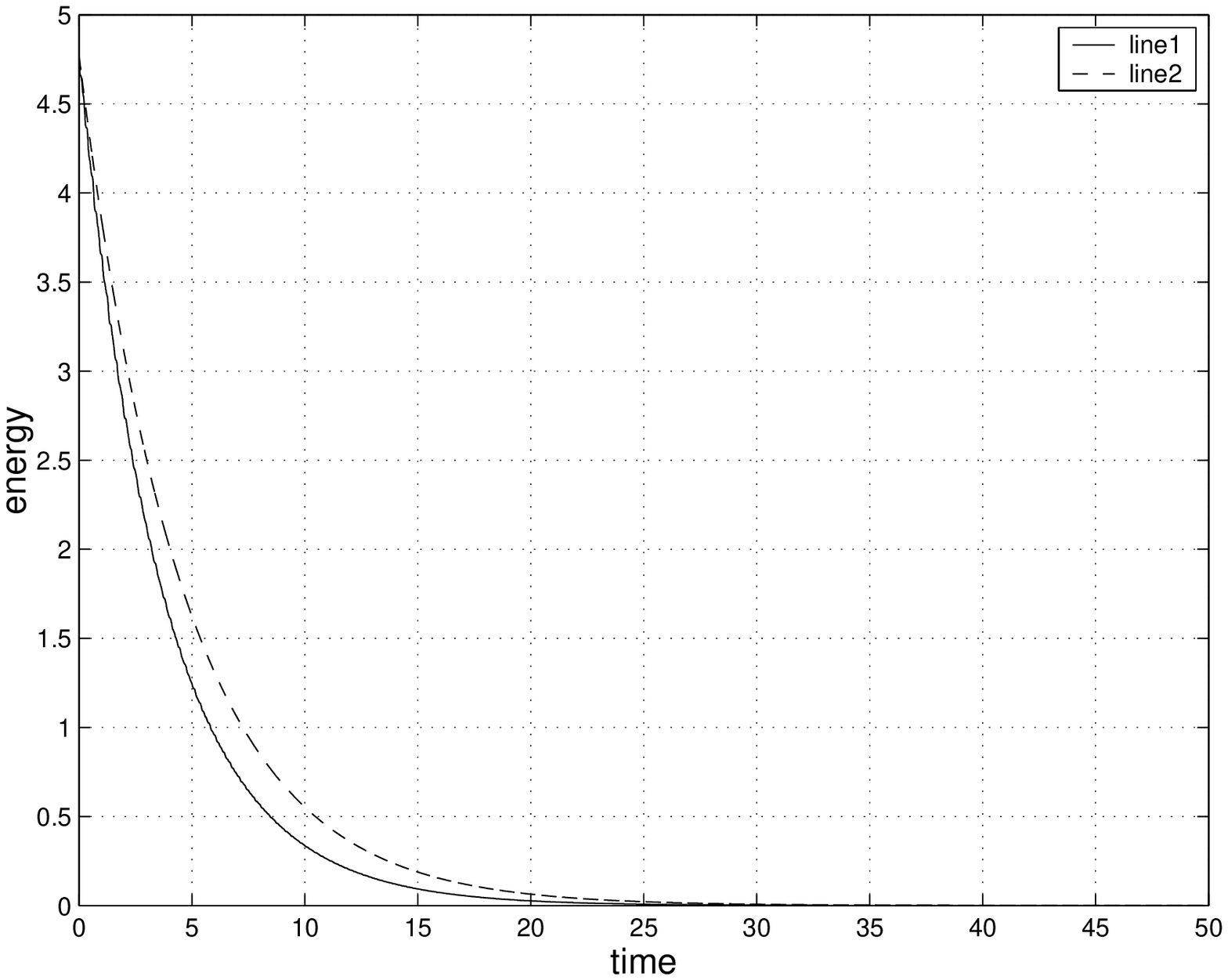}
\includegraphics[angle=0,width=6.5cm]{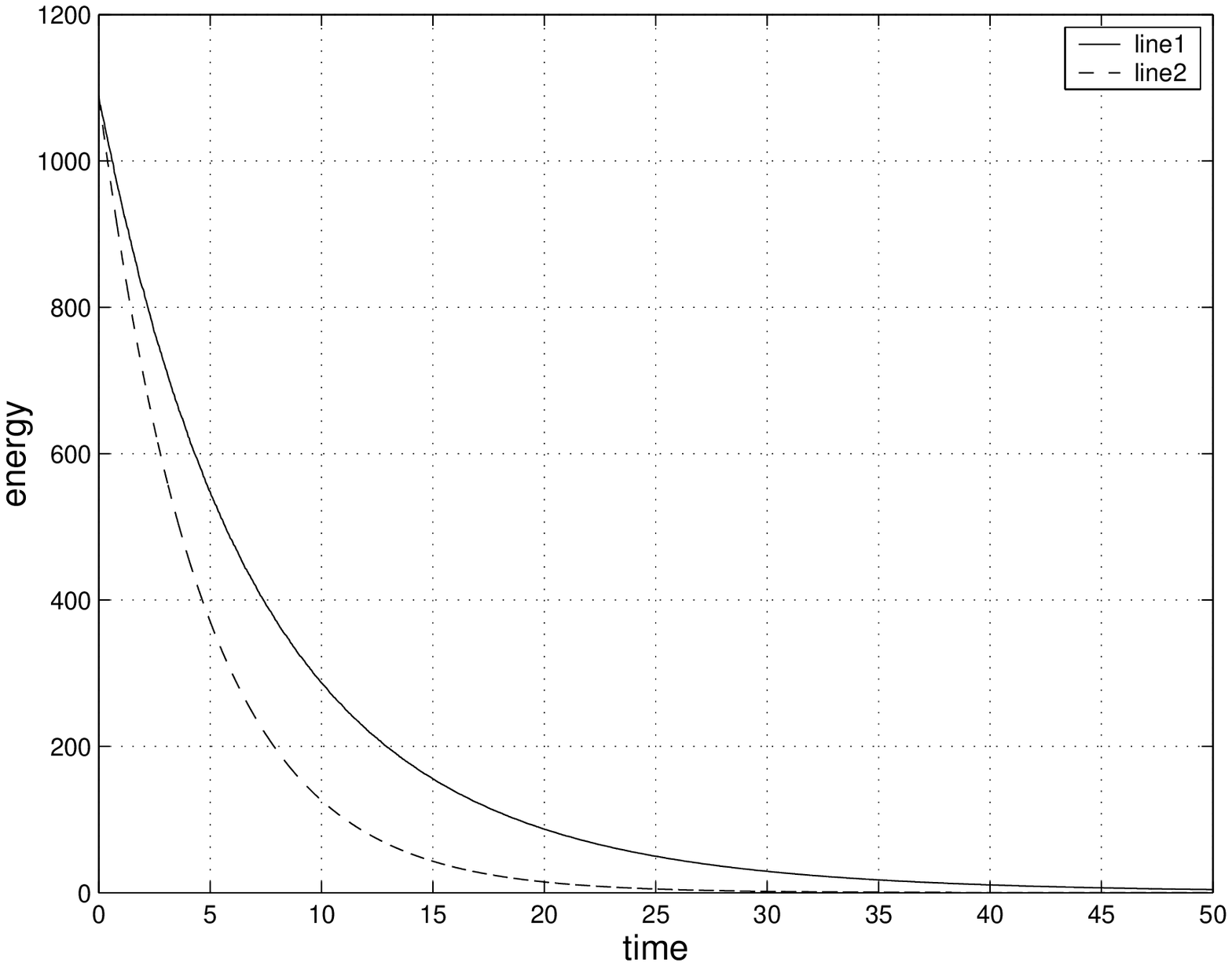}
\end{center}
\begin{center}
figure 3 : energy of the damped plate, $a(x)=1, \omega = (0 , \frac{1}{2}) \times (0, \frac{1}{2})$, 
$n=m=3$ (left), and $n=m=8$ (right).
\end{center}
In figures 4. 5. and 6. we consider a damping in the domain $\omega = (0 , \frac{2}{5}) \times (0, \frac{2}{5})$ and $\omega = (0 , \frac{3}{5}) \times (0, \frac{3}{5})$ :
\begin{center}
\includegraphics[angle=0,width=6.5cm]{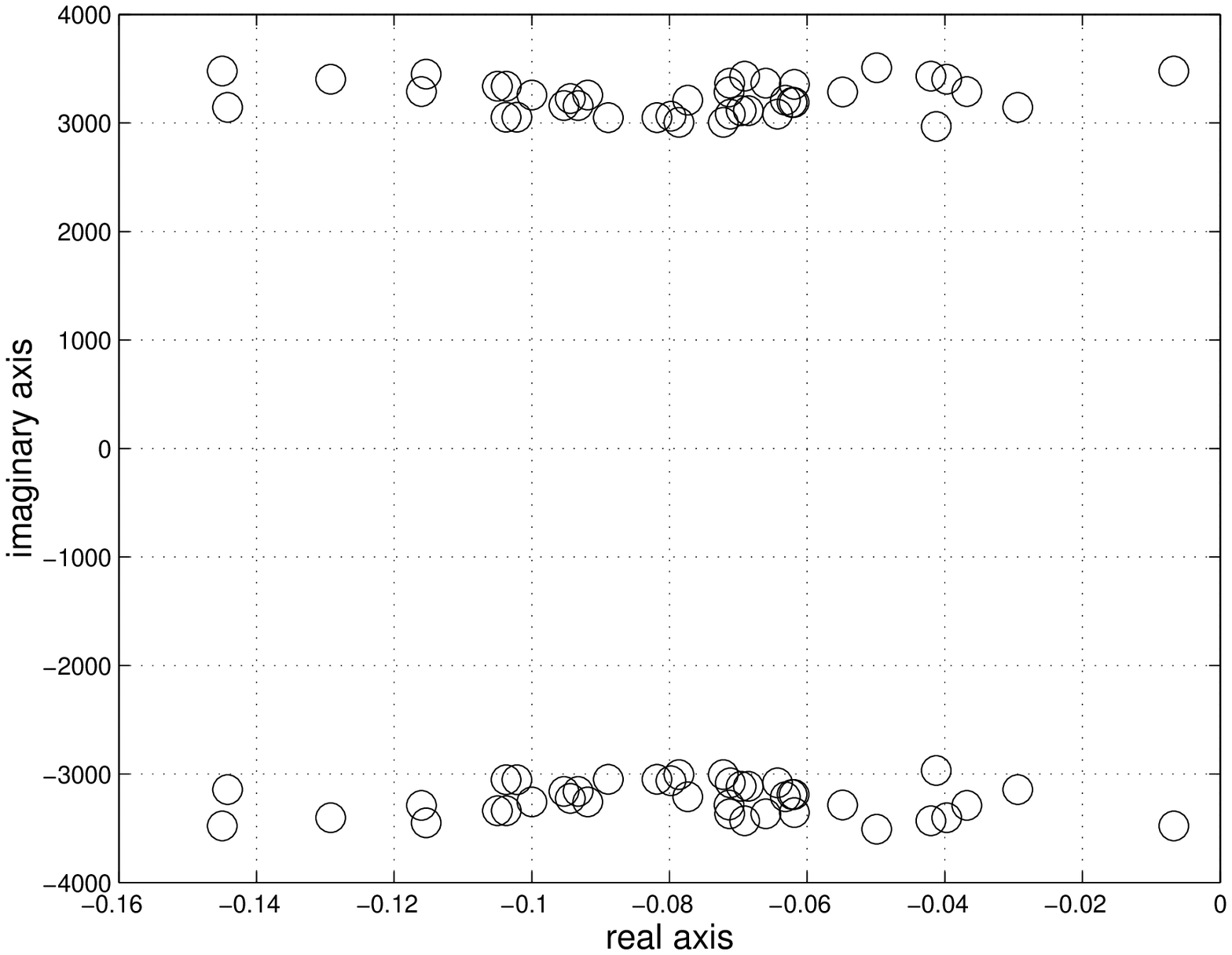}
\includegraphics[angle=0,width=6.5cm]{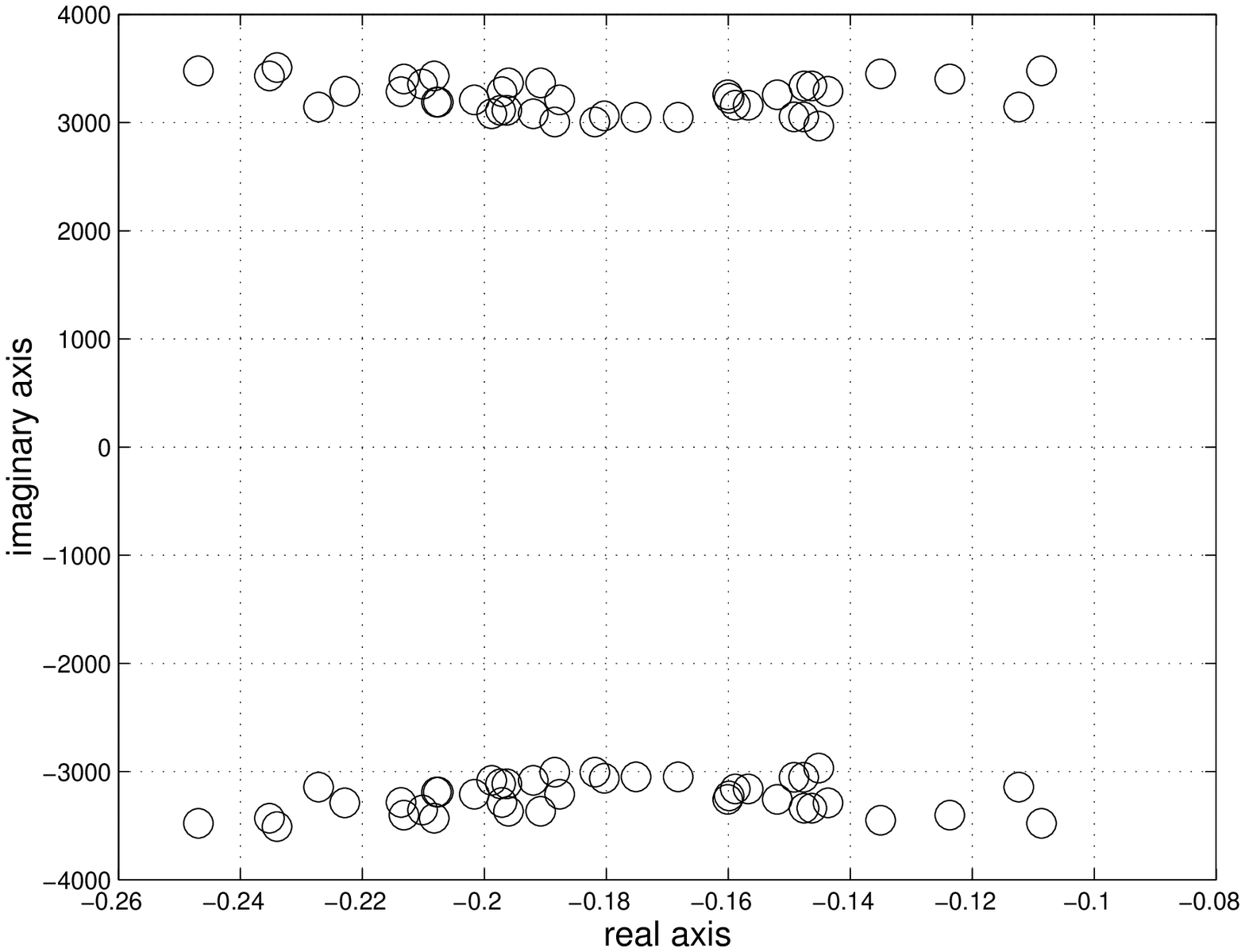}
\end{center}
\begin{center}
figure 4 : spectrum of the damped plate, $a(x)=1$ ,$\omega  =(0,0.4) \times (0,0.4)$ (left), $\omega = (0 , \frac{3}{5}) \times (0, \frac{3}{5})$ (right).
\end{center}
\begin{center}
\includegraphics[angle=0,width=6.5cm]{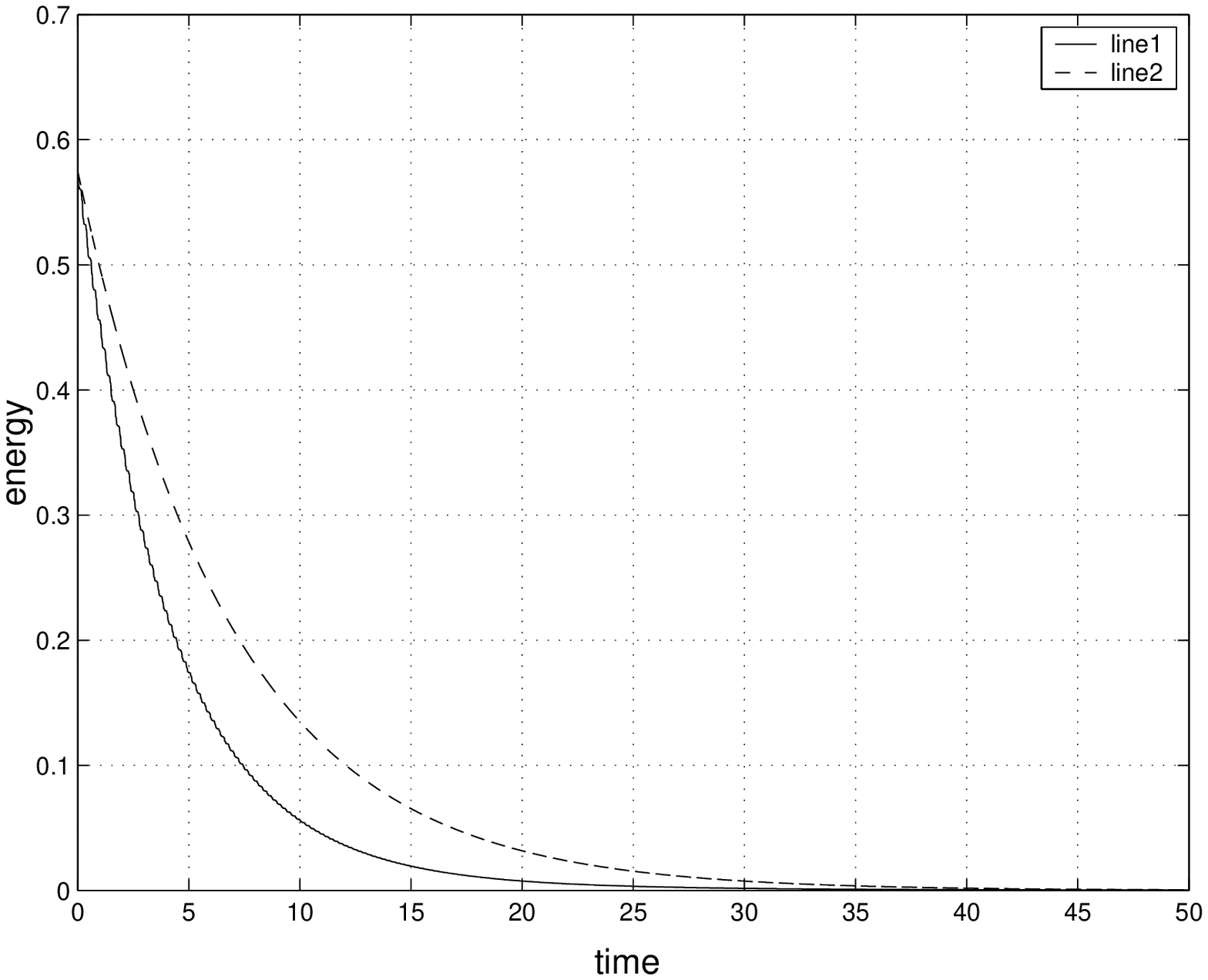}
\includegraphics[angle=0,width=6.5cm]{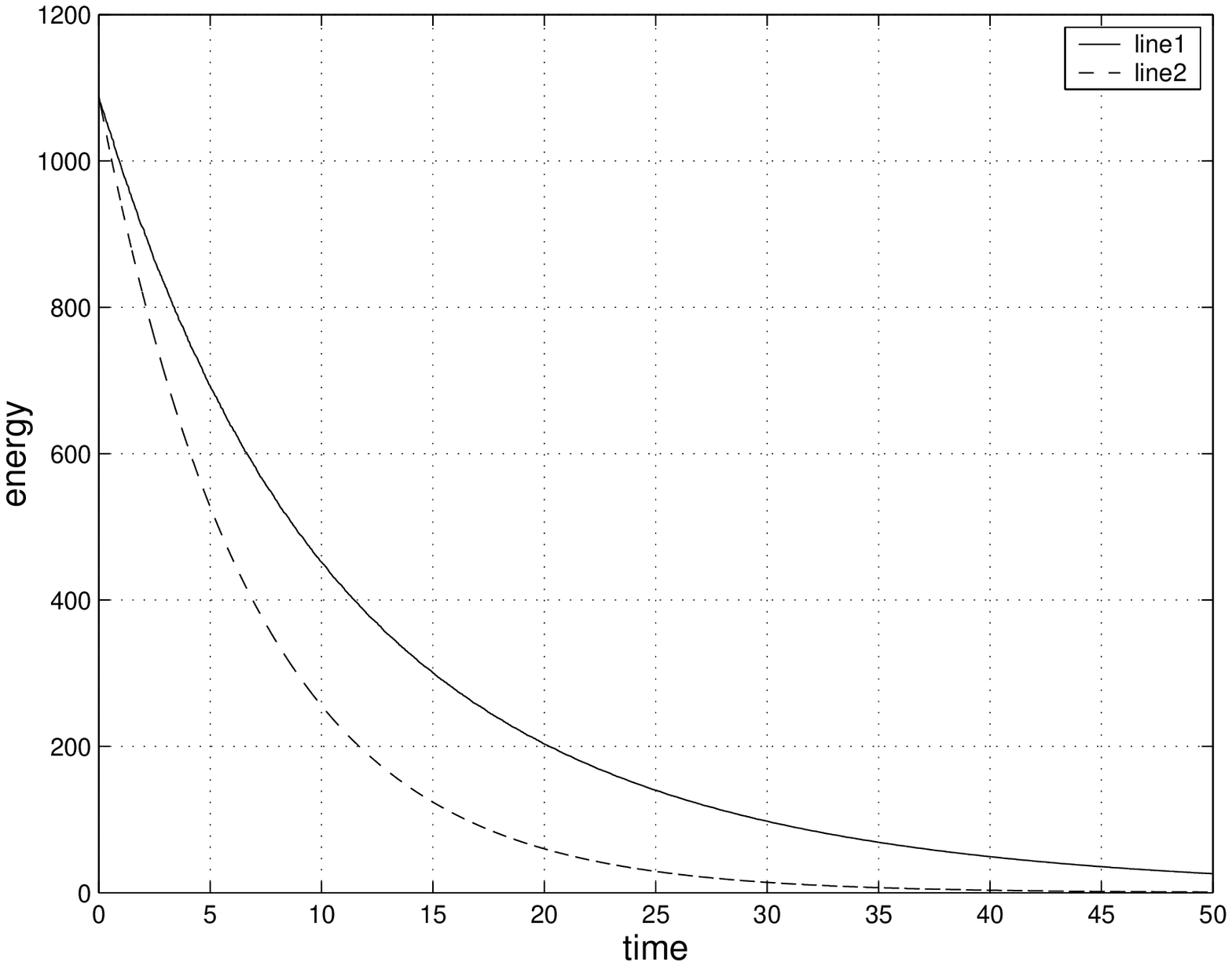}
\end{center}
\begin{center}
figure 5 : energy of the damped plate, $a(x)=1,\omega  =(0,0.4) \times (0,0.4), n=m=2$ (left), and $n=m=8$  (right).
\end{center}
%
\begin{center}
\includegraphics[angle=0,width=6.5cm]{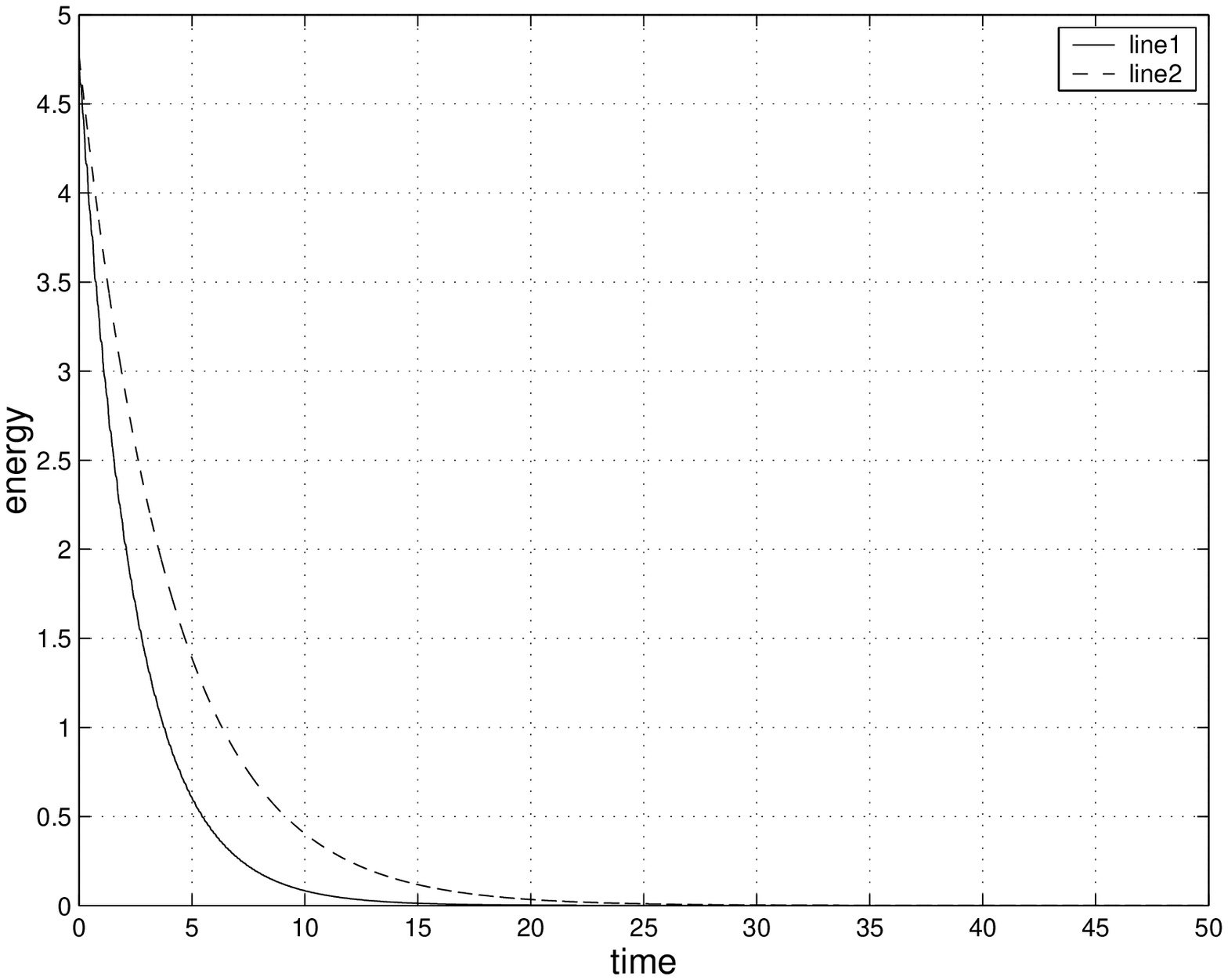}
\includegraphics[angle=0,width=6.5cm]{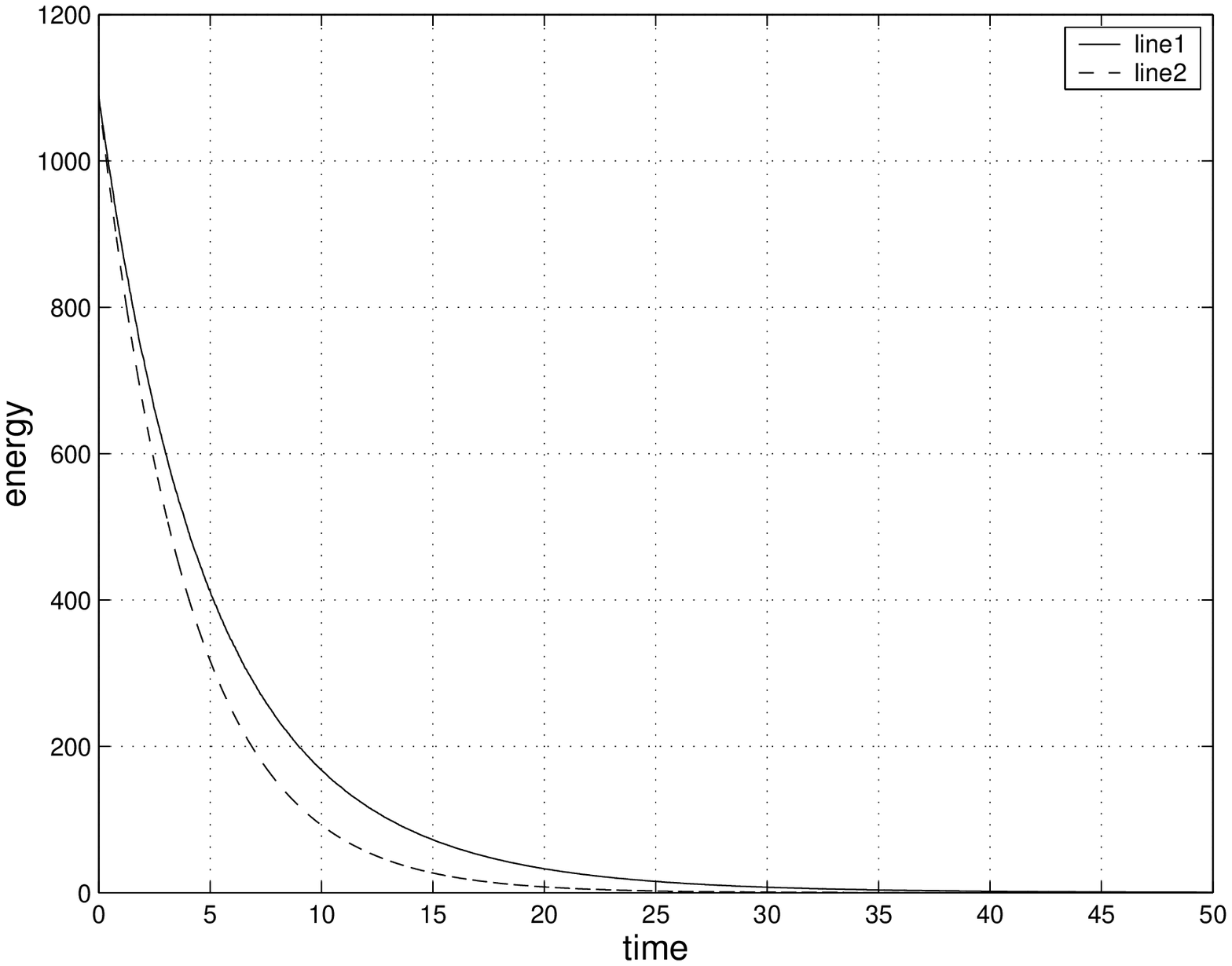}
\end{center}
\begin{center}
figure 6 : energy of the damped plate, $a(x)=1,\omega  =(0,0.6) \times (0,0.6), n=m=3$ (left), and $n=m=8$ (right).
\end{center}
In figures 7., 8. and 9. we consider a damping in the domain $\omega = (0 , \frac{1}{4}) \times (0, \frac{1}{4})$ and $\omega = (0 , \frac{1}{5}) \times (0, \frac{1}{5})$ :
\begin{center}
\includegraphics[angle=0,width=6.5cm]{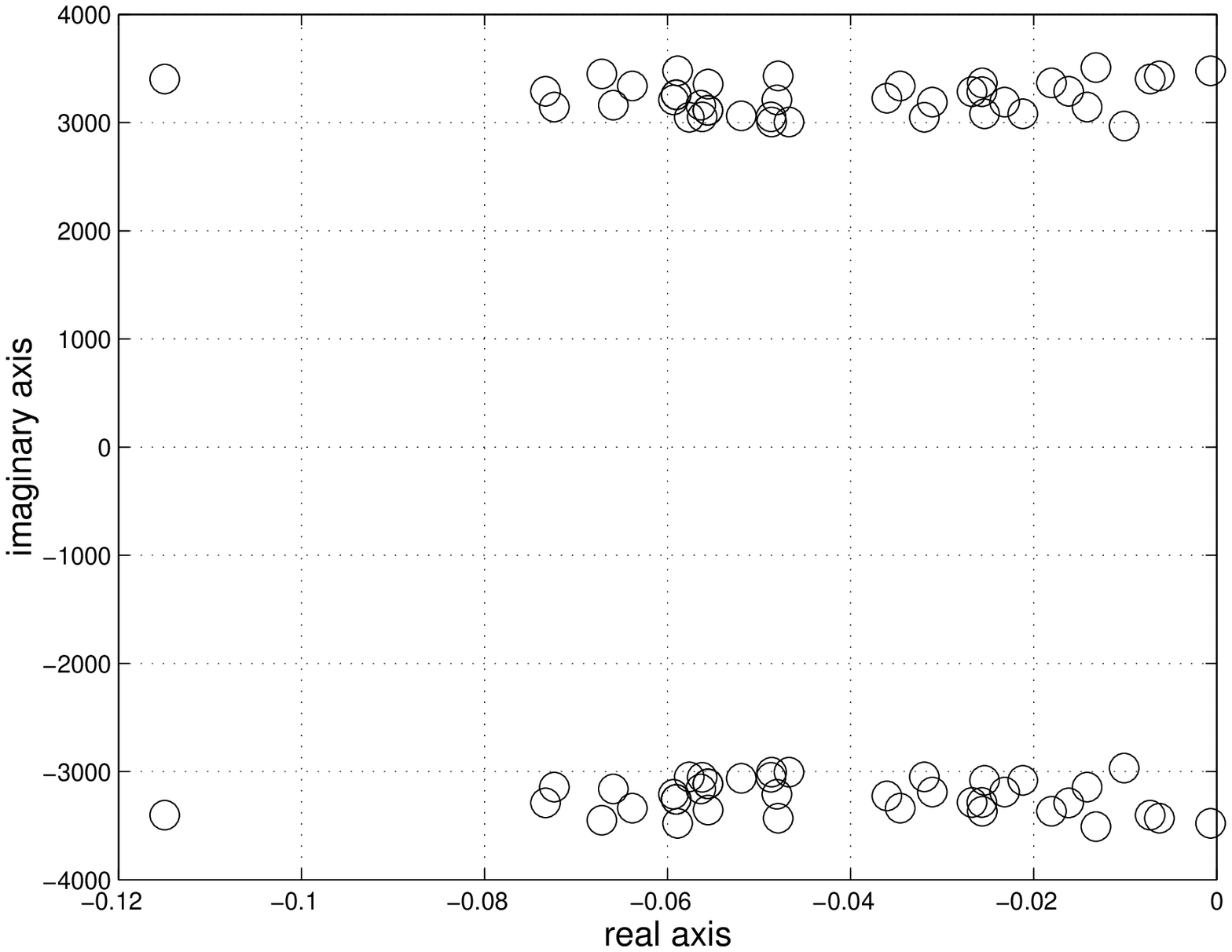}
\includegraphics[angle=0,width=6.5cm]{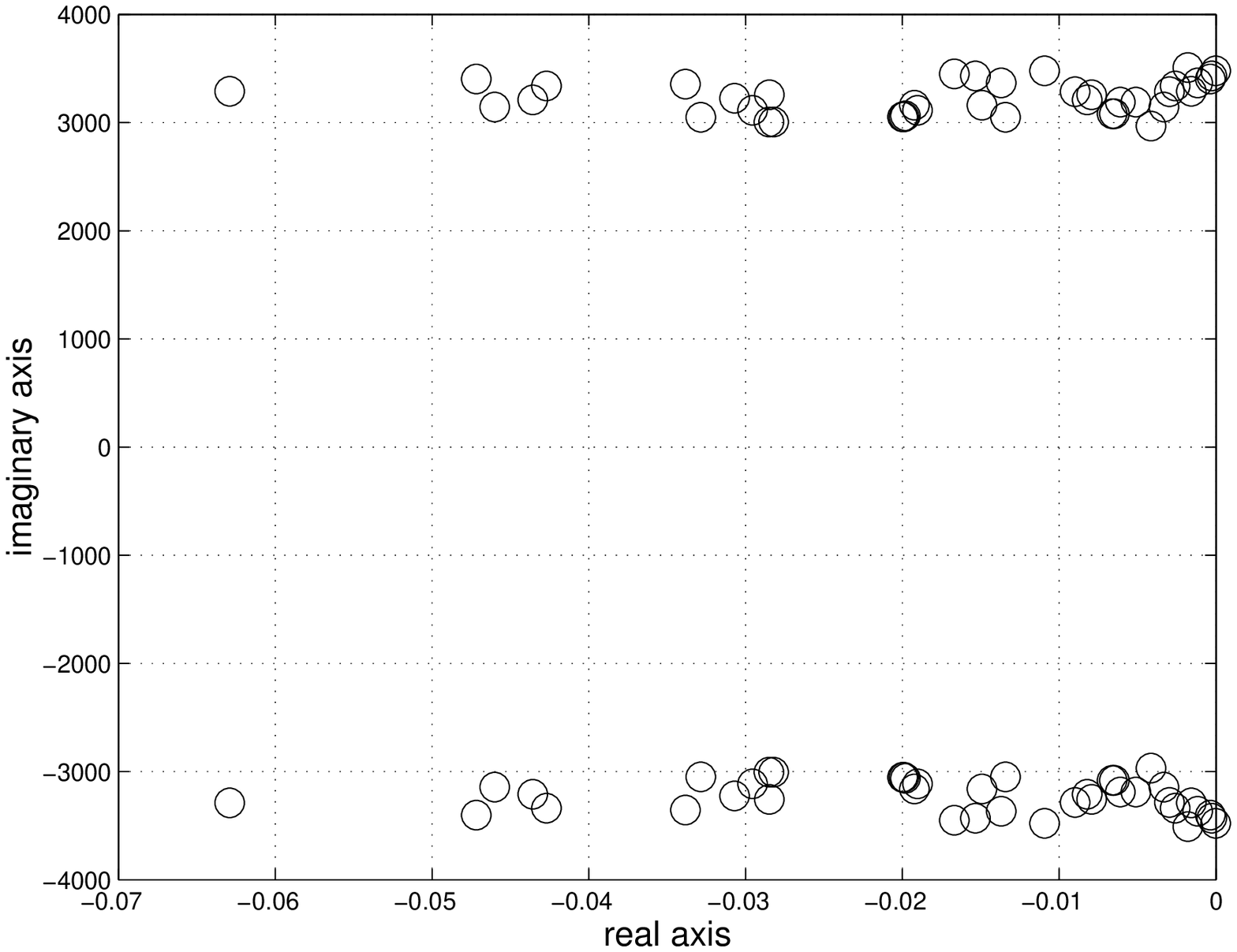}
\end{center}
\begin{center}
figure 7 : spectrum of the damped plate, $a(x)=1$, $\omega = (0 , \frac{1}{4}) \times (0, \frac{1}{4})$ (left) and  $\omega = (0 , \frac{1}{5}) \times (0, \frac{1}{5})$ (right).
\end{center}
\begin{center}
\includegraphics[angle=0,width=6.5cm]{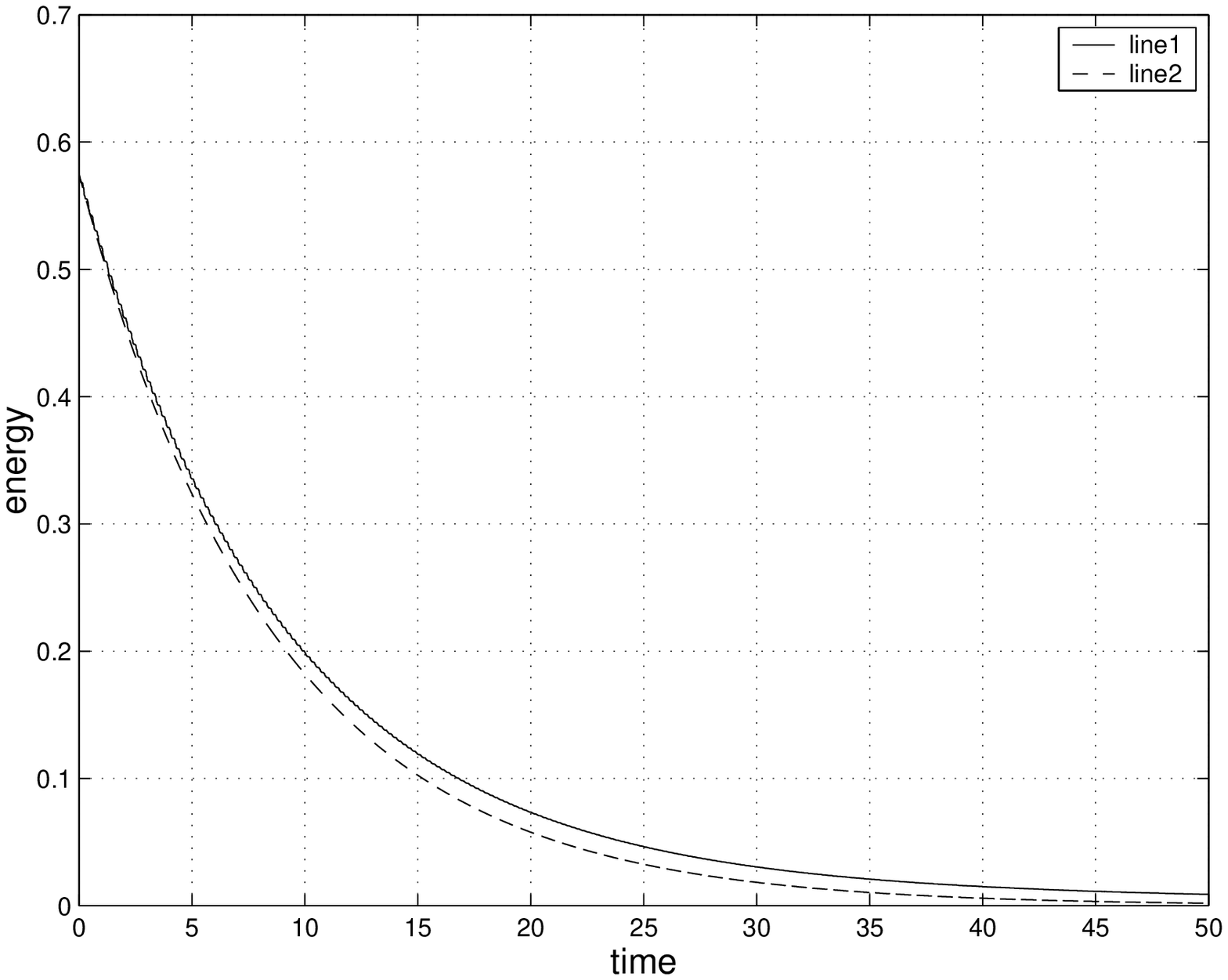}
\includegraphics[angle=0,width=6.5cm]{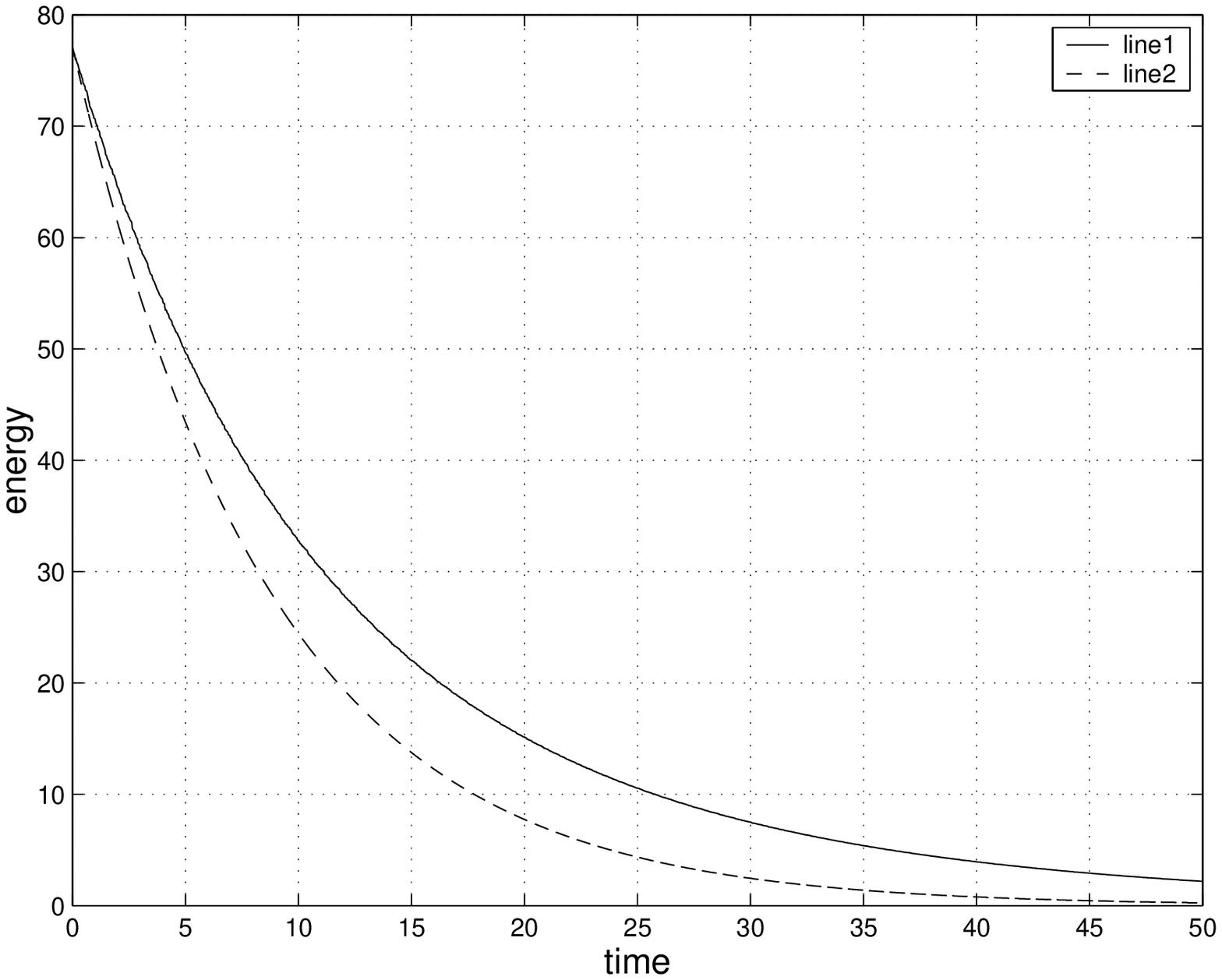}
\end{center}
\begin{center}
figure 8 : energy of the damped plate, $a(x)=1,\omega  =(0,0.25) \times (0,0.25),  n=m=2$ (left), and $n=m=5$ (right).
\end{center}
%
\begin{center}
\includegraphics[angle=0,width=6.5cm]{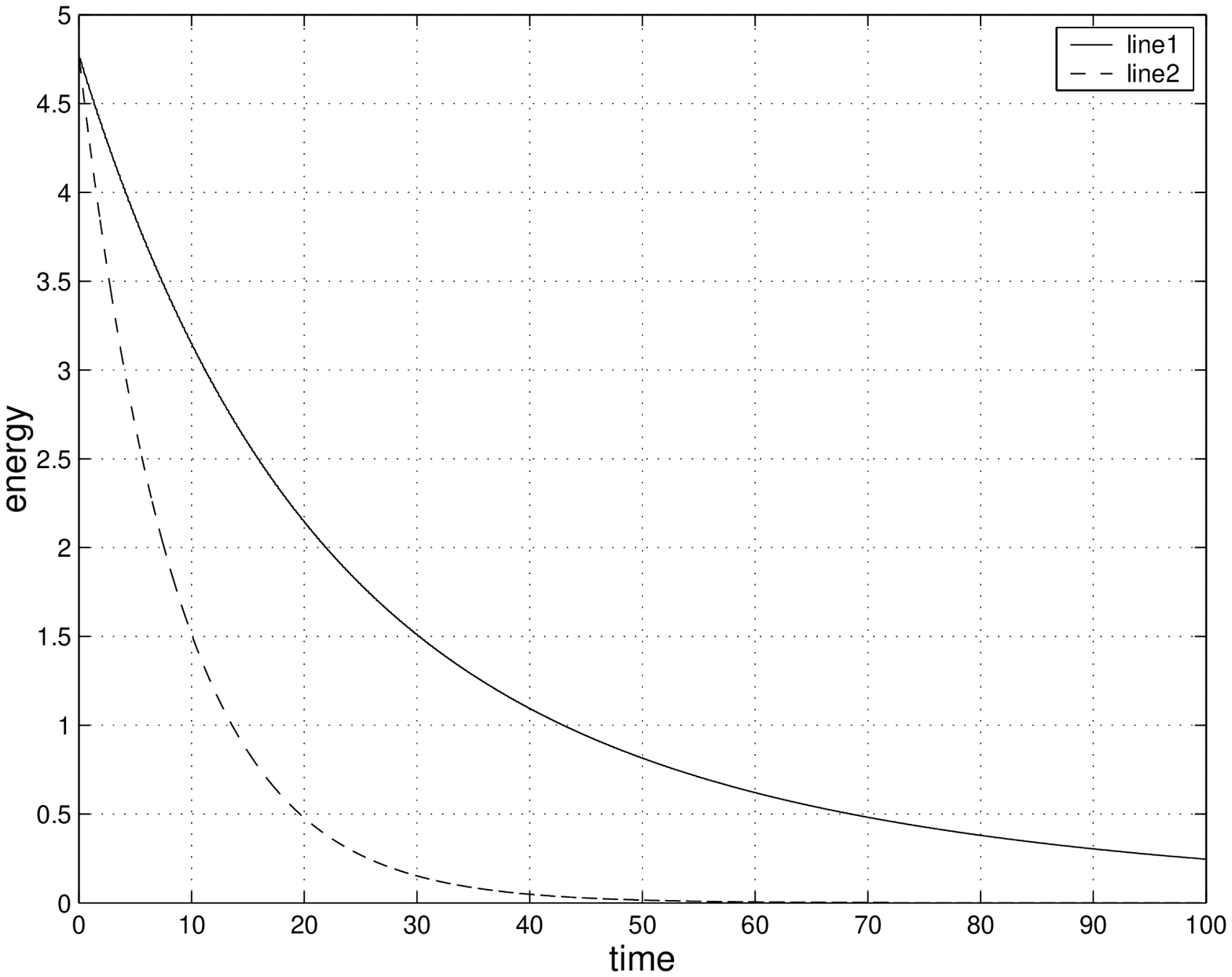}
\includegraphics[angle=0,width=6.5cm]{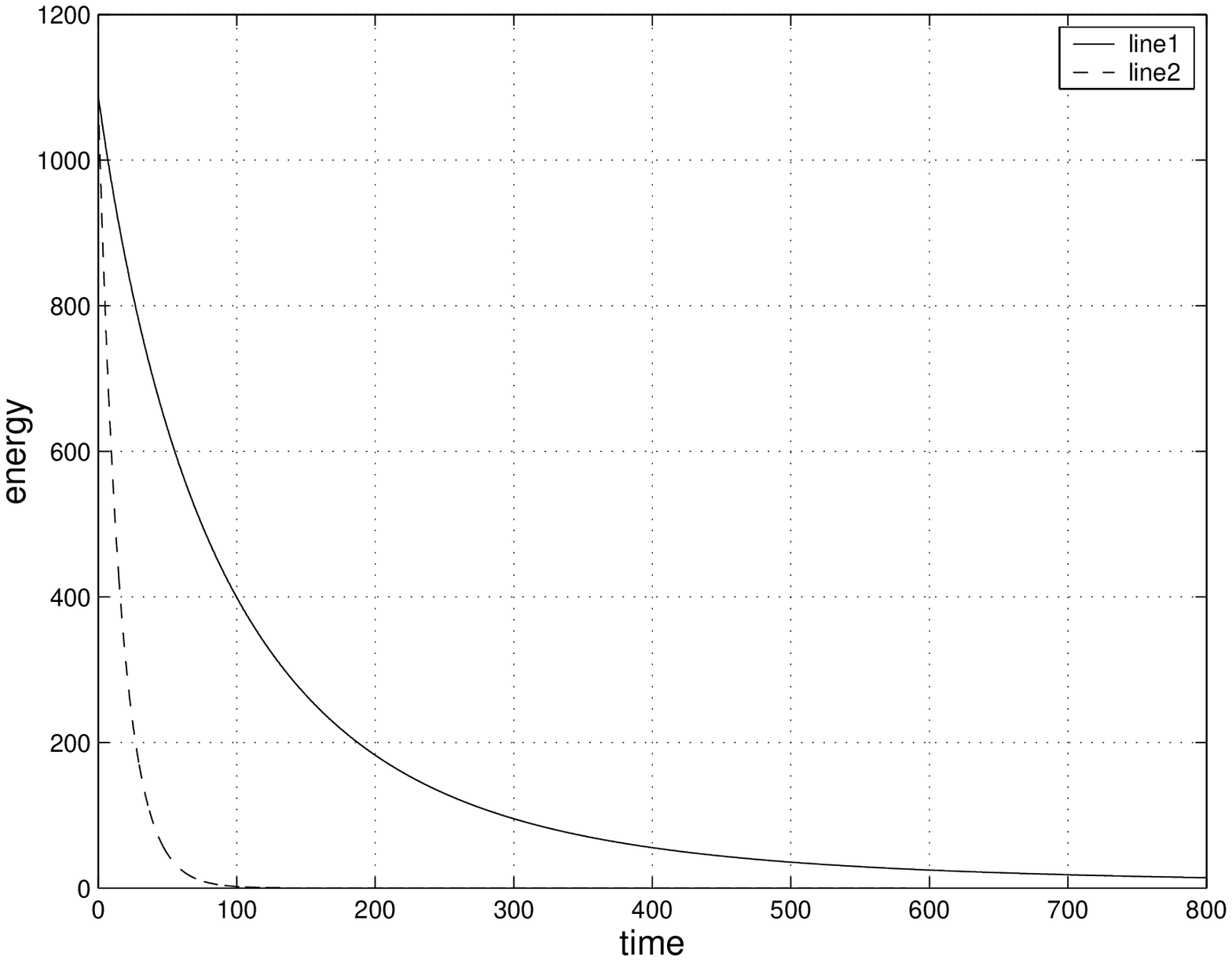}
\end{center}
\begin{center}
figure 9 : energy of the damped plate, $a(x)=1,\omega  =(0,0.2) \times (0,0.2), n=m=3$ (left), and $n=m=8$ (right).
\end{center}
In figures 10., 11. and 12. we consider a damping in the domain $\omega = (0 , 1) \times (0, \frac{1}{2})$ and  $\omega = (0 , 1) \times (0, \frac{1}{4})$ :
\begin{center}
\includegraphics[angle=0,width=6.5cm]{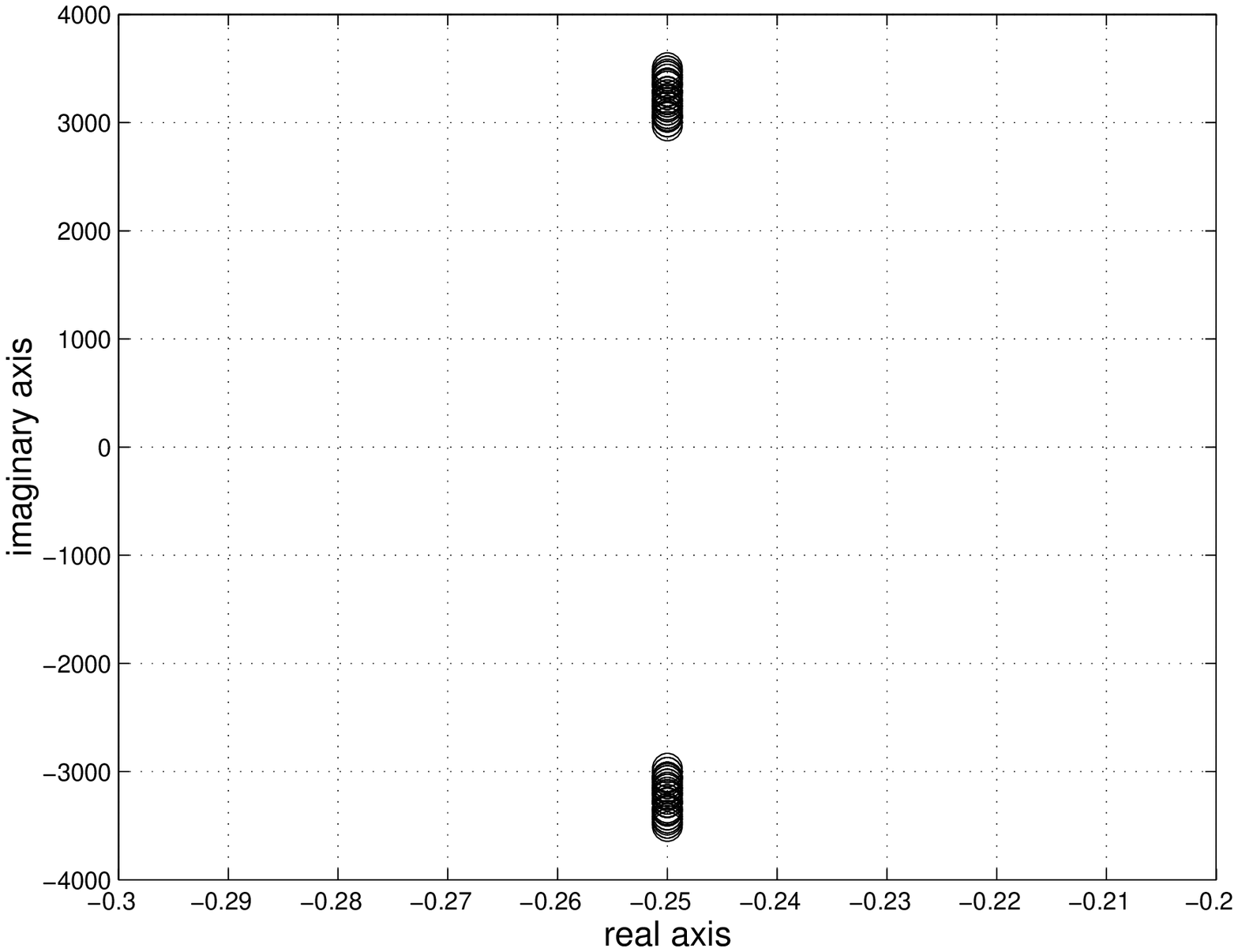}
\includegraphics[angle=0,width=6.5cm]{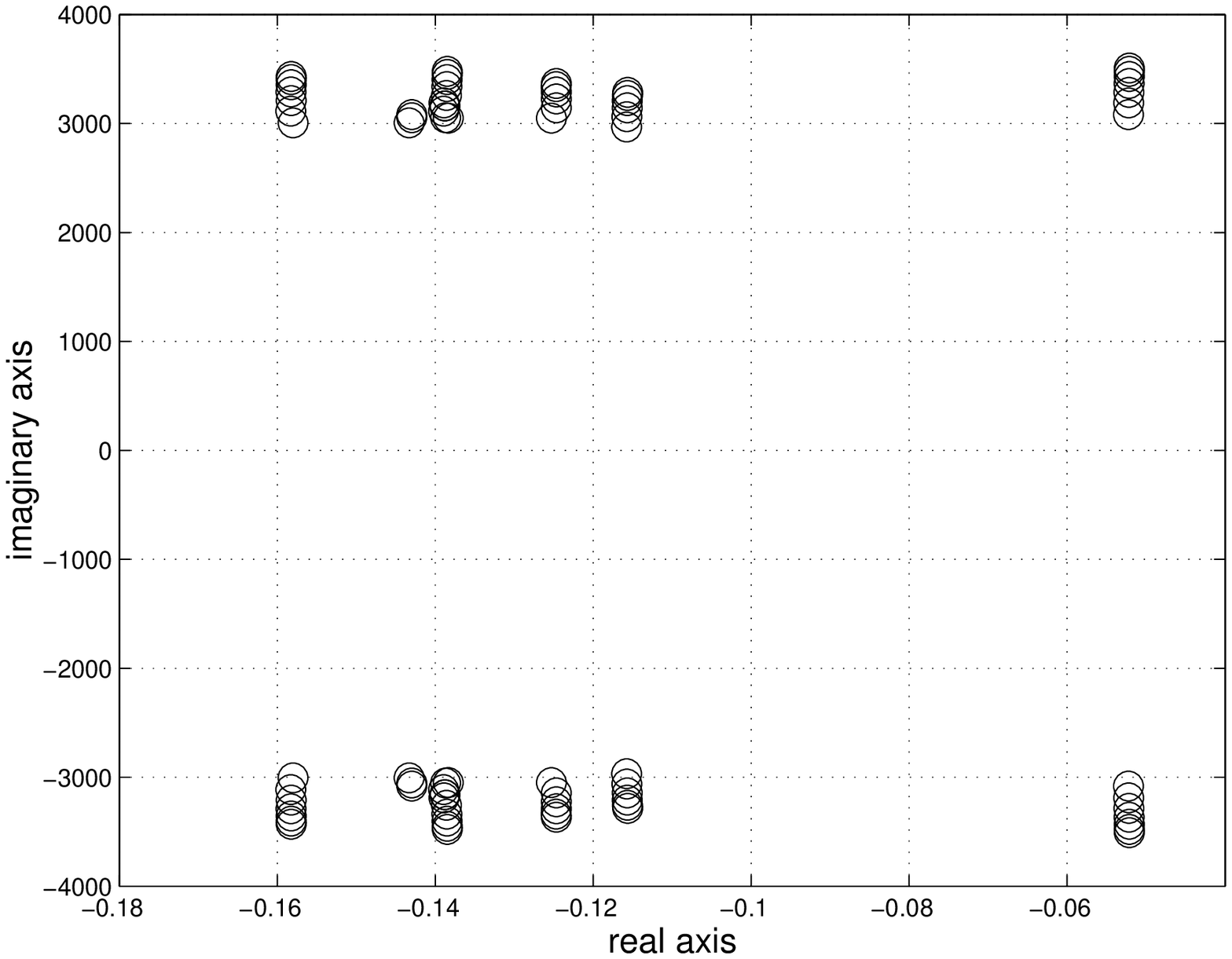}
\end{center}
\begin{center}
figure 10 : spectrum of the damped plate, $a(x)=1$ and $\omega = (0 , 1) \times (0, \frac{1}{2})$ (left) , $a(x)=1$ and $\omega = (0 , 1) \times (0, \frac{1}{4})$ (right).
\end{center}
\begin{center}
\includegraphics[angle=0,width=6.5cm]{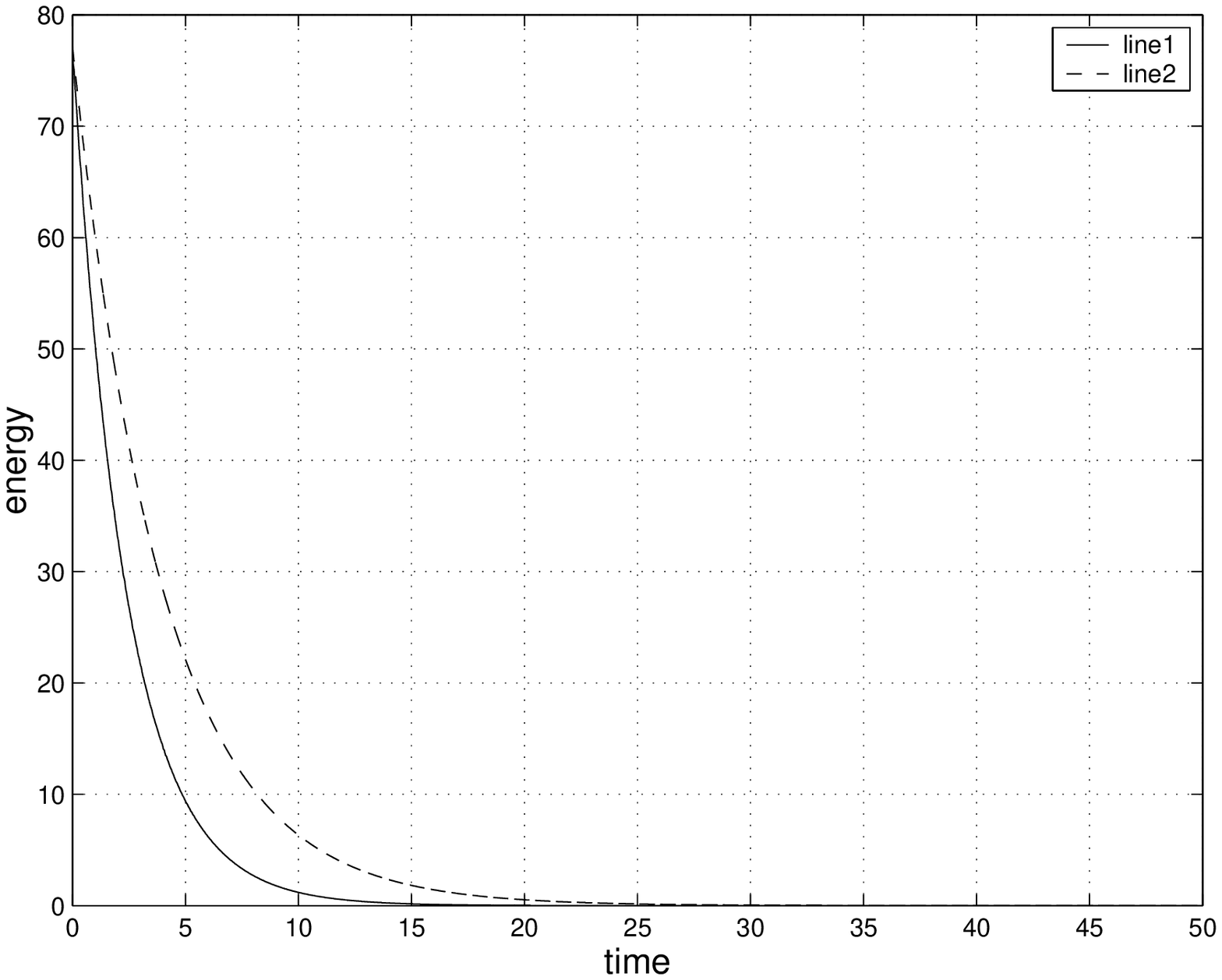}
\includegraphics[angle=0,width=6.5cm]{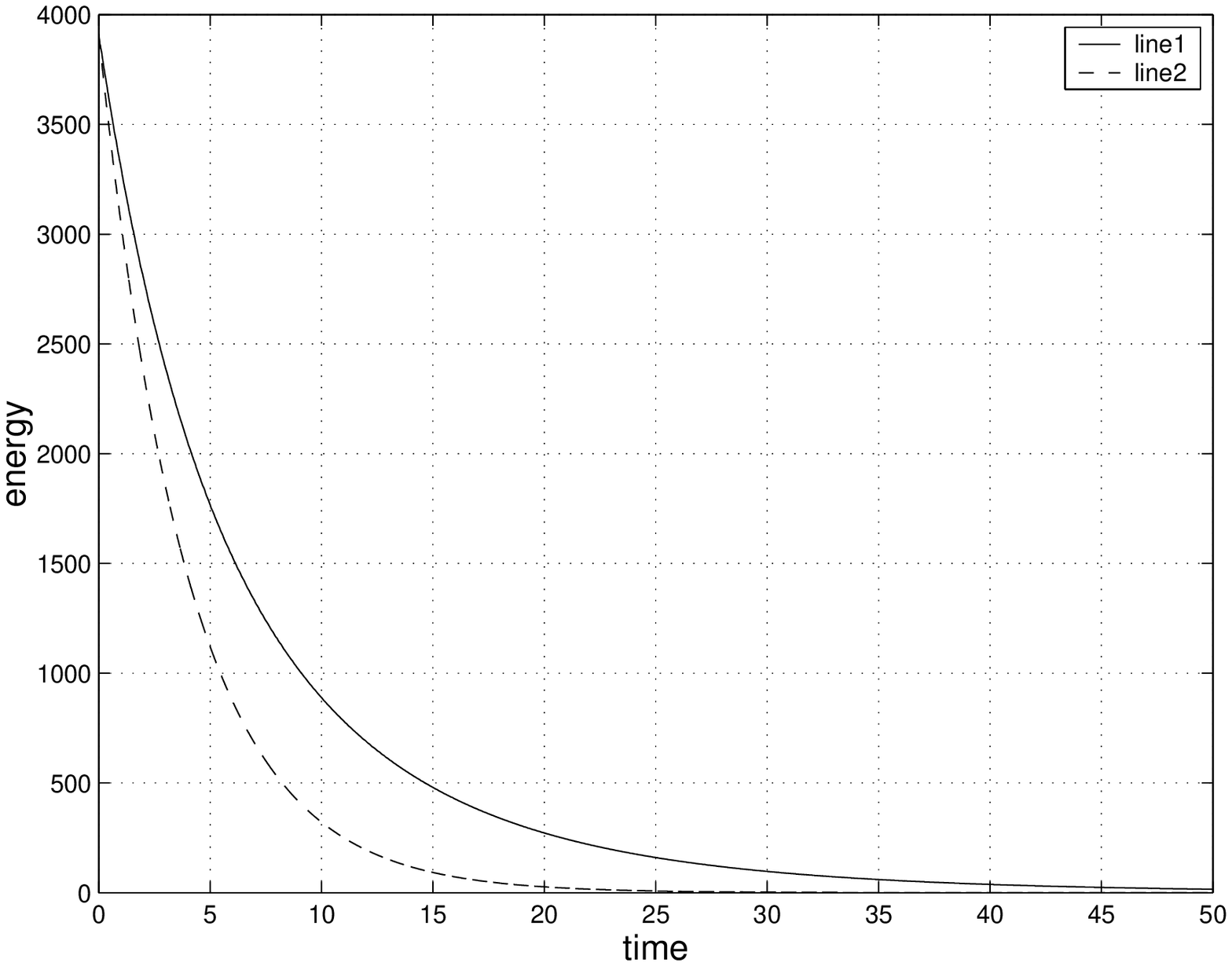}
\end{center}
\begin{center}
figure 11 : energy of the damped plate, $a(x)=1,\omega  =(0,1) \times (0,0.5),  n=m=5$ (left), and $n=m=10$ (right).
\end{center}
%
\begin{center}
\includegraphics[angle=0,width=6.5cm]{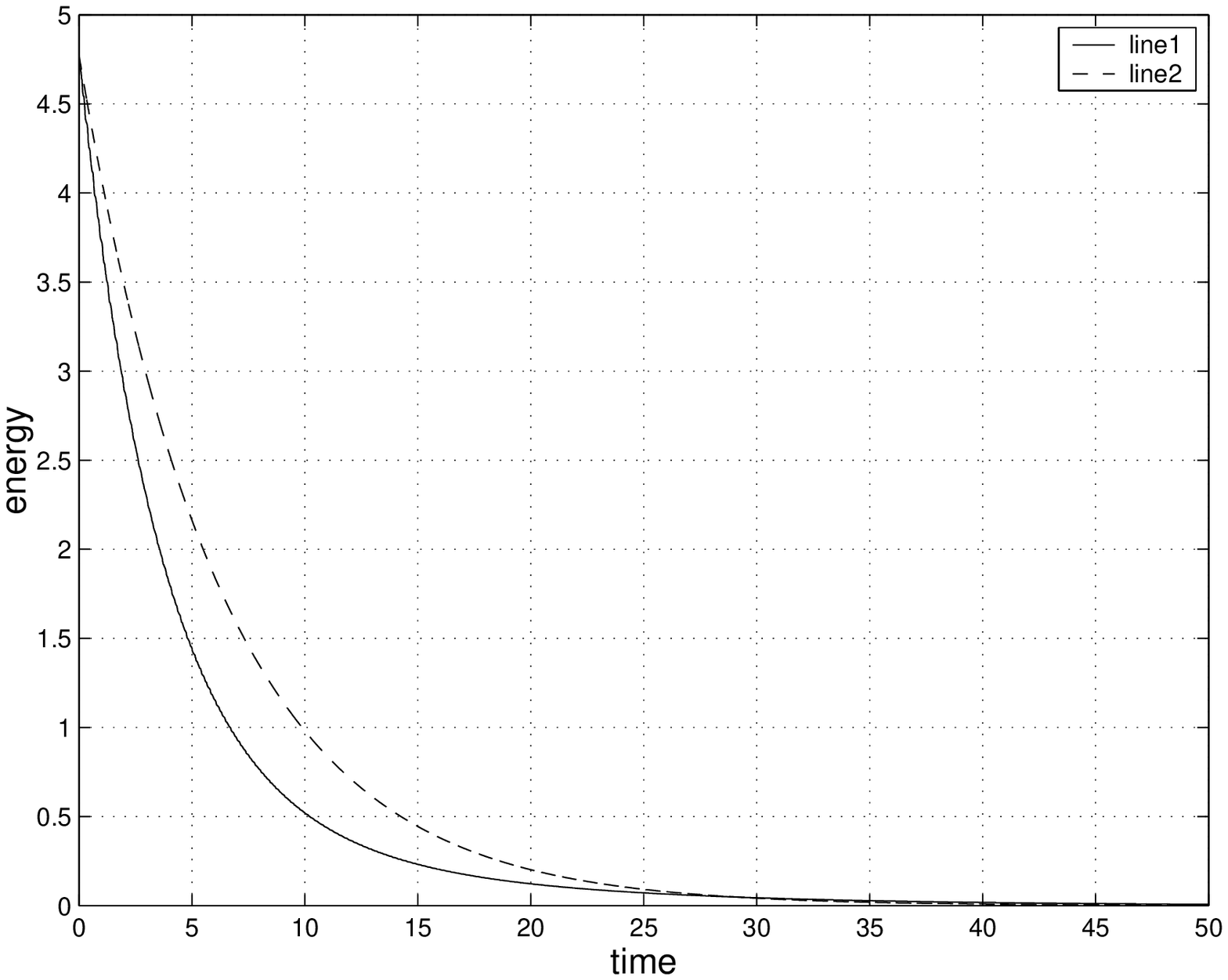}
\includegraphics[angle=0,width=6.5cm]{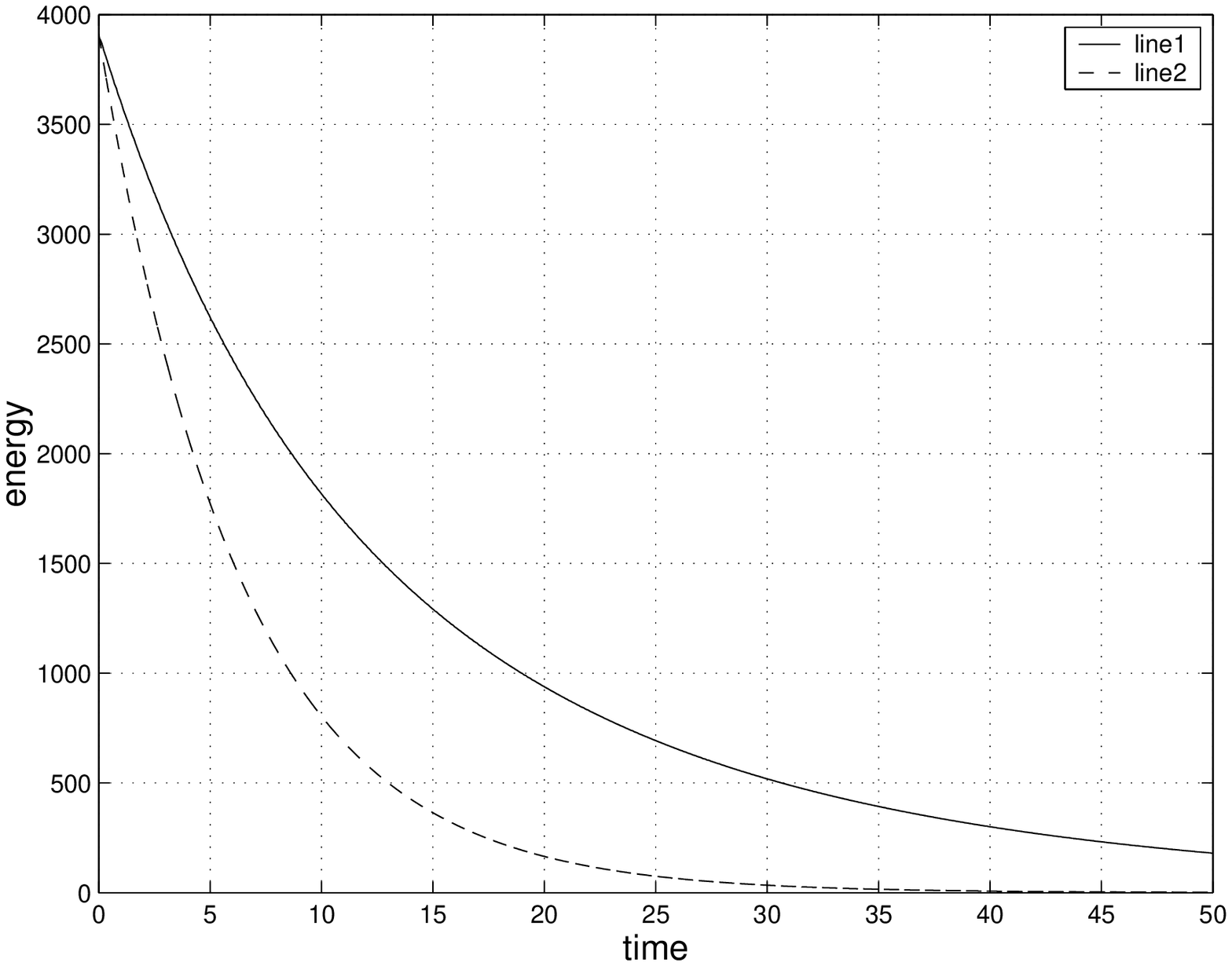}
\end{center}
\begin{center}
figure 12 : energy of the damped plate, $a(x)=1,\omega  =(0,1) \times (0, \frac{1}{4})$, $n=m=3$ (left), and $n=m=10$ (right).
\end{center}
%
\begin{center}
\includegraphics[angle=0,width=6.5cm]{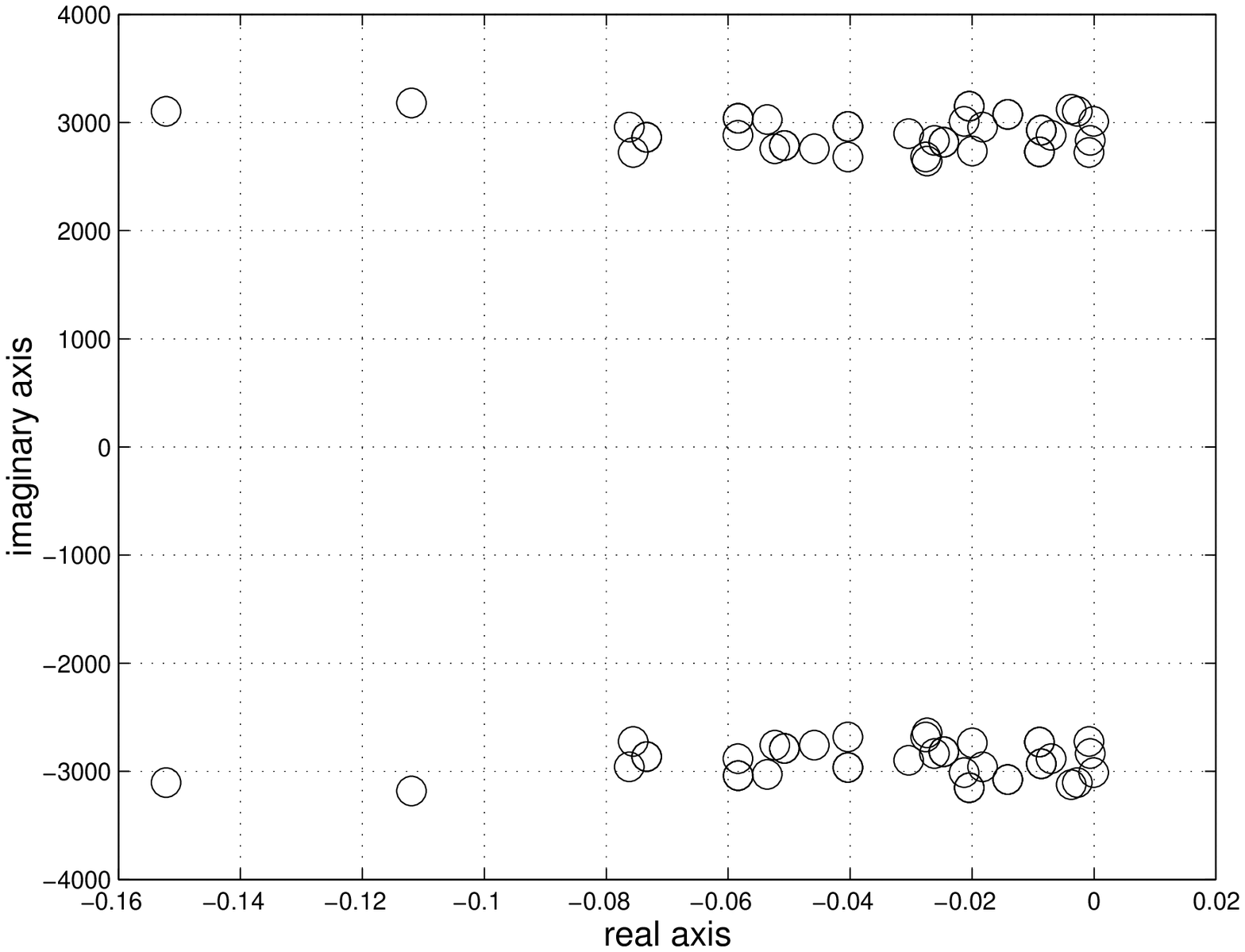}
\includegraphics[angle=0,width=6.5cm]{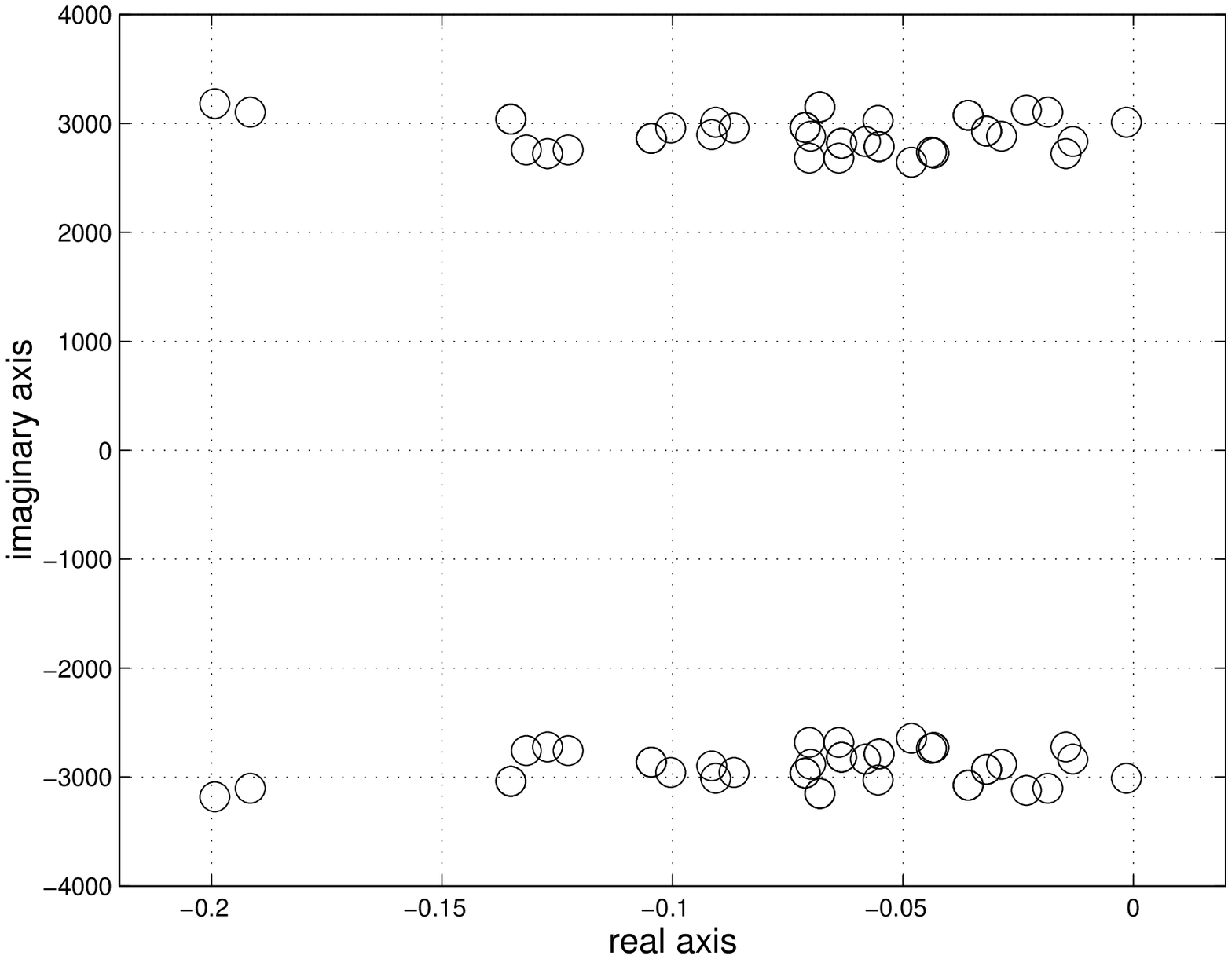}
\end{center}
\begin{center}
figure 13 : spectrum of the damped plate, $a(x)=1,\omega  = ( 0.4,0.6 ) \times  ( 0.4,0.6 )$ (left), $a(x)=1,\omega  = ( 0.35,0.65 ) \times  ( 0.35,0.65 )$ (right).
\end{center}
\begin{center}
\includegraphics[angle=0,width=6.5cm]{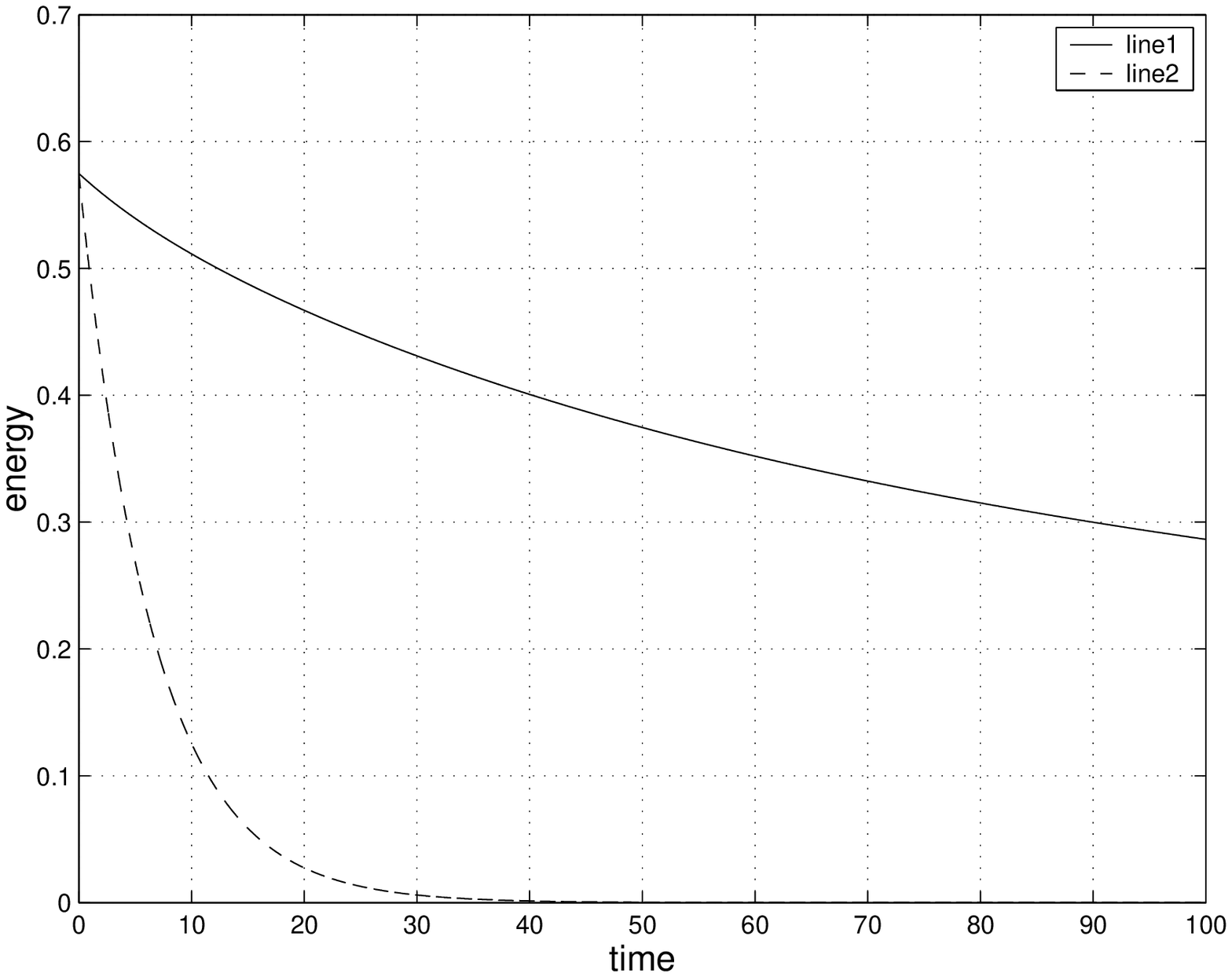}
\includegraphics[angle=0,width=6.5cm]{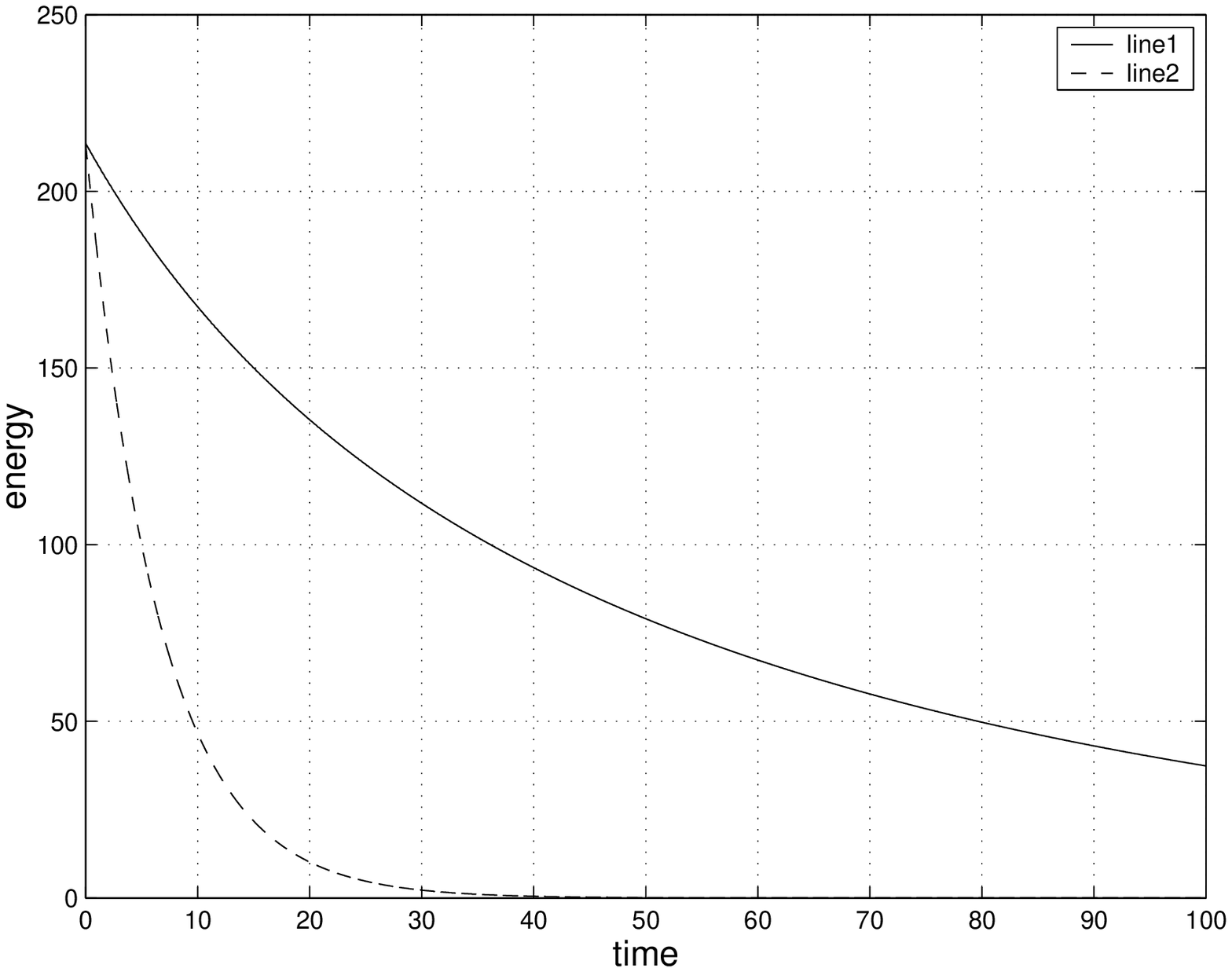}
\end{center}
\begin{center}
figure 14 : energy of the damped plate, $a(x)=1,\omega  = ( 0.4,0.6 ) \times  ( 0.4,0.6 )$, $n=m=2$ (left),  and $n=m=6$ (right).
\end{center}
%
\begin{center}
\includegraphics[angle=0,width=6.5cm]{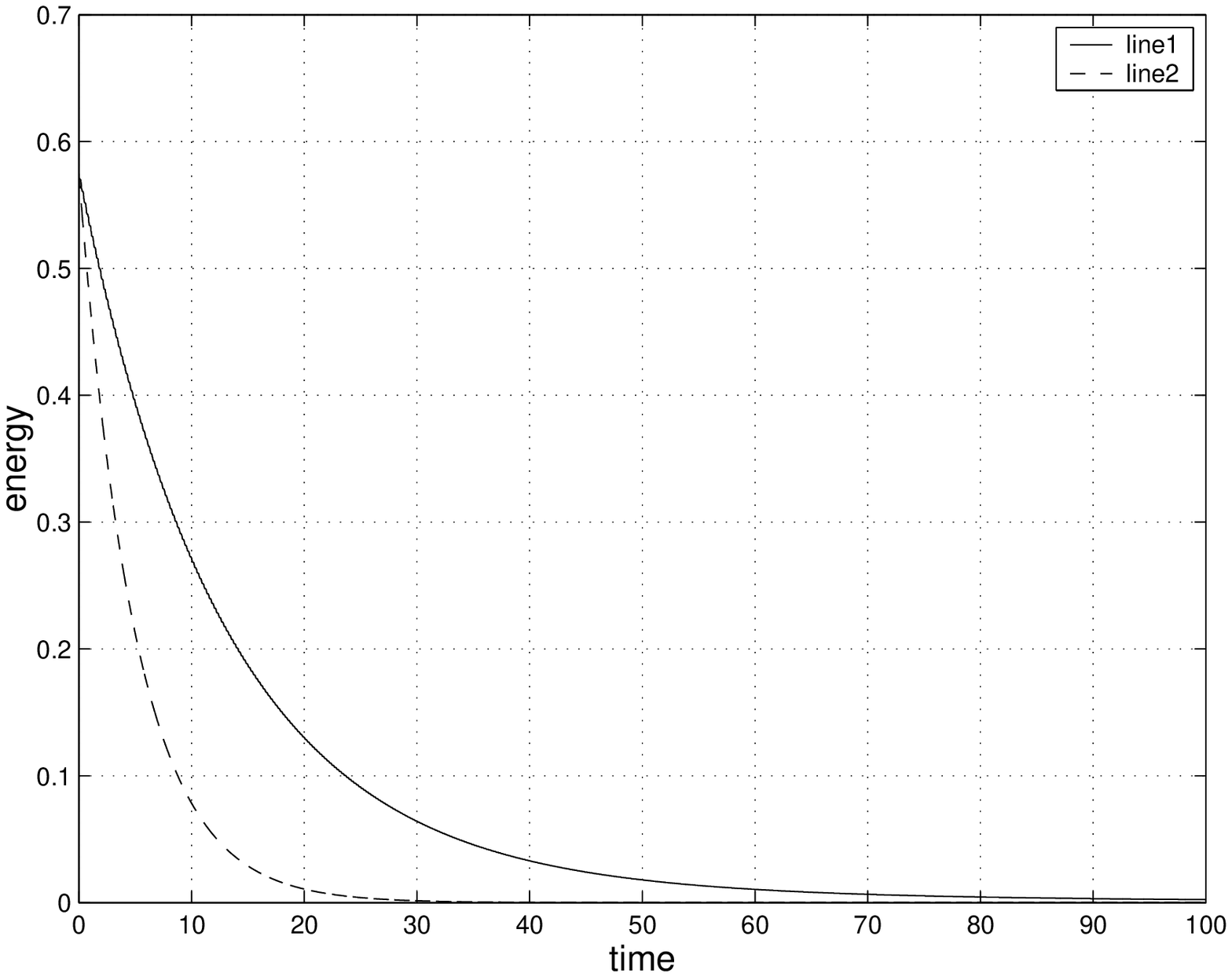}
\includegraphics[angle=0,width=6.5cm]{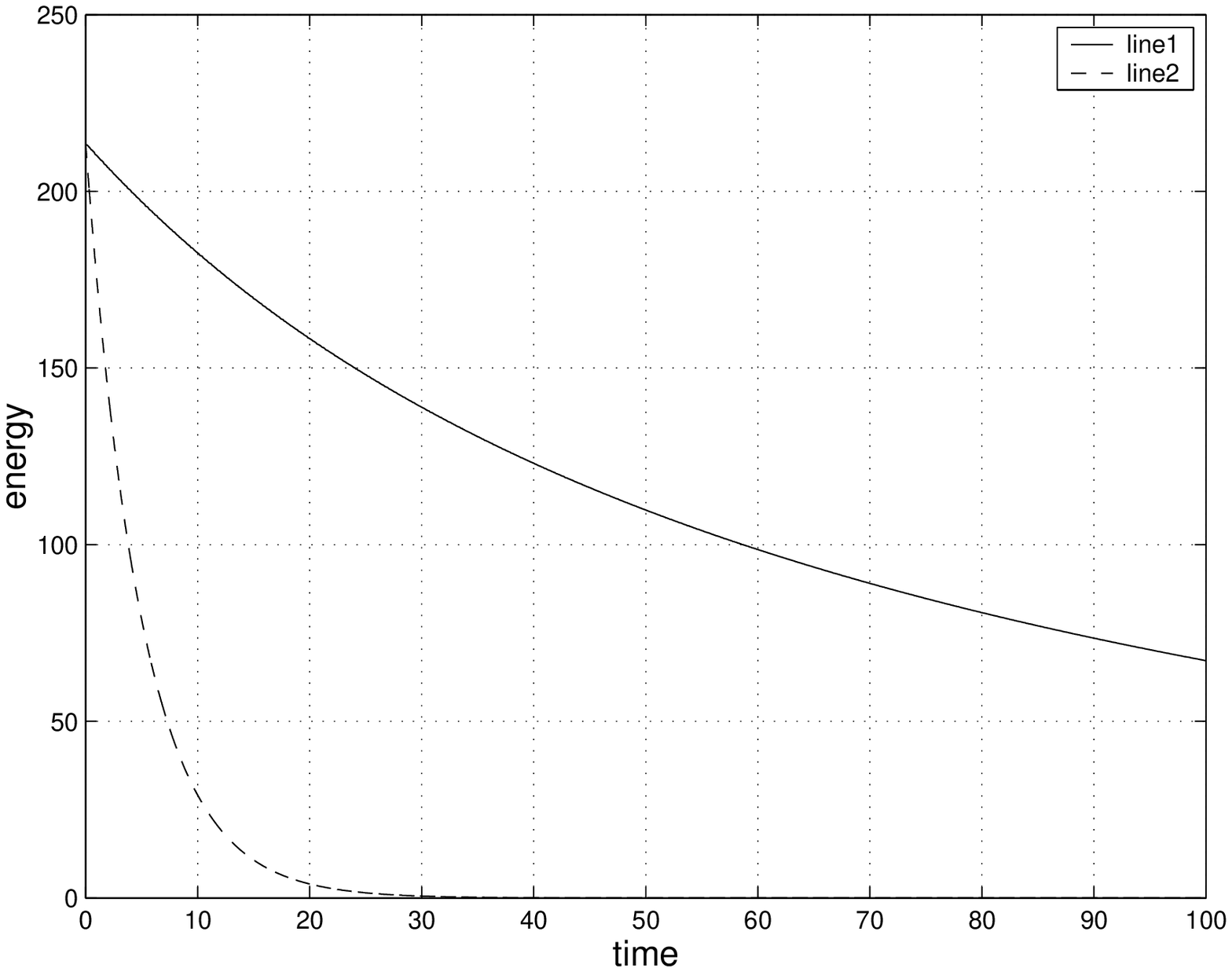}
\end{center}
\begin{center}
figure 15 : energy of the damped plate,  $a(x)=1,\omega  = ( 0.35,0.65 ) \times  ( 0.35,0.65 )$, $n=m=2$ (left),  and $n=m=6$ (right).
\end{center}
%
%
\begin{center}
\includegraphics[angle=0,width=6.5cm]{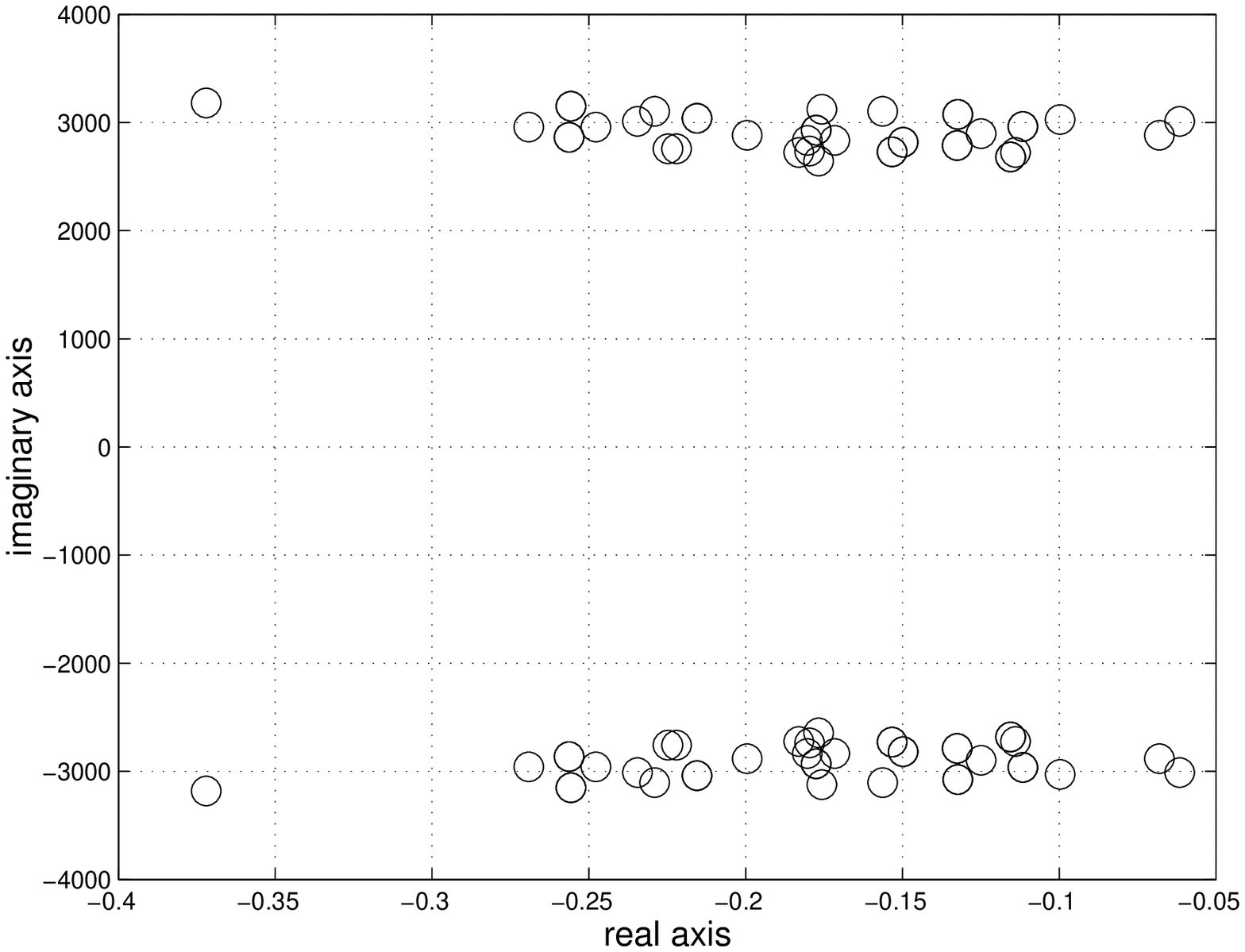}
\includegraphics[angle=0,width=6.5cm]{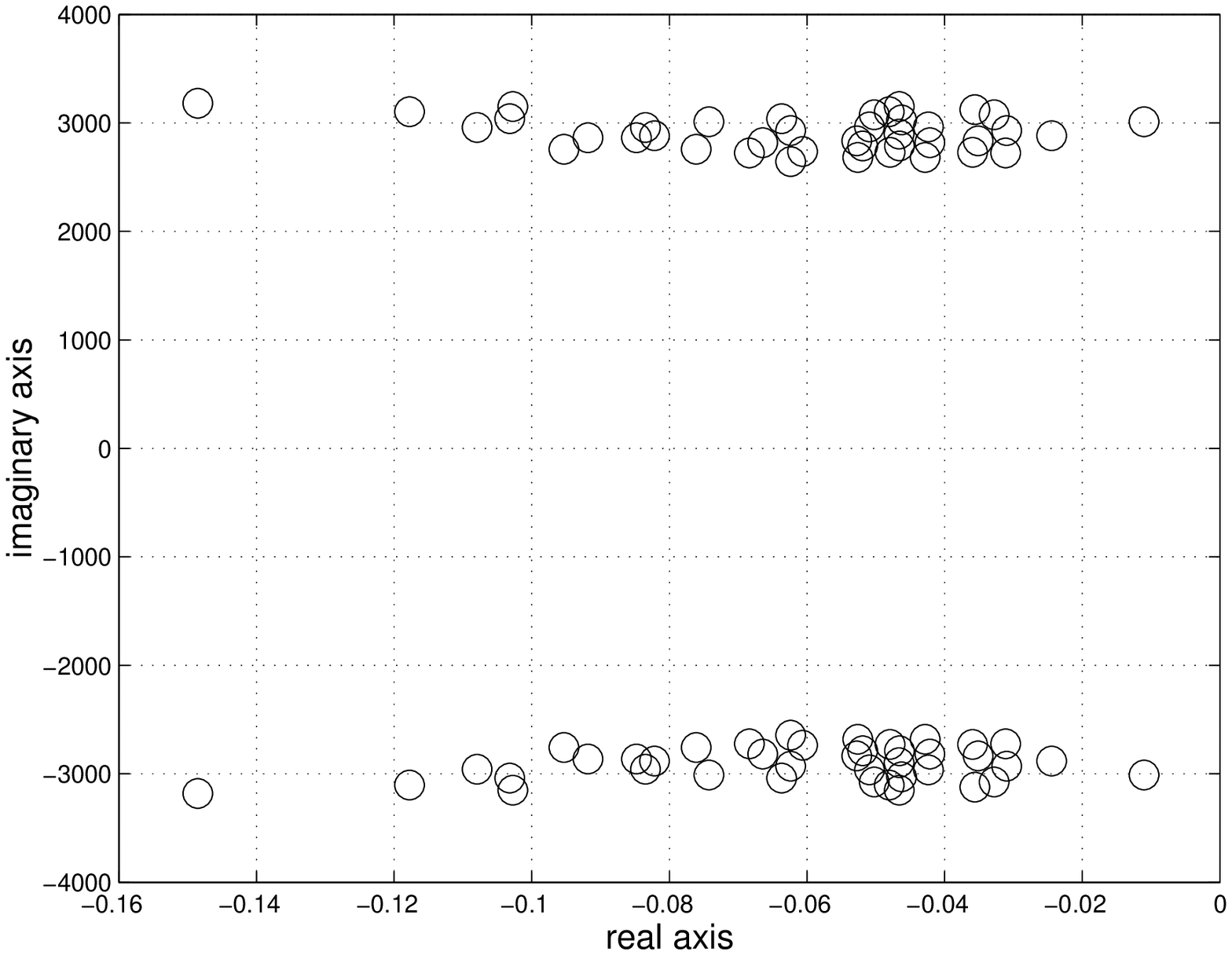}
\end{center}
\begin{center}
figure 16 : spectrum of the damped plate,  $a(x)=1,\omega  = ( 0.25,0.75 ) \times  ( 0.25,0.75 )$ (left), $a(x)=1,\omega  = ( 0.5,0.75 ) \times  ( 0.3,0.6 )$ (right).
\end{center}
\begin{center}
\includegraphics[angle=0,width=6.5cm]{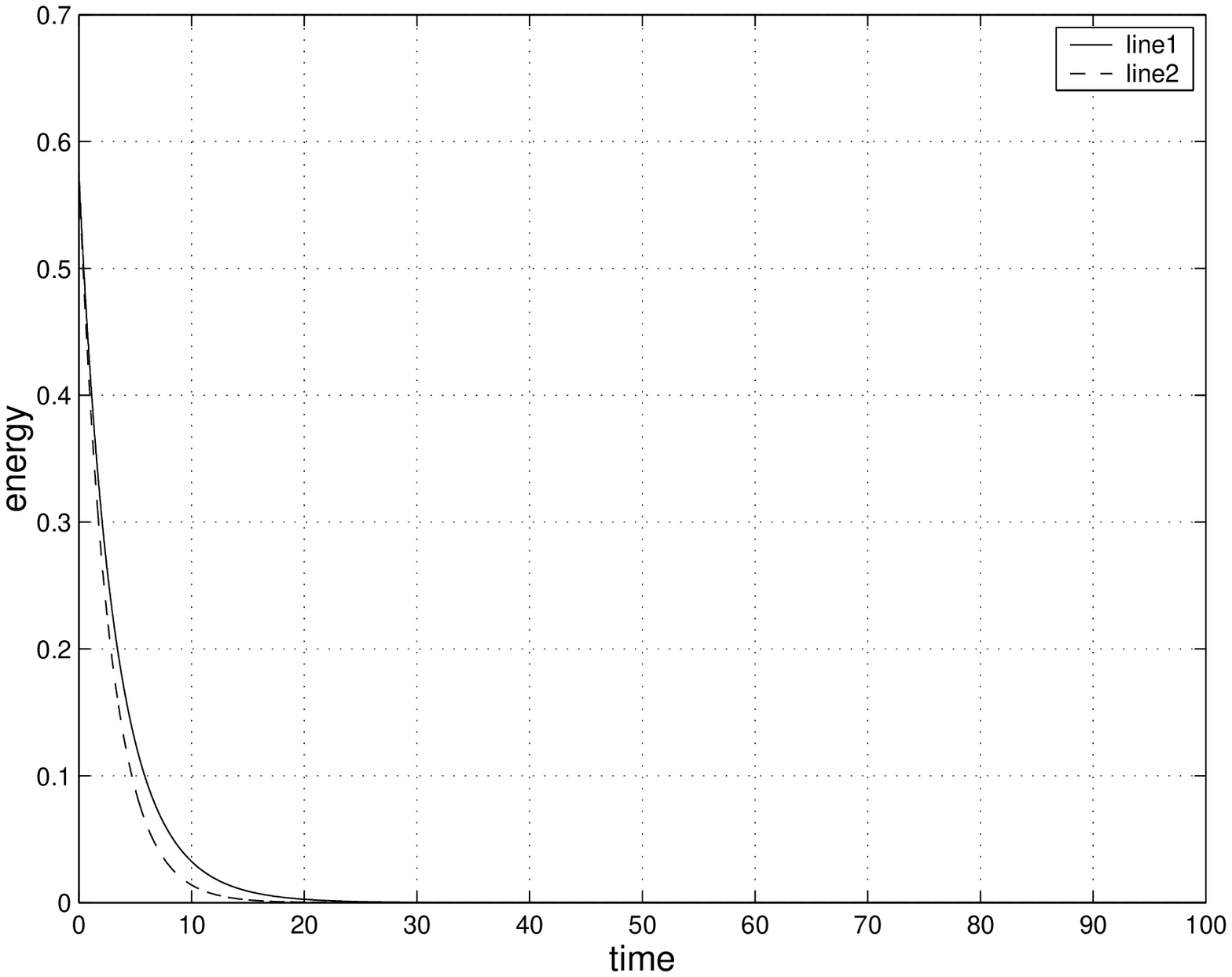}
\includegraphics[angle=0,width=6.5cm]{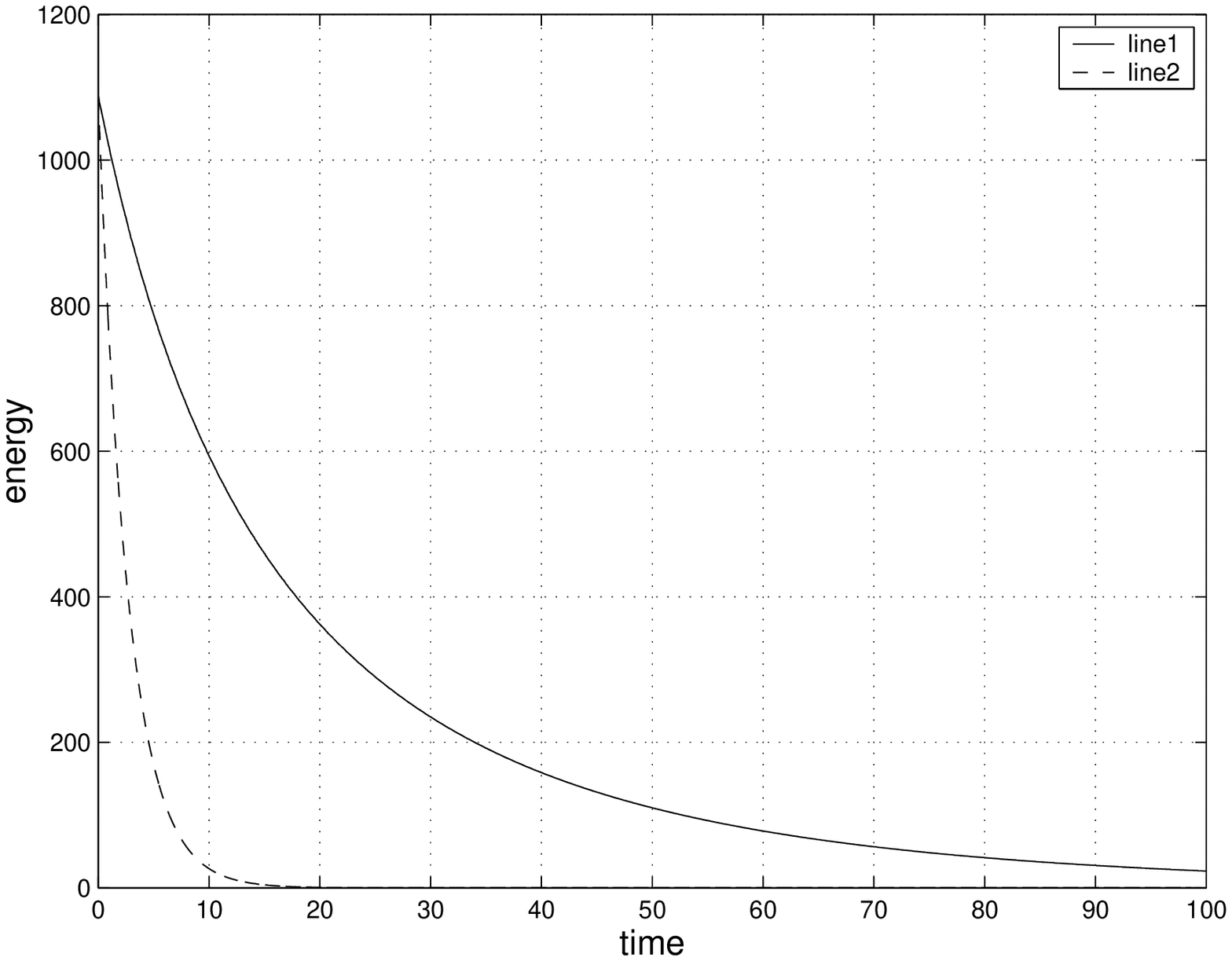}
\end{center}
\begin{center}
figure 17 : energy of the damped plate,  $a(x)=1,\omega  = ( 0.25,0.75 ) \times  ( 0.25,0.75 )$, $n=m=2$ (left), and $n=m=8$ (right).
\end{center}
In figures 13., 14. and 15., for the value  $a(x)=1$ we consider a damping in the domain  $\omega  = ( 0.4,0.6 ) \times  ( 0.4,0.6 )$ and $\omega  = ( 0.35,0.65 ) \times  ( 0.35,0.65 )$. In figures 16., 17. and 18. , we consider a damping in the domain $\omega  = ( 0.25,0.75 ) \times  ( 0.25,0.75 )$, and  $\omega  = ( 0.5,0.75 ) \times  ( 0.3,0.6 )$.
%
\begin{center}
\includegraphics[angle=0,width=6.5cm]{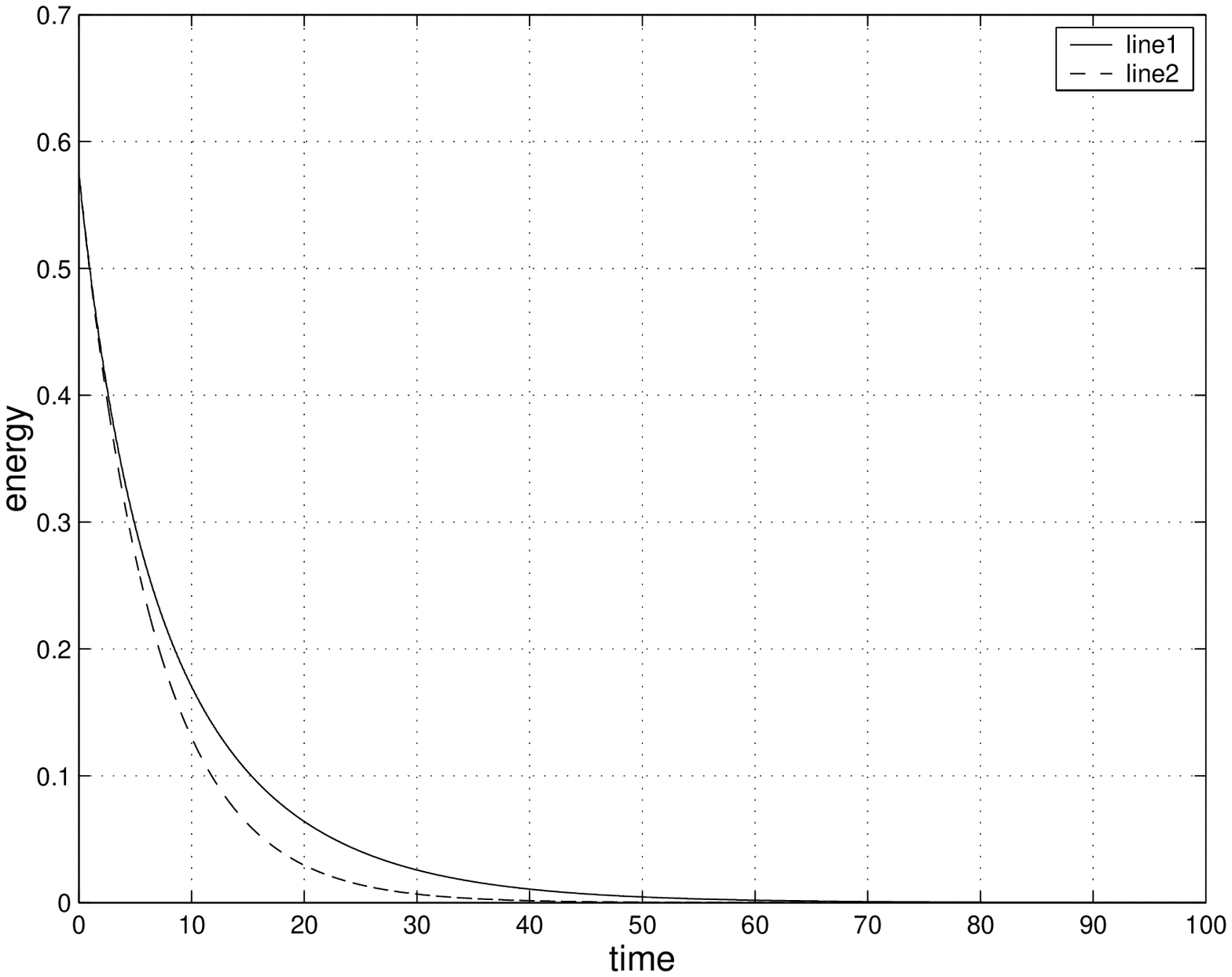}
\includegraphics[angle=0,width=6.5cm]{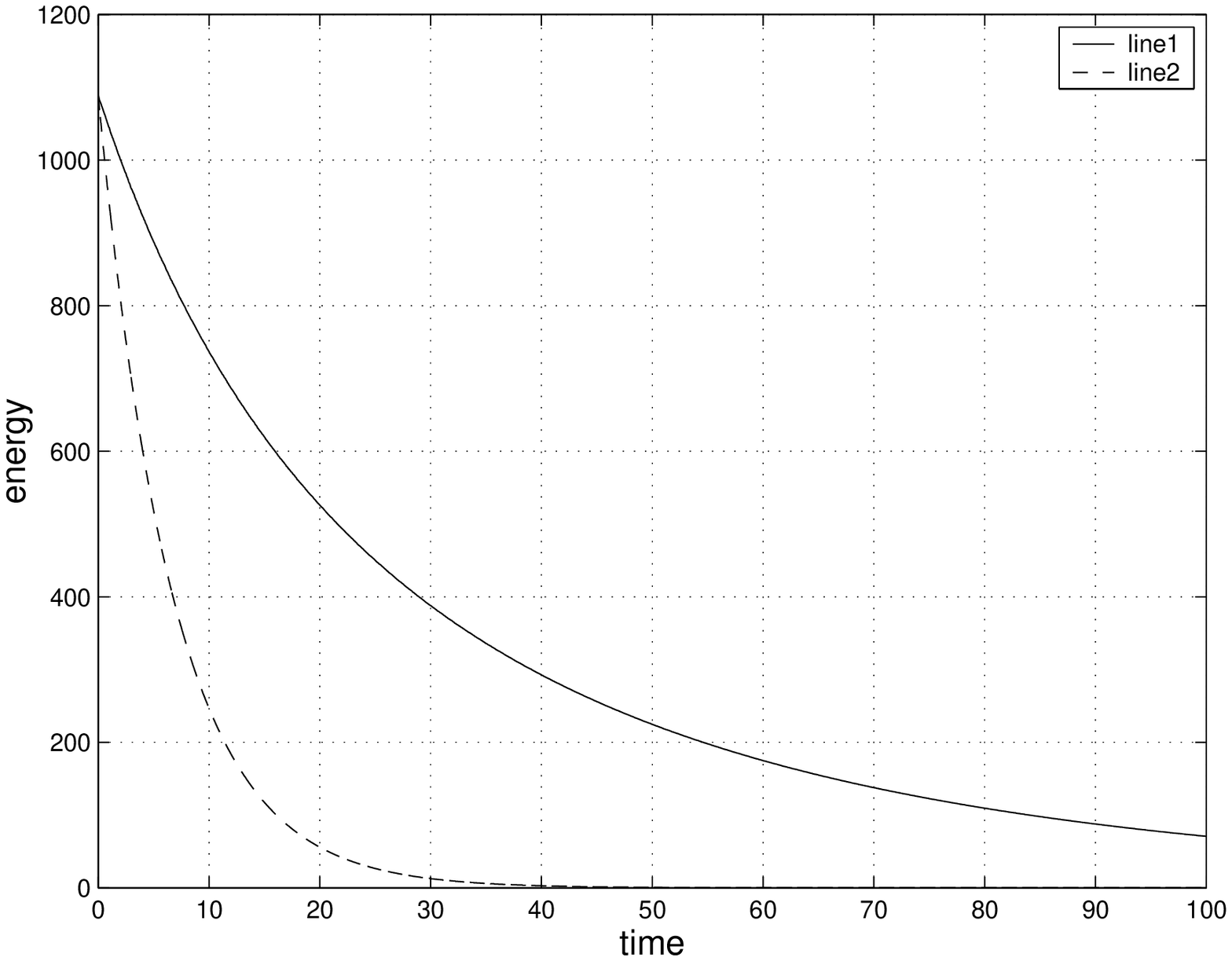}
\end{center}
\begin{center}
figure 18 : energy of the damped plate,  $a(x)=1,\omega  = ( 0.5,0.75 ) \times  ( 0.3,0.6 )$, $n=m=2$ (left),   and $n=m=8$ (right).
\end{center}
%
%
%
%
%
%
%

In the following examples the damping is a function given by $a(x_1,x_2)=x_1 x_2$. In the case of a square the spectrum of the damped plate is computed for three position of damping domain $\omega=\Omega$,$\omega=(0,0.5) \times (0,0.5)$ and $\omega=(0.35,0.75) \times (0.35,0.75)$
\begin{center}
\includegraphics[angle=0,width=7cm]{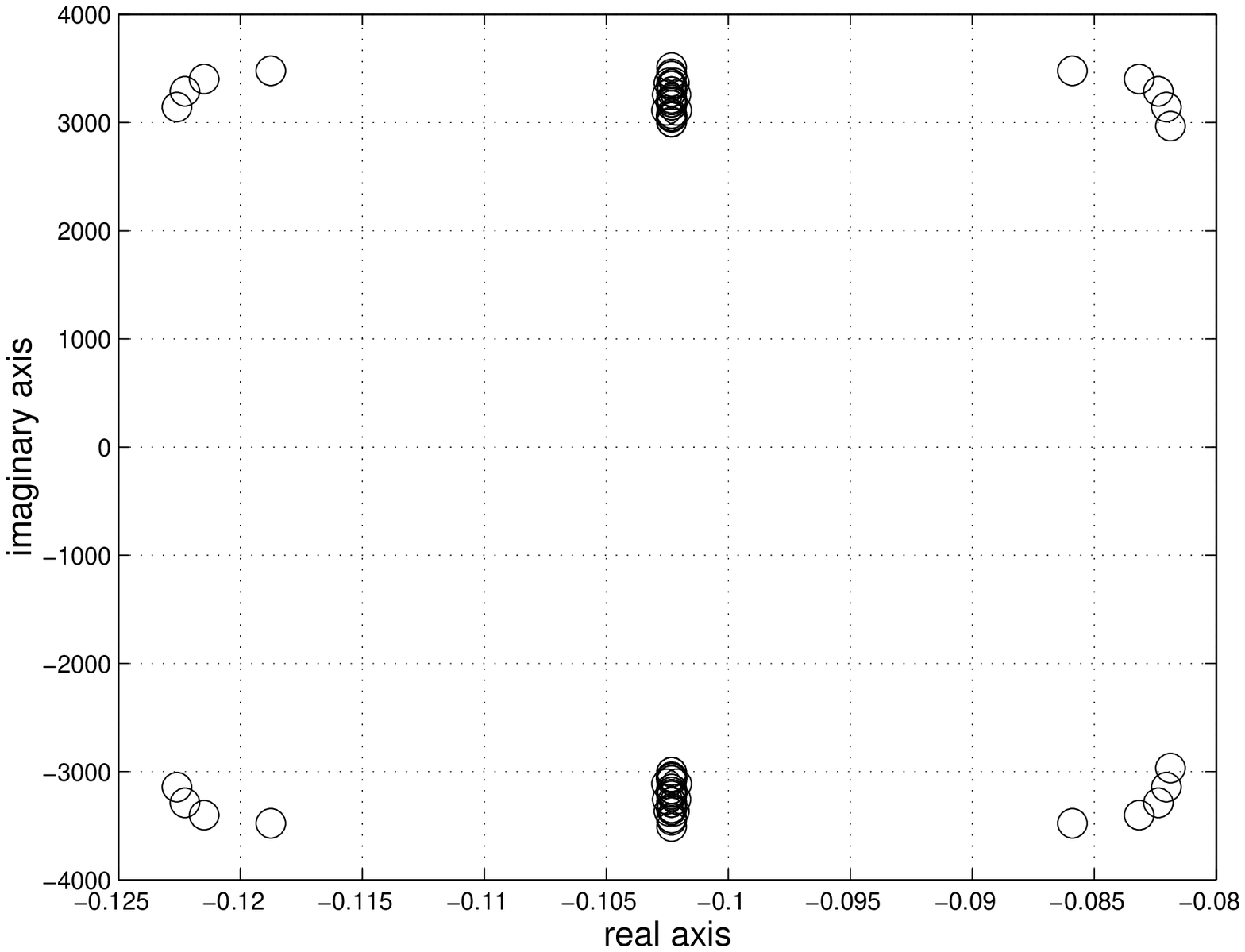}
\end{center}
\begin{center}
figure 19 : spectrum of the damped plate, $a(x_1,x_2)=x_1 x_2$, $\omega=\Omega$ .
\end{center}
\begin{center}
\includegraphics[angle=0,width=6.5cm]{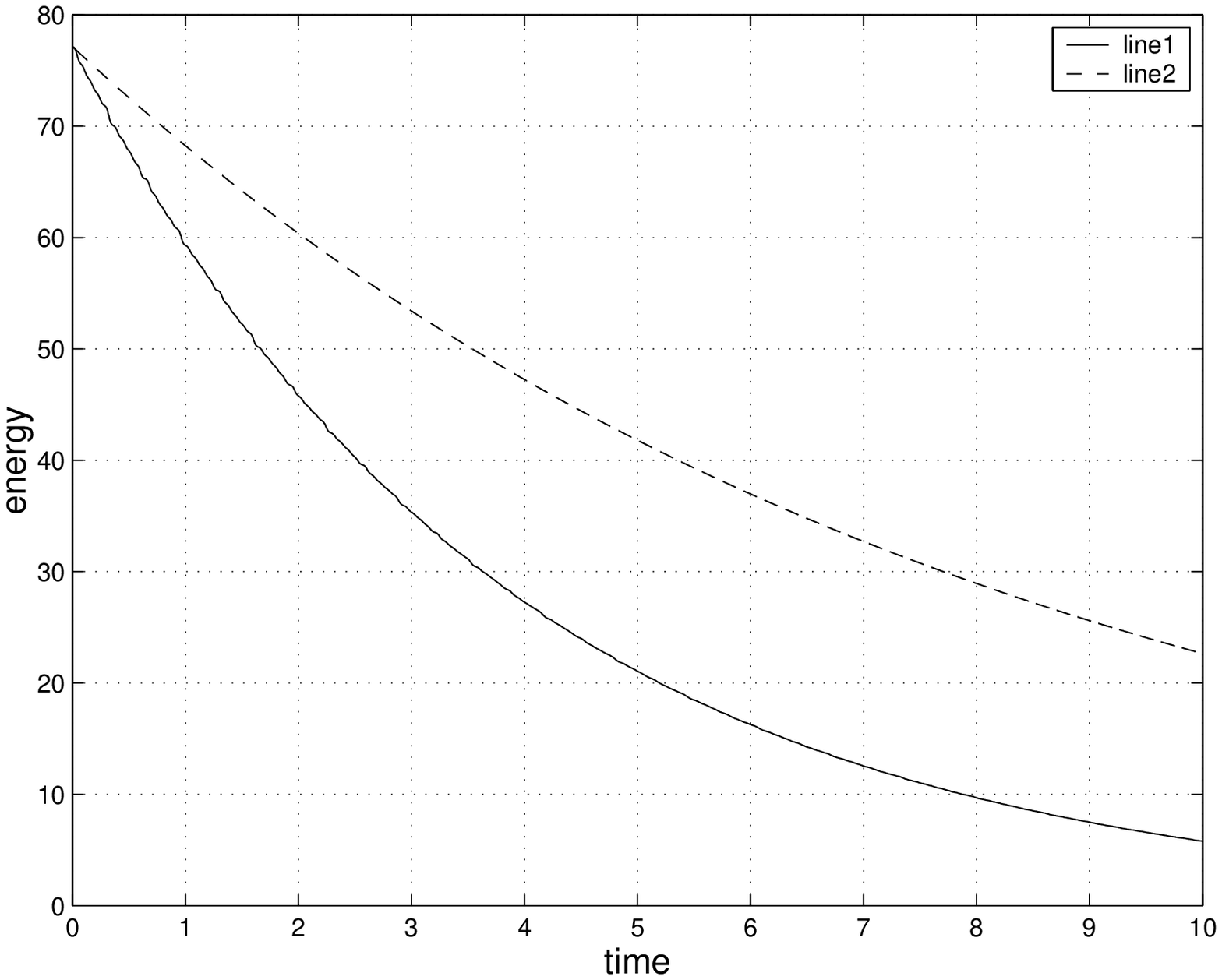}
\includegraphics[angle=0,width=6.5cm]{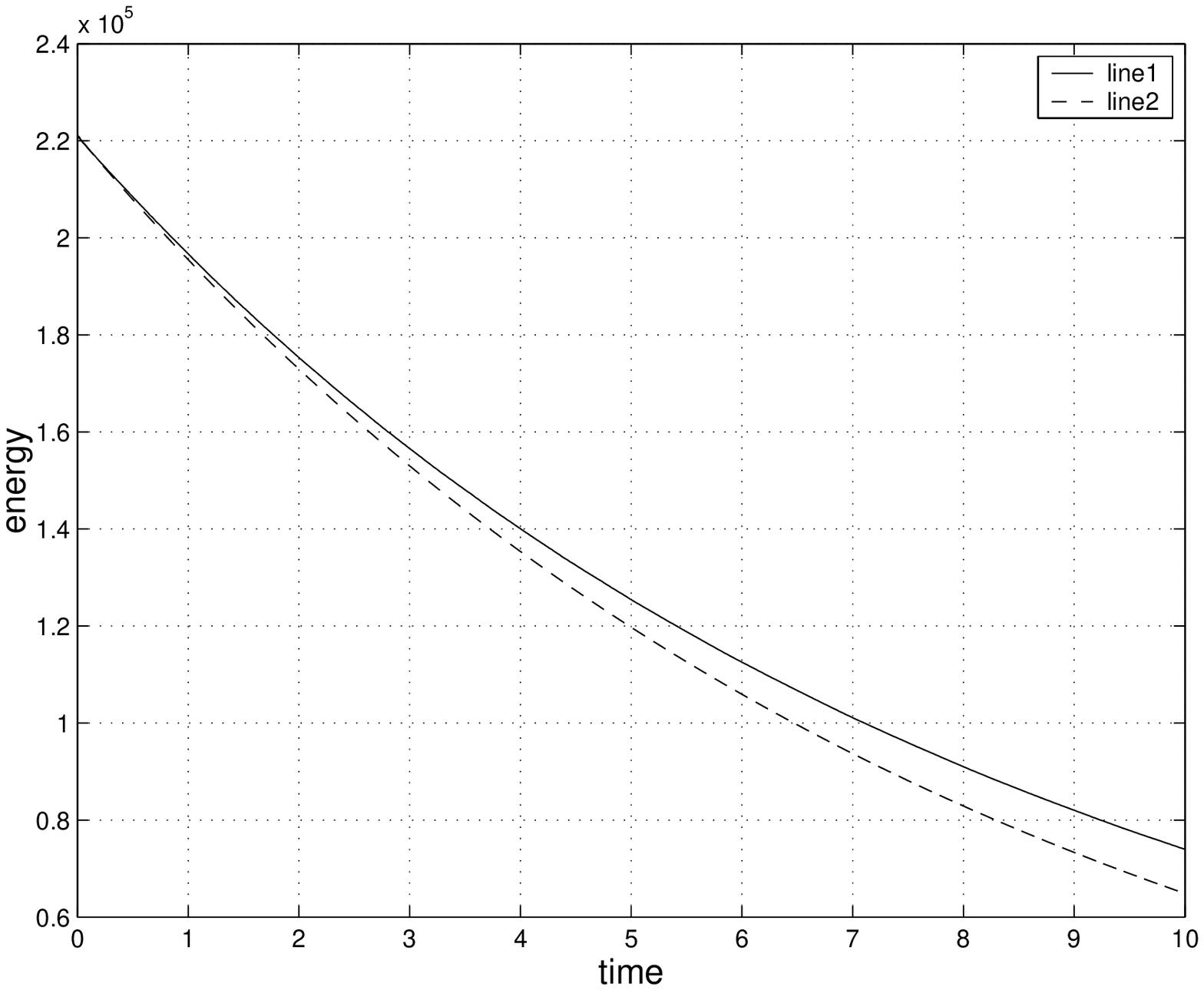}
\end{center}
\begin{center}
figure 20 : energy of the damped plate, $a(x_1,x_2)=x_1 x_2$, damping in all the domain $n=m=5$ (left), and 
 $n=m=20$ (right).
\end{center}
\begin{center}
\includegraphics[angle=0,width=6.5cm]{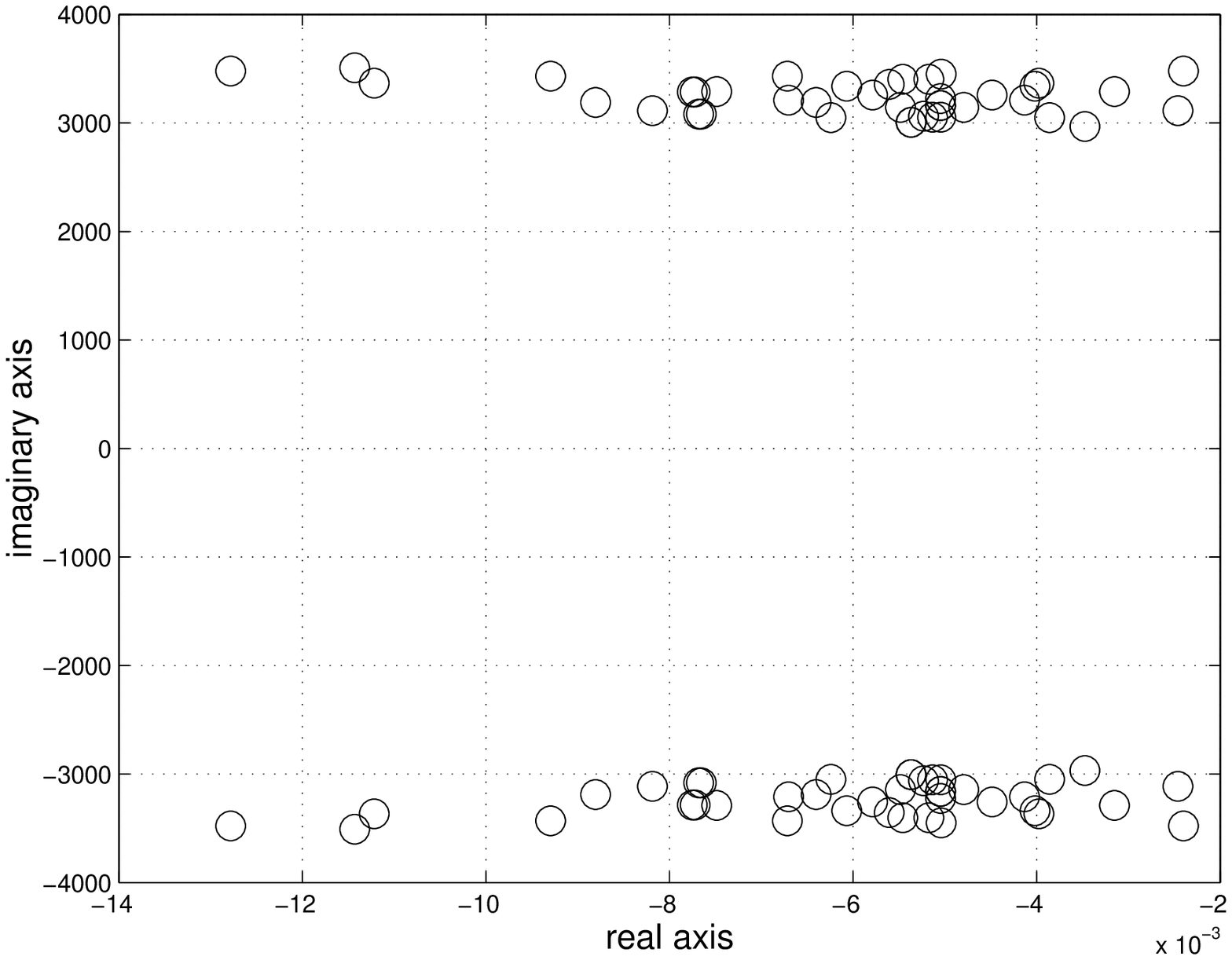}
\includegraphics[angle=0,width=6.5cm]{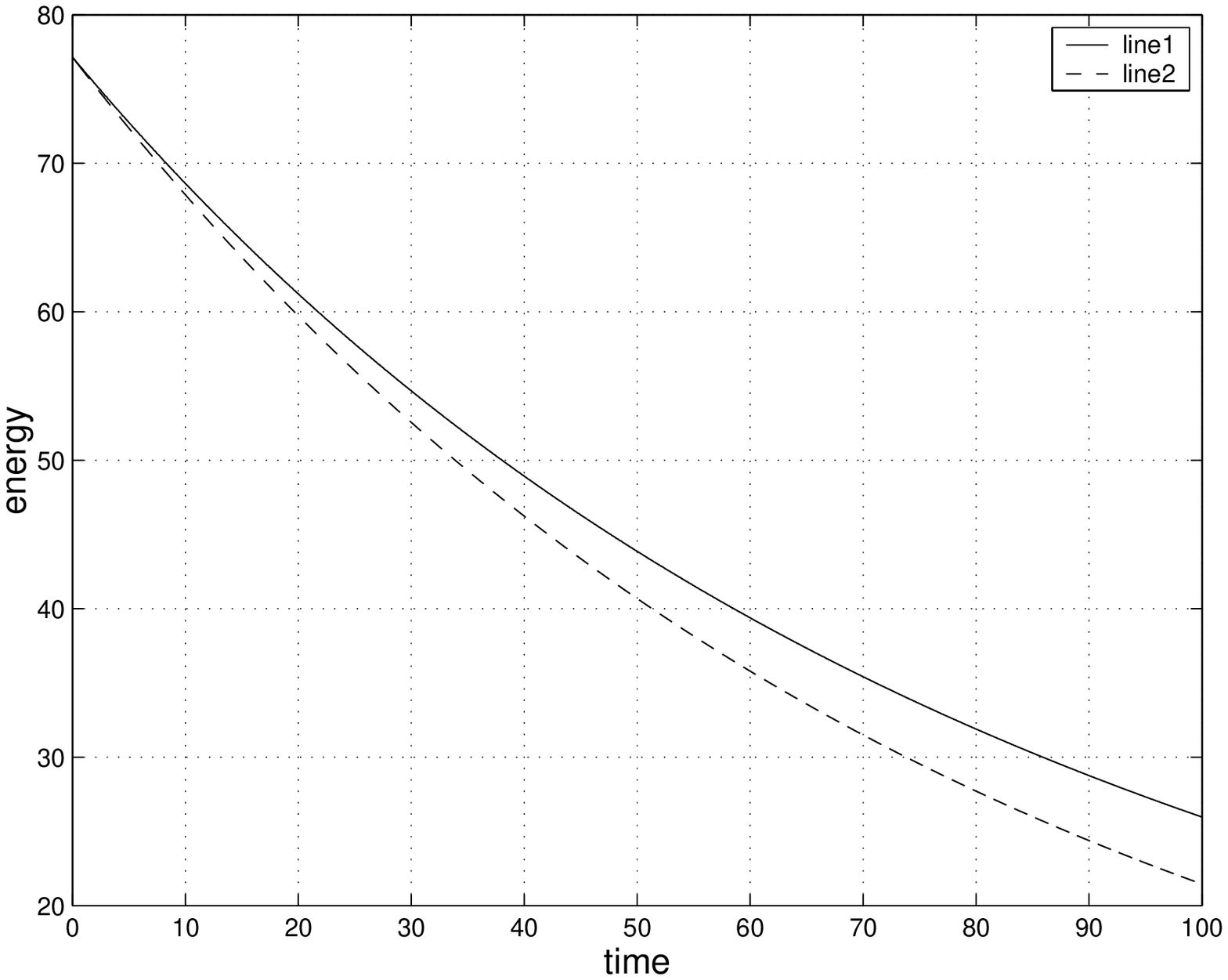}
\end{center}
\begin{center}
figure 21 : $a(x_1,x_2)=x_1 x_2$, $\omega = (0,0.5) \times (0,0.5)$ spectrum of the damped plate,  (left), energy of the damped plate,   $n=m=5$,  (right) .
\end{center}
\begin{center}
\includegraphics[angle=0,width=6.5cm]{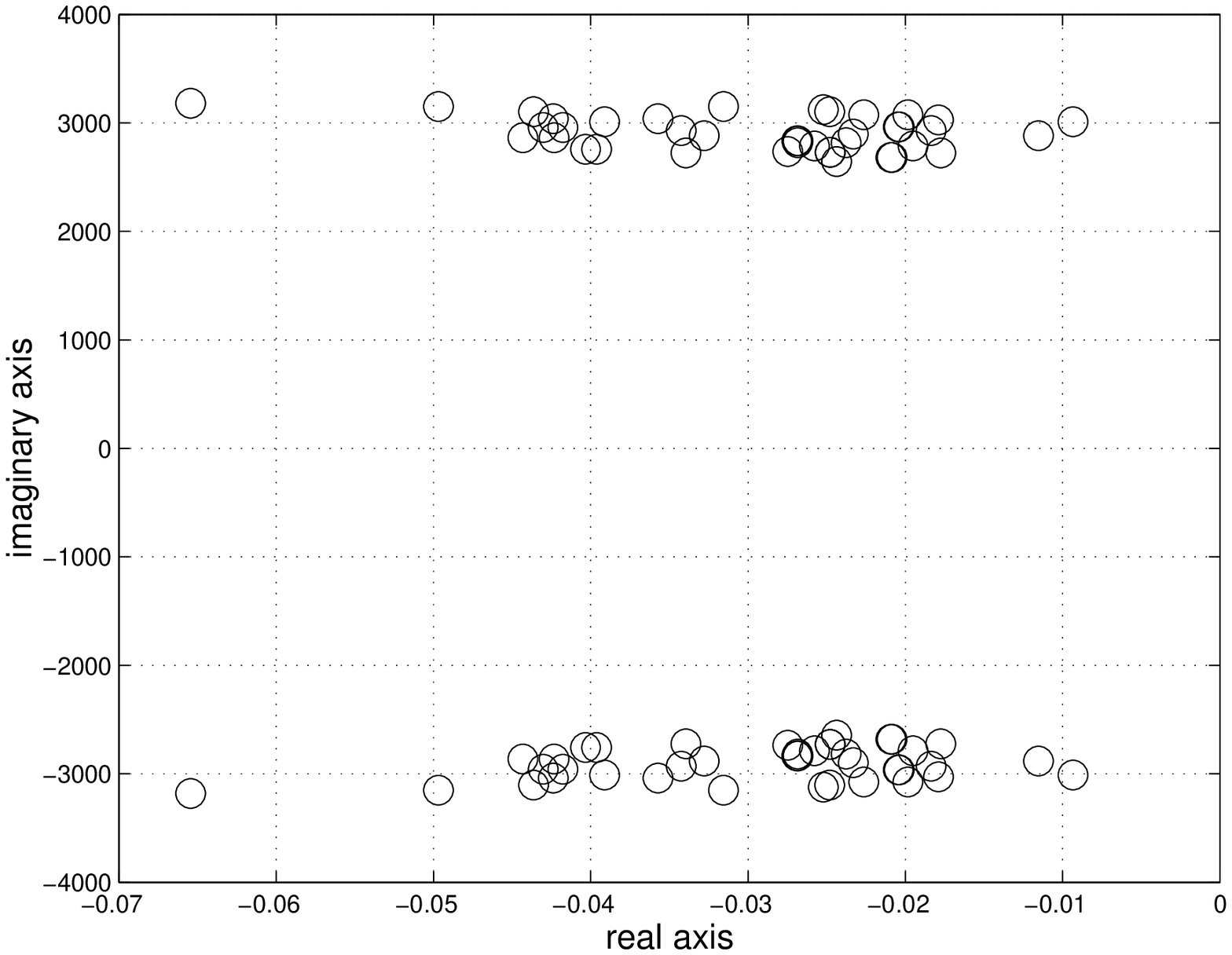}
\includegraphics[angle=0,width=6.5cm]{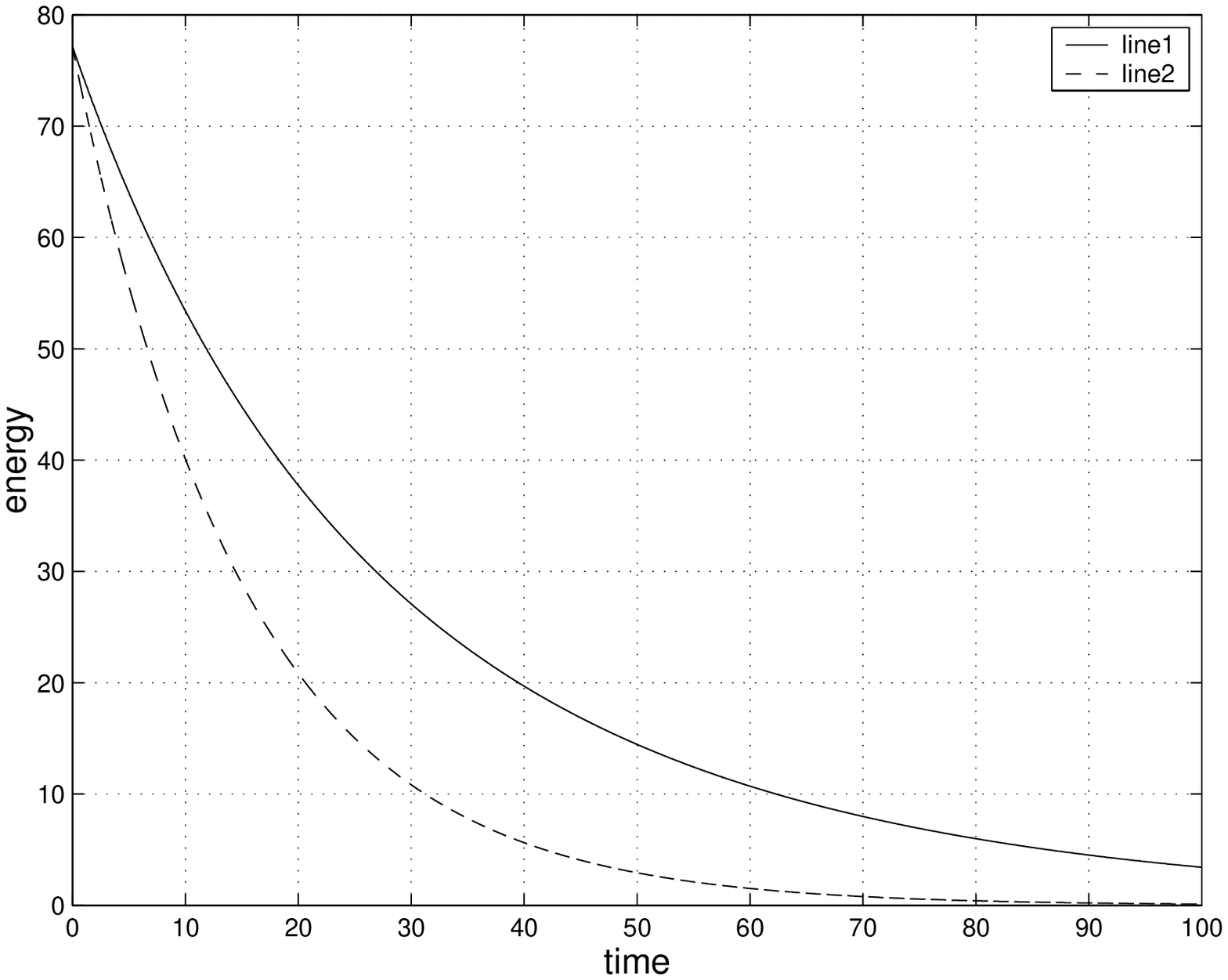}
\end{center}
\begin{center}
figure 22 :  $a(x_1,x_2)=x_1 x_2$,  $\omega = (0.35,0.75) \times (0.35,0.75)$ , spectrum of the damped plate,  (left) , energy of the damped plate,   $n=m=5$,  (right).
\end{center}
%
%
%
%
%
%
%
\vspace{1cm}
In the next example we consider a damping defined by $a(x_1,x_1) = sin(x_1) cos(x_2)$ for three region $\omega = \Omega$,
$\omega = (0,0.5) \times (0,0.5)$ and $\omega= (0.3 , 0.7) \times (0.3 , 0,7)$. We plot the spectrum in the three situations and the energy of the damped plate compared to $E_0 e^{Re(\lambda_0) t}$ (line 2 : dashed line). figures 23., 24., and figure 25.
\vspace{1cm}
\begin{center}
\includegraphics[angle=0,width=6.5cm]{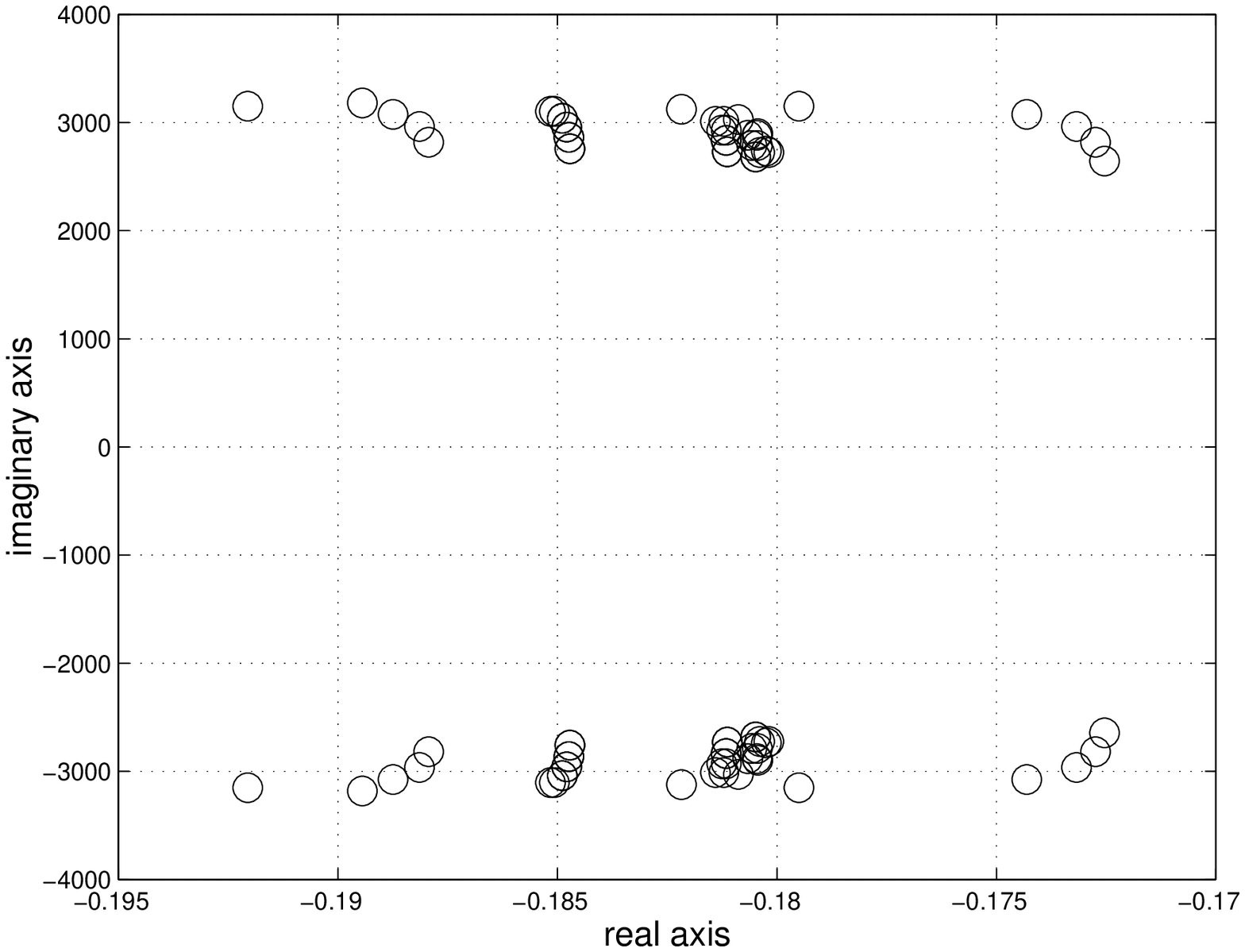}
\includegraphics[angle=0,width=6.5cm]{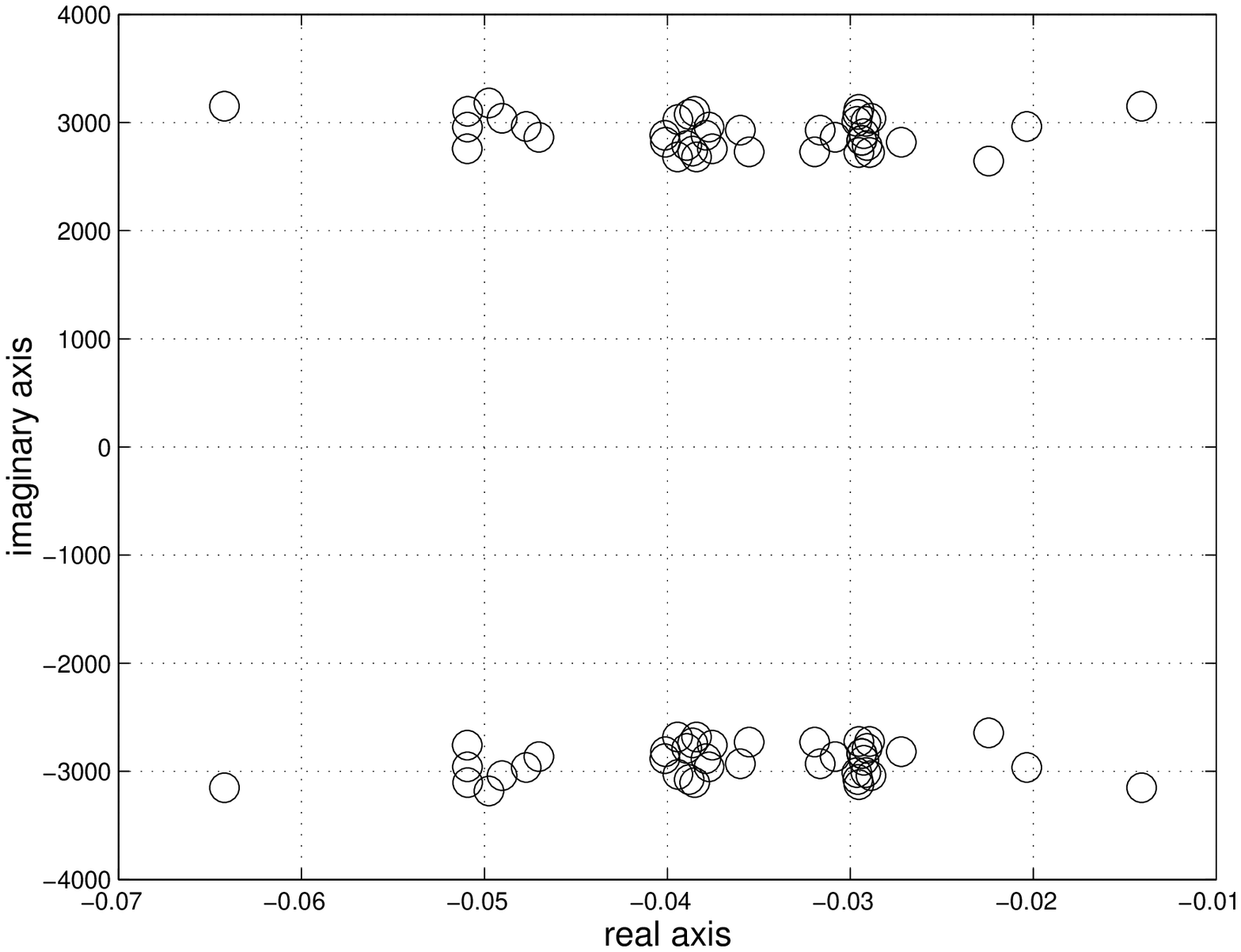}
\end{center}
\begin{center}
figure 23 : spectrum of the damped plate, $a(x_1,x_2)=sin(x_1) cos(x_2)$,$\omega = \Omega$ (left) , $\omega = (0,0.5) \times (0,0.5)$ (right).
\end{center}
\begin{center}
\includegraphics[angle=0,width=6.5cm]{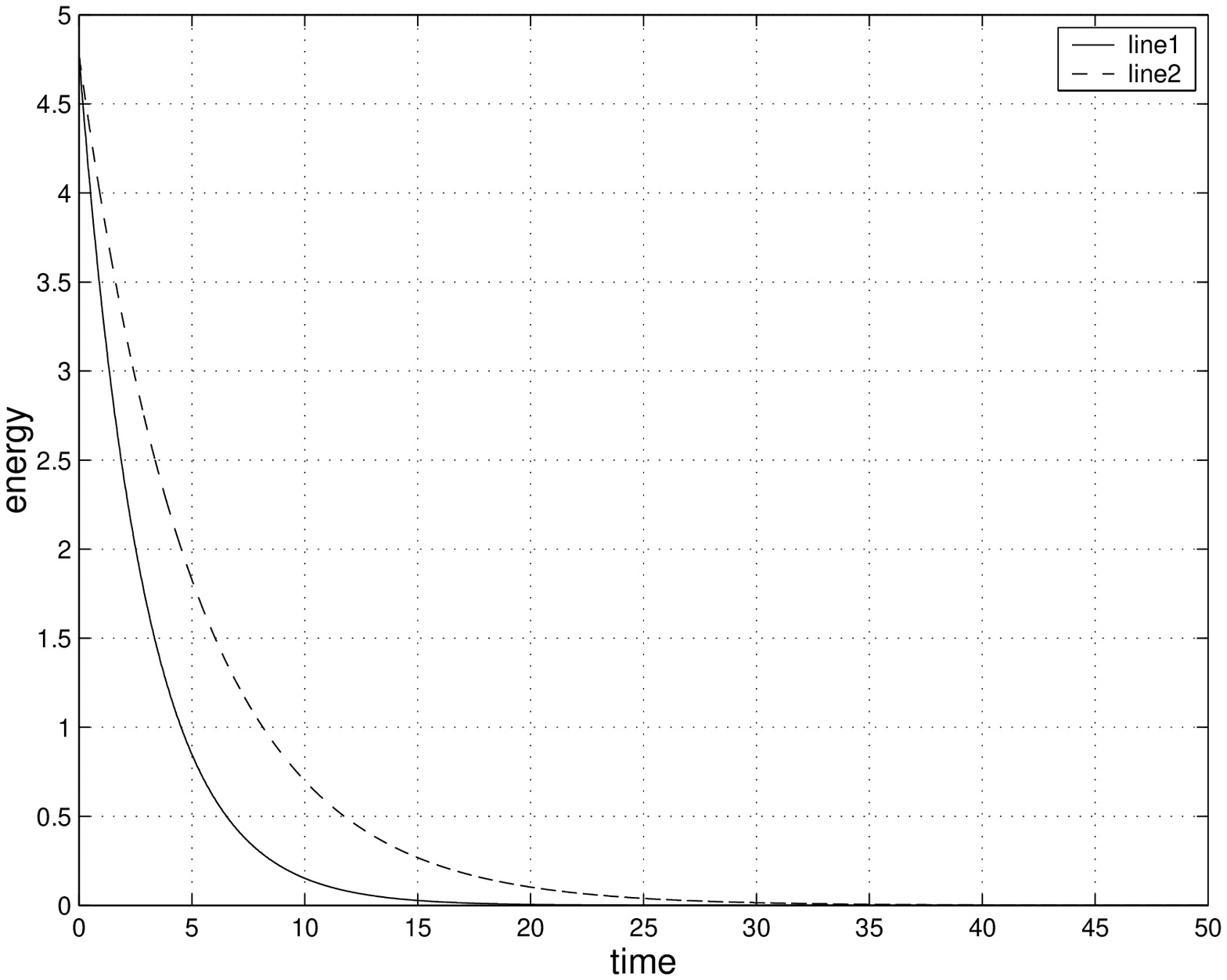}
\includegraphics[angle=0,width=6.5cm]{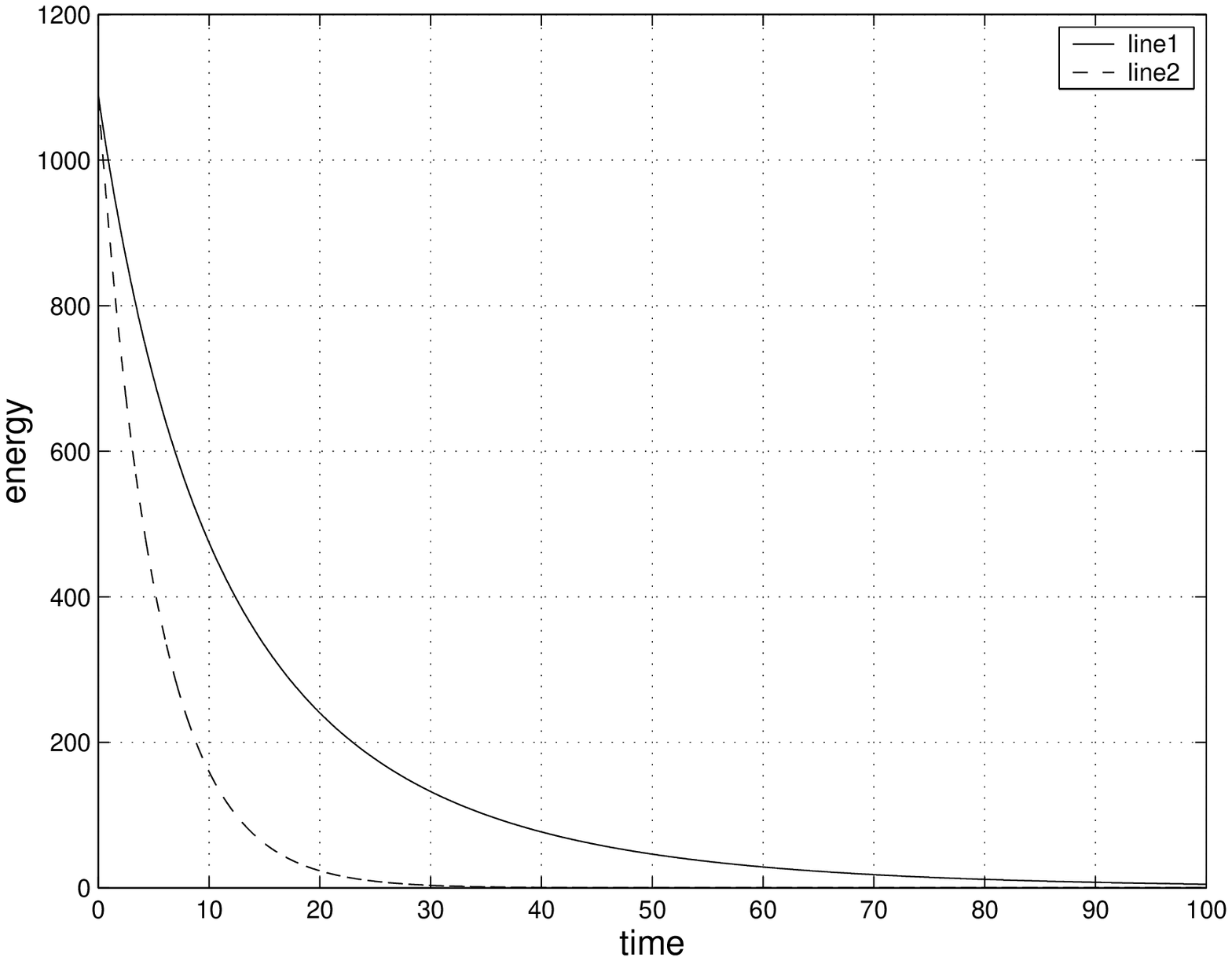}
\end{center}
\begin{center}
figure 24 : energy of the damped plate, $a(x_1,x_2)=sin(x_1) cos(x_2)$,$\omega = \Omega$, $n=m=3$ (left), and $n=m=8$ (right).
\end{center}
%
%
%
%
%
\begin{center}
\includegraphics[angle=0,width=6.5cm]{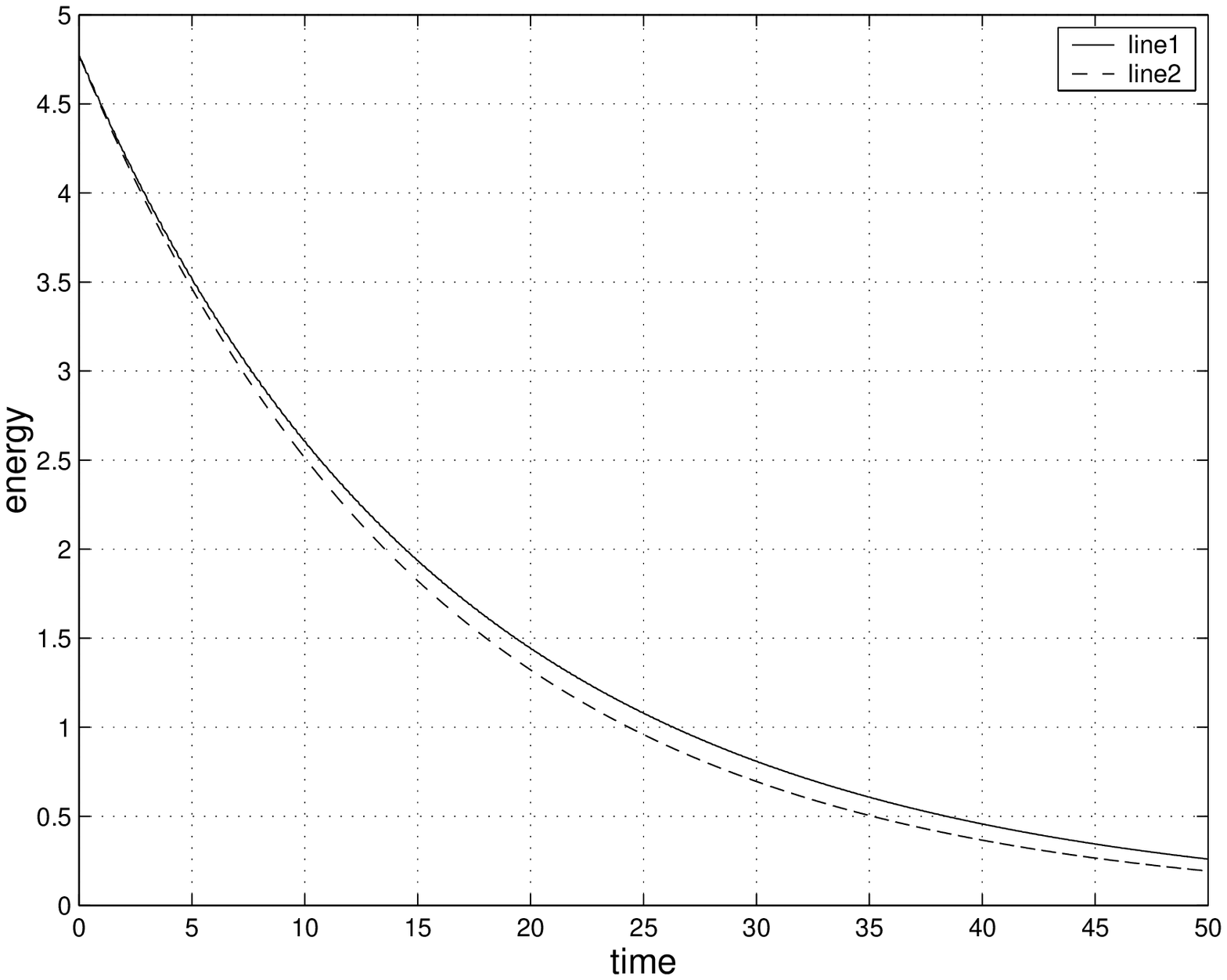}
\includegraphics[angle=0,width=6.5cm]{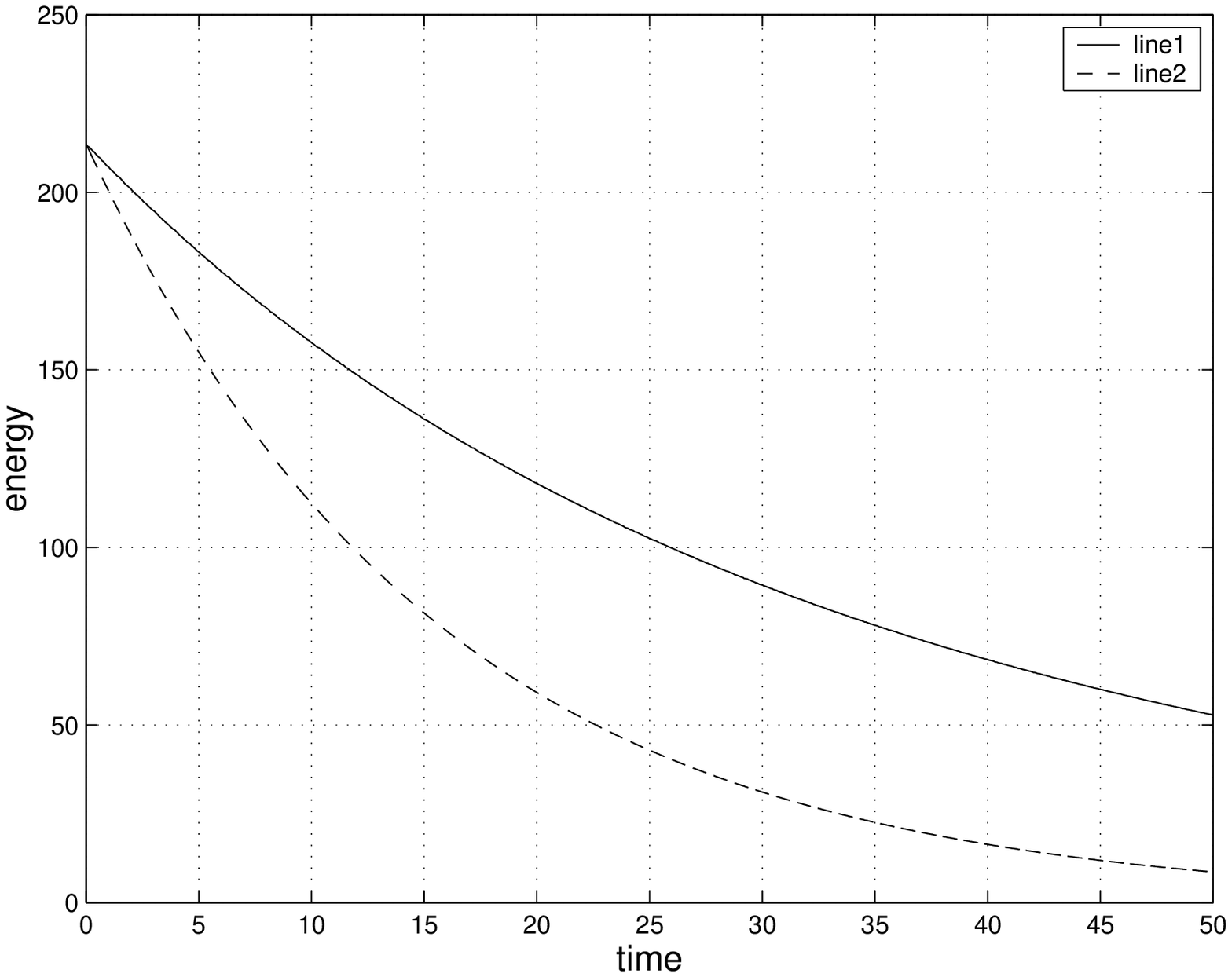}
\end{center}
\begin{center}
figure 25: energy of the damped plate, $a(x_1,x_2)=sin(x_1) cos(x_2)$,$\omega = (0,0.5) \times (0,0.5)$, $n=m=3$(left), and $n=m=6$ (right).
\end{center}
%
%
%

\section{The case of rectangle}
%
In this case we consider rectangular a domain $\Omega=(0,a) \times (0,b)$ with a set of normalized eigenfunctions $\Phi_{n,m}(x,y) = \sqrt \frac{ 2}{ a } sin(\frac{n \pi x}{a}) \sqrt \frac{ 2}{ b } sin(\frac{m \pi y}{b})$,
for all $(x,y) \in \Omega $ the solution $u(x,y,t)$
of the problem (\ref{eq1})-(\ref{eq2}) is given by :
\be \label{eq4}
u(x,y,t)=\sum_{n,m} \alpha_{n,m}(t) \Phi_{n,m}(x,y)
\ee
with :
\be
\frac{d^2 \alpha_{n,m}(t)}{dt^2} +  ( (\frac{n \pi}{a} )^2 + (\frac{m \pi}{b} )^2 )^2 \alpha_{n,m} (t) + \frac{d \alpha_{n,m}(t)}{dt} = 0 
\ee
as in the case of square we compute the spectrum and energy of the first eigenmodes for different domains $\omega$. We compare in each case the energy of the solution (line 1) to $E_0 e^{Re(\lambda_0) t}$  (line 2), with $Re(\lambda_0) = inf \{Re(\lambda)/ \lambda \in \sigma(A) \}$.
\begin{center}
\includegraphics[angle=0,width=6.5cm]{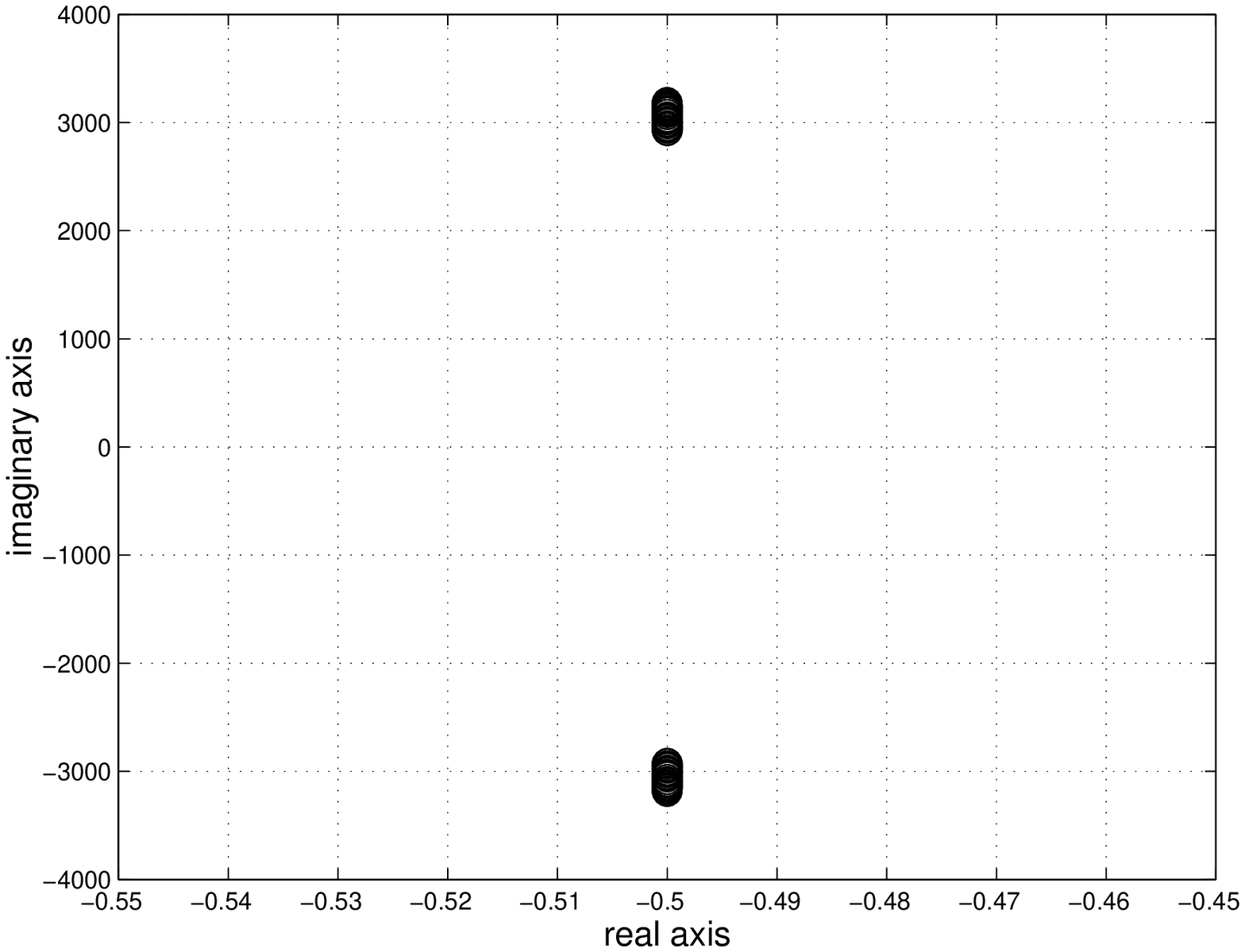}
\includegraphics[angle=0,width=6.5cm]{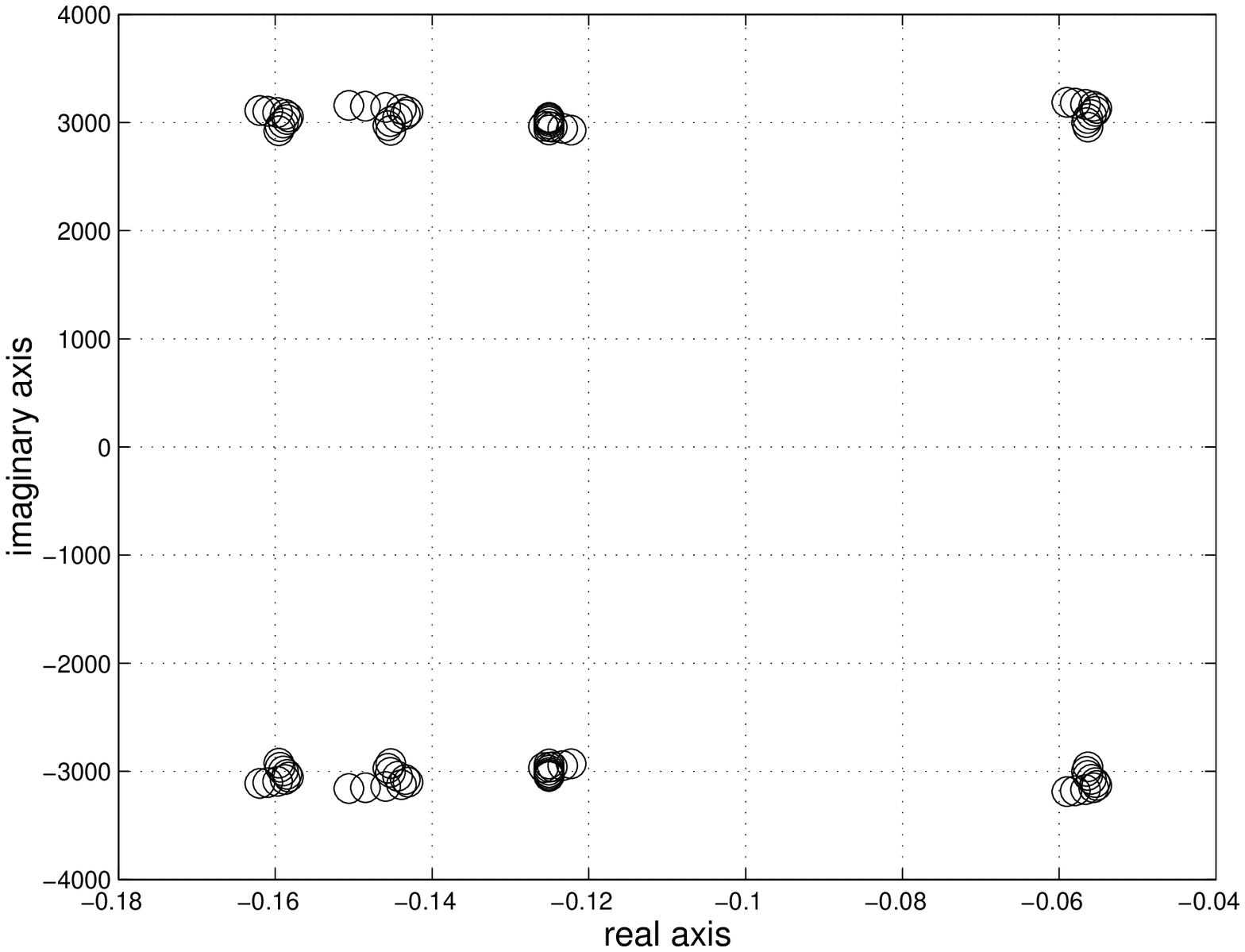}
\end{center}
\begin{center}
figure 26 : spectrum of the damped plate $a(x)=1$,$\omega = \Omega $ (left), and $\omega  = (0, \frac{1}{2}) \times  (0,1)$ (right).
\end{center}
\begin{center}
\includegraphics[angle=0,width=6.5cm]{figtest00ener3.3.ps}
\includegraphics[angle=0,width=6.5cm]{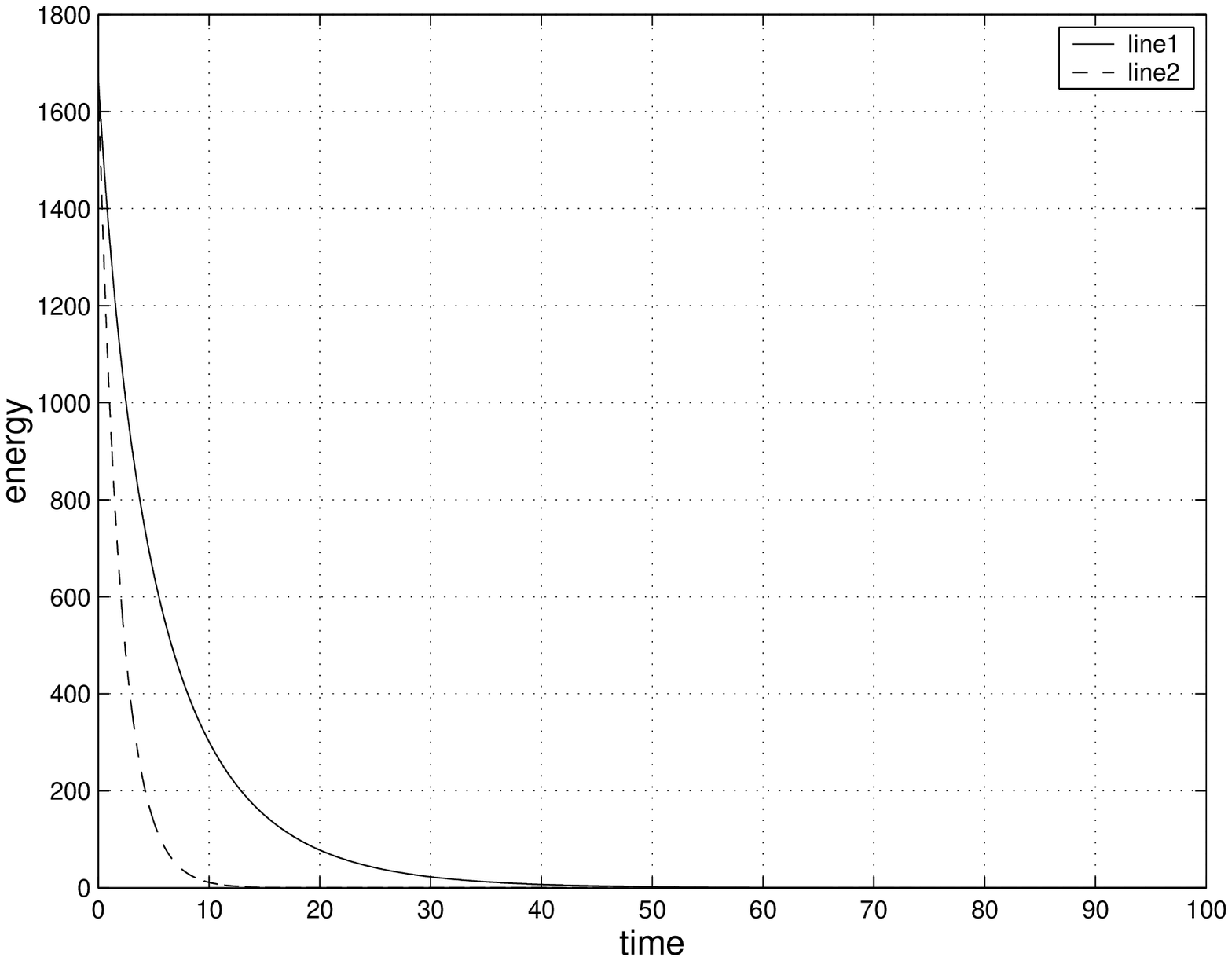}
\end{center}
\begin{center}
figure 27 : energy of the damped plate ,$a(x)=1$, $n=m=3$ (left),  and $n=m=12$ (right).
\end{center}
\begin{center}
\includegraphics[angle=0,width=6.5cm]{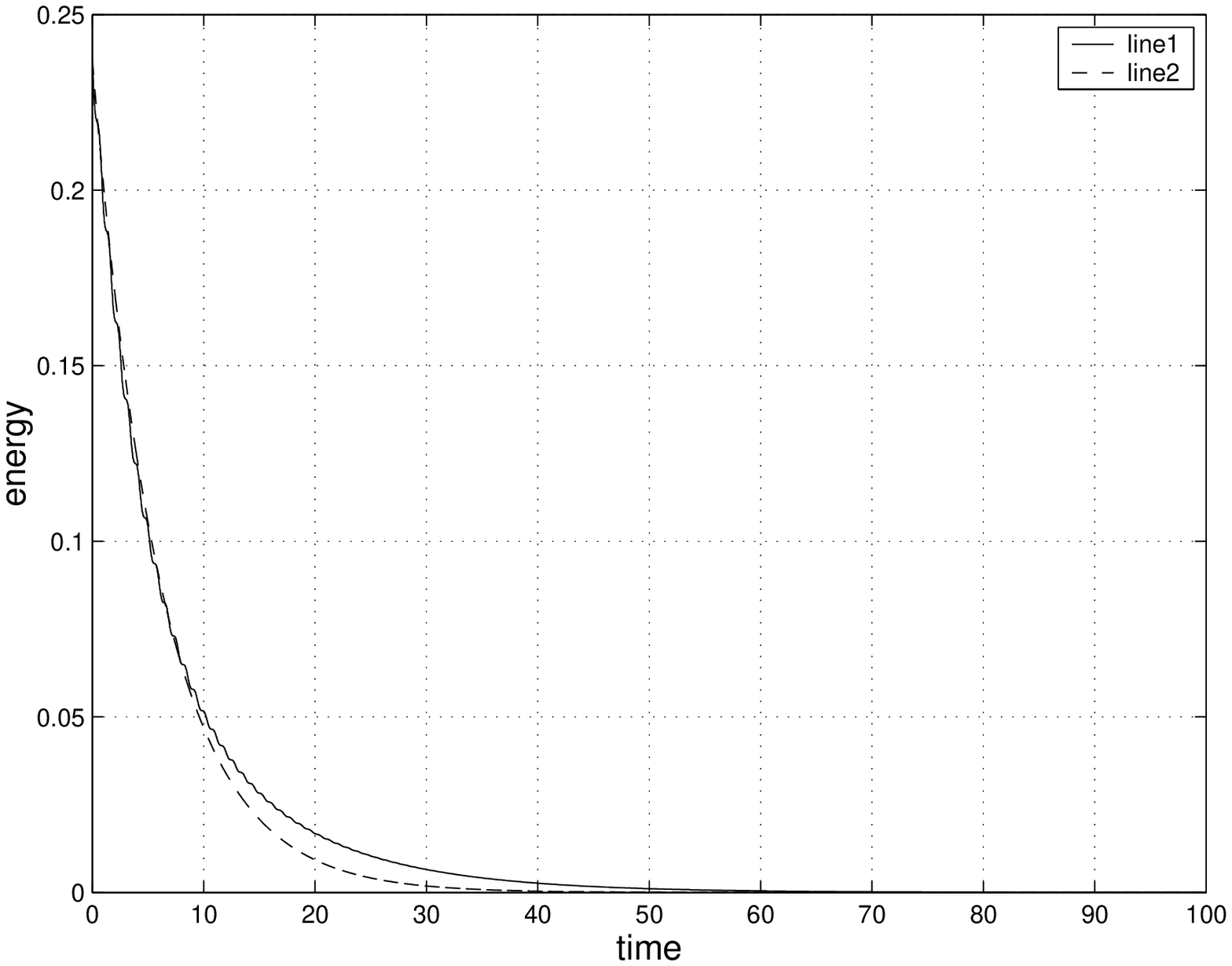}
\includegraphics[angle=0,width=6.5cm]{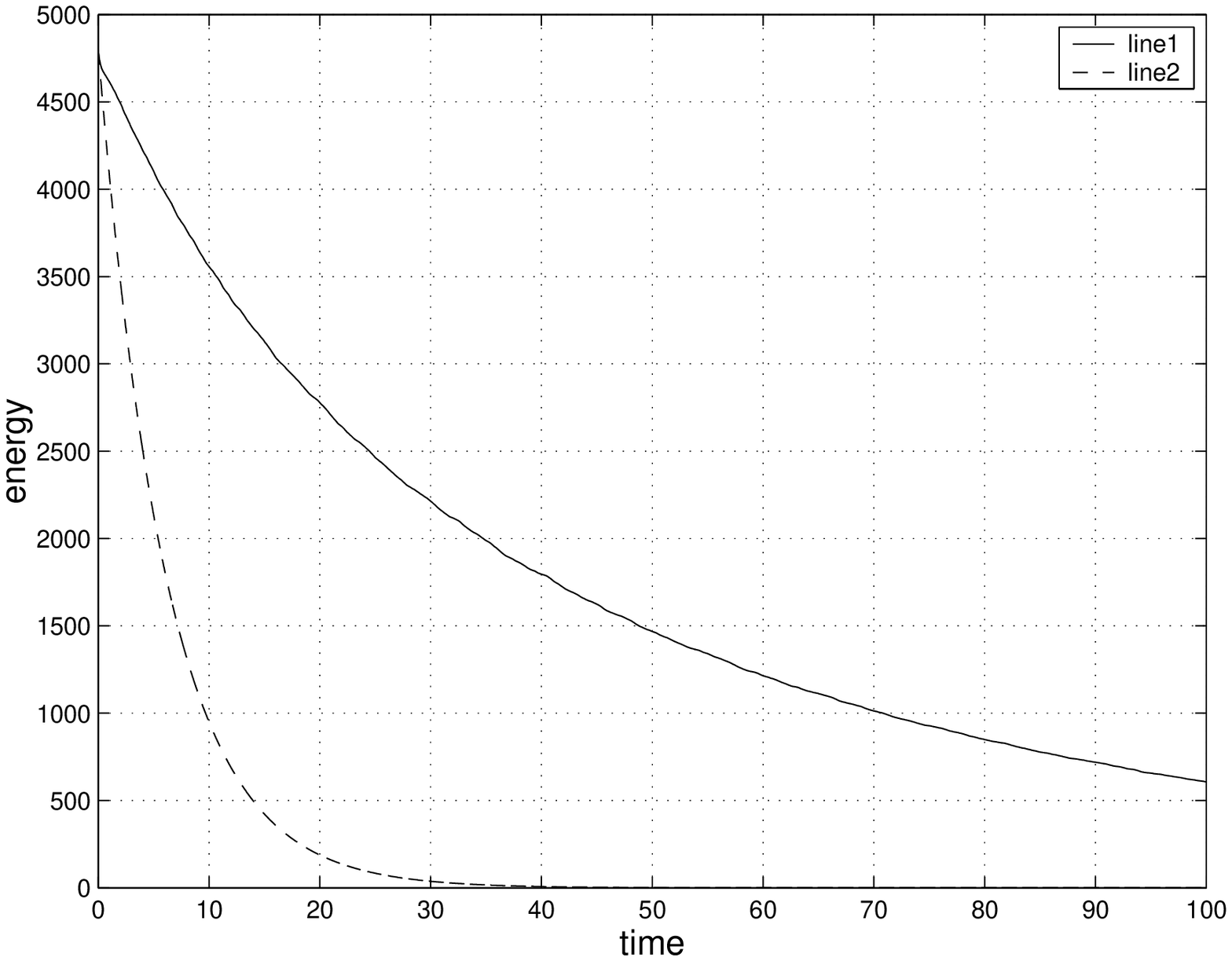}
\end{center}
\begin{center}
figure 28 : energy of the damped plate, $a(x)=1,\omega  = (0, \frac{1}{2}) \times  (0,1)$, $n=m=2$ (left),  and $n=m=10$ (right).
\end{center}
\begin{center}
\includegraphics[angle=0,width=6.5cm]{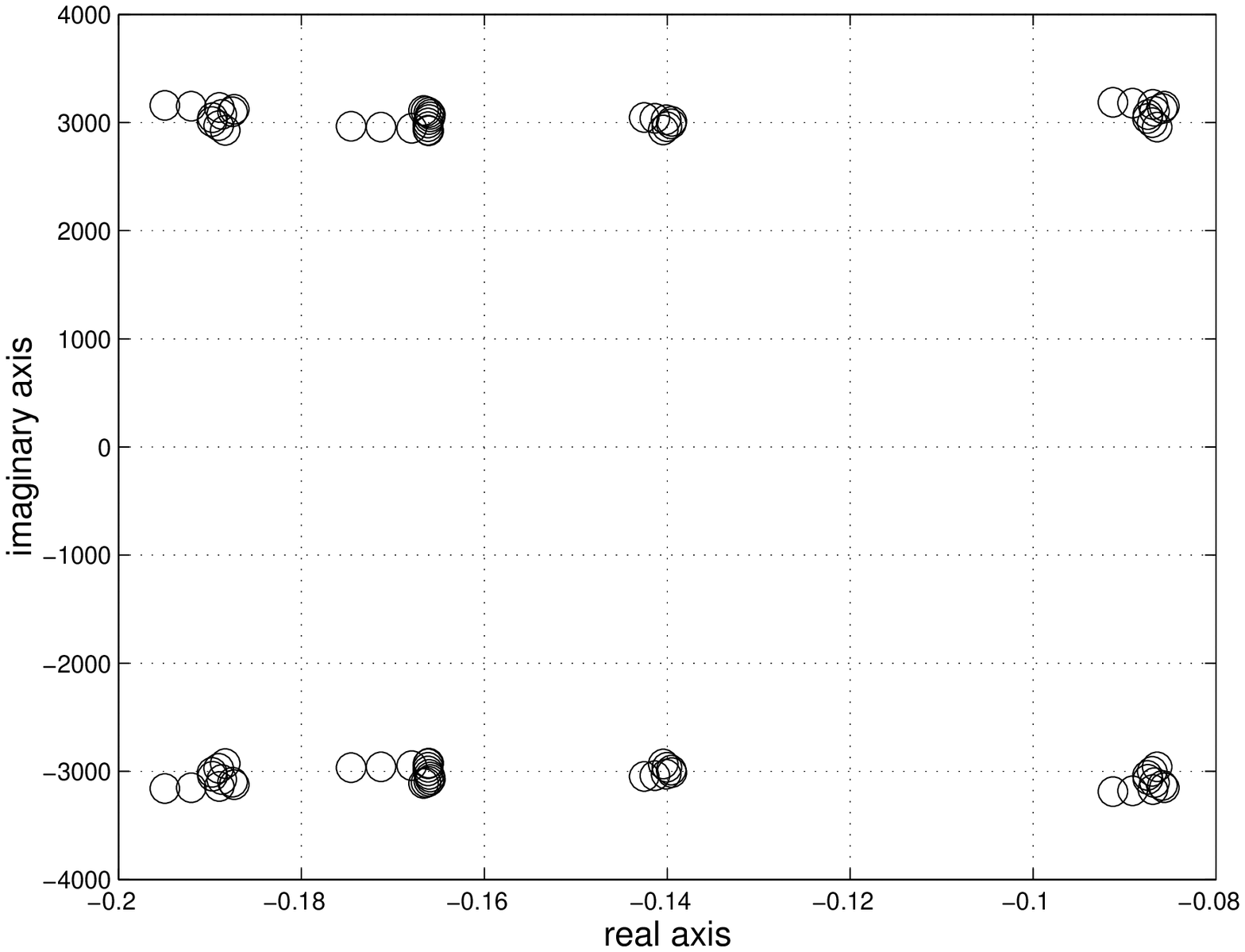}
\includegraphics[angle=0,width=6.5cm]{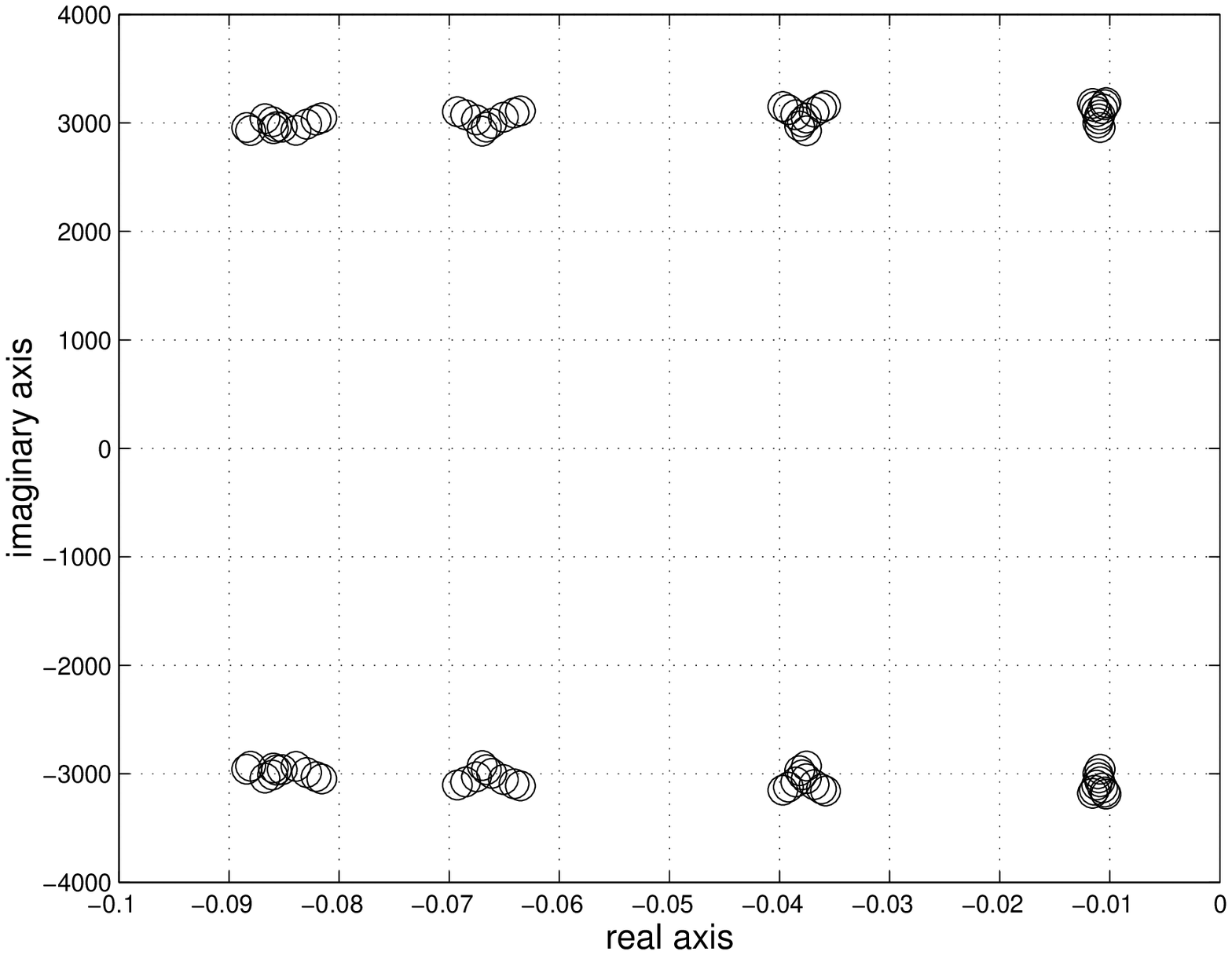}
\end{center}
\begin{center}
figure 29 : spectrum of the damped plate  $a(x)=1$,$,\omega  = (0, 0.6) \times  (0,1)$, (left) , $\omega  = (0, \frac{1}{5}) \times  (0,1)$ (right).
\end{center}
\begin{center}
\includegraphics[angle=0,width=6.5cm]{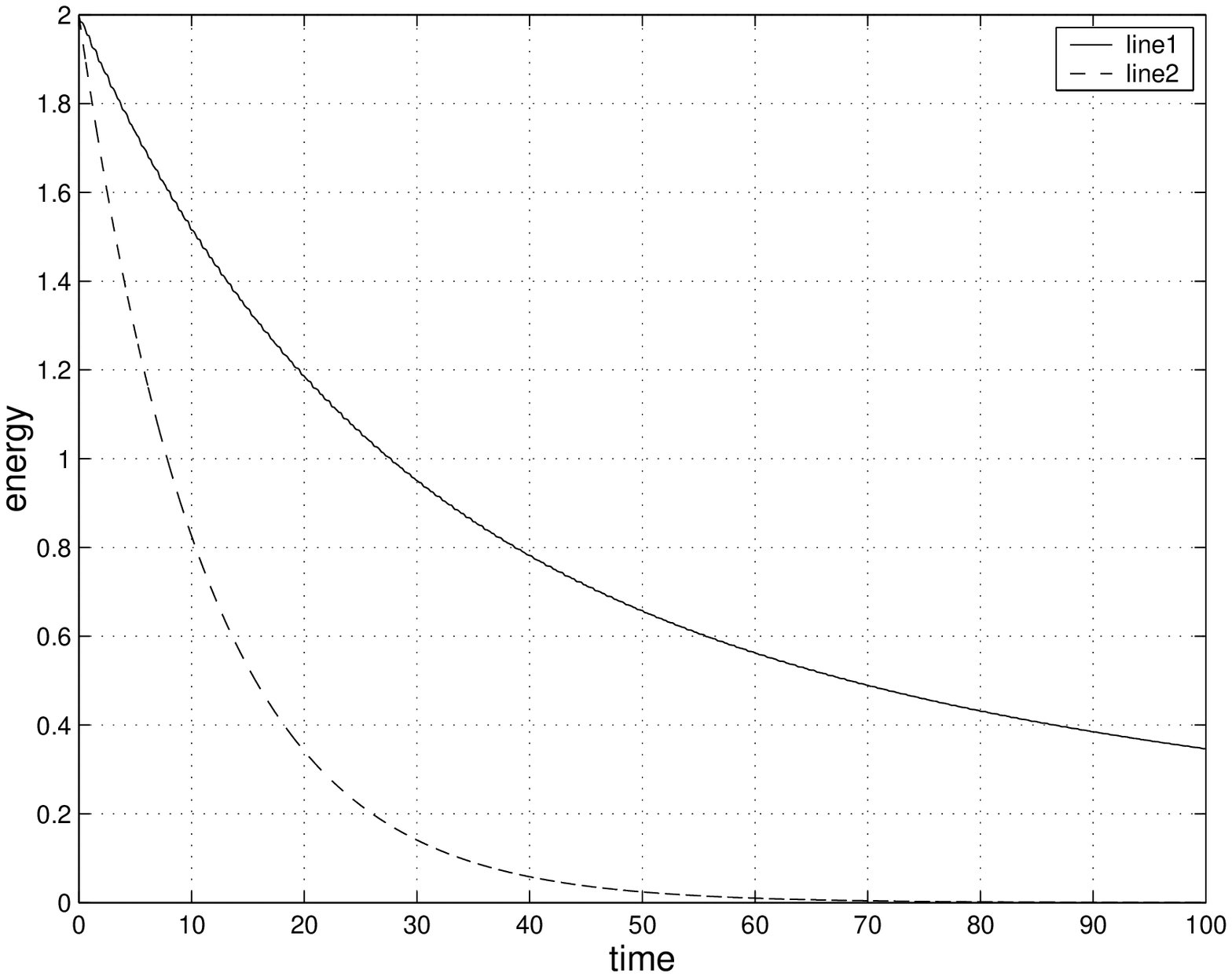}
\includegraphics[angle=0,width=6.5cm]{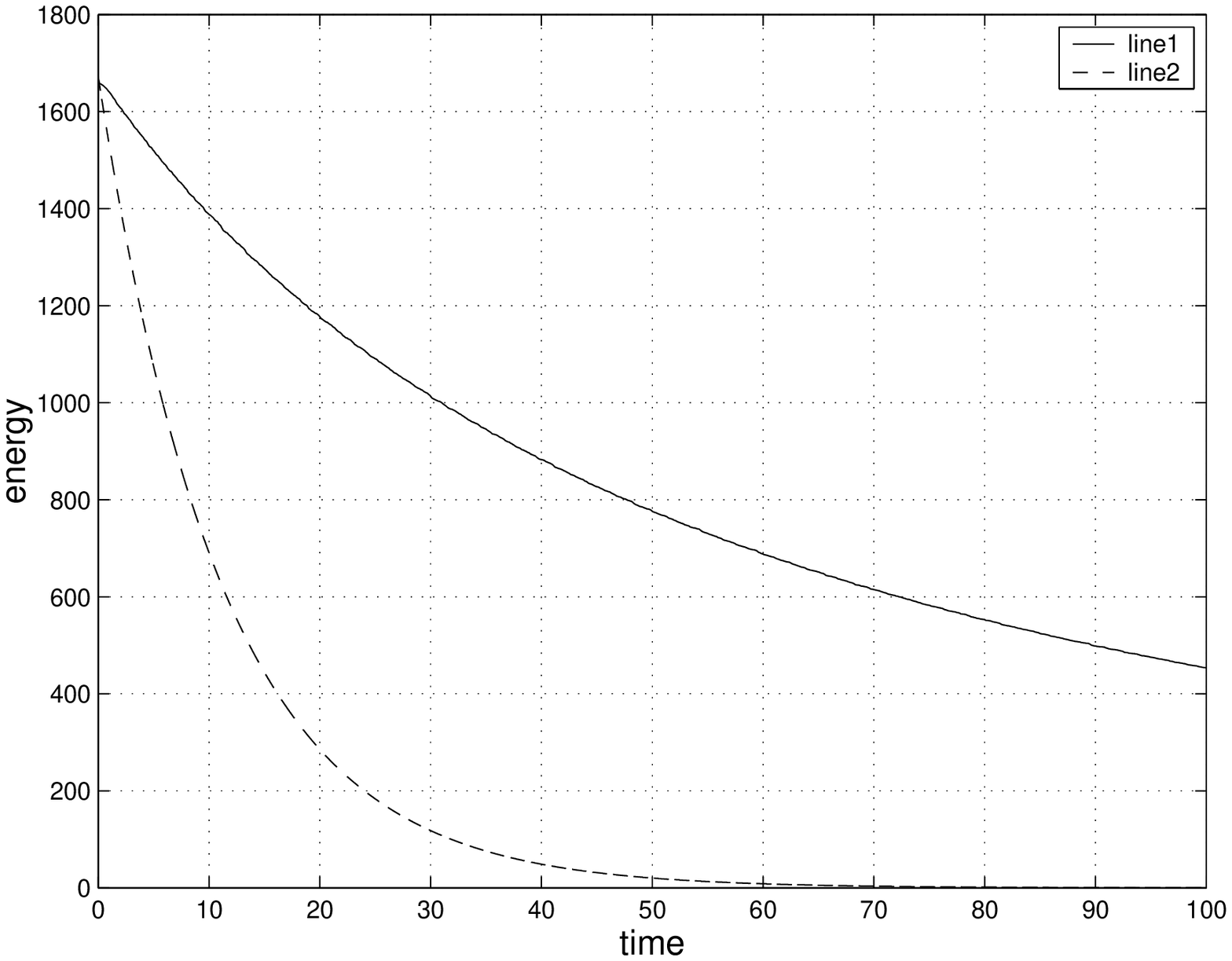}
\end{center}
\begin{center}
figure 30 : energy of the damped plate, $a(x)=1,\omega  = (0, \frac{1}{5}) \times  (0,1)$, $n=m=3$ (left),  and $n=m=10$ (right).
\end{center}


\begin{center}
\includegraphics[angle=0,width=6.5cm]{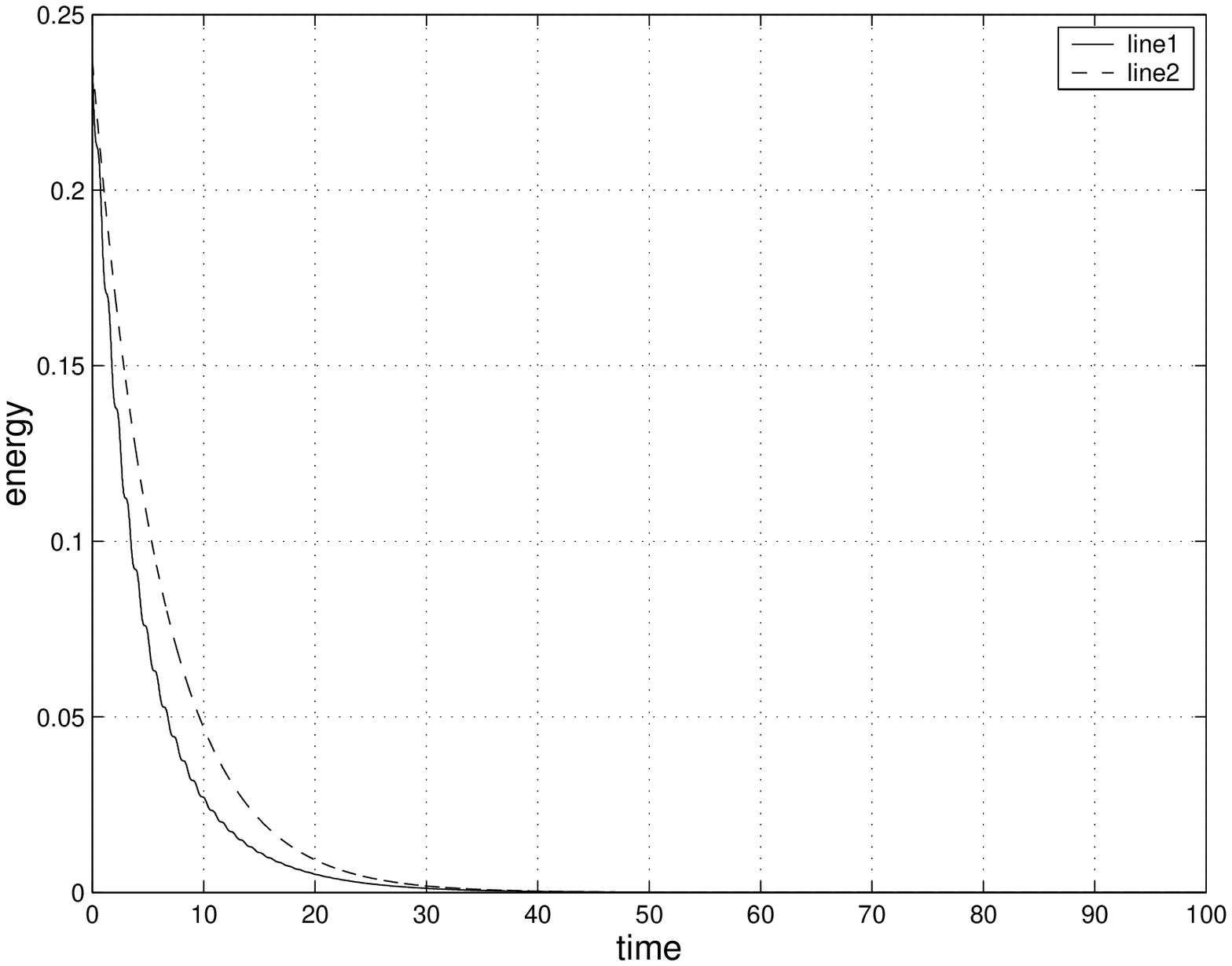}
\includegraphics[angle=0,width=6.5cm]{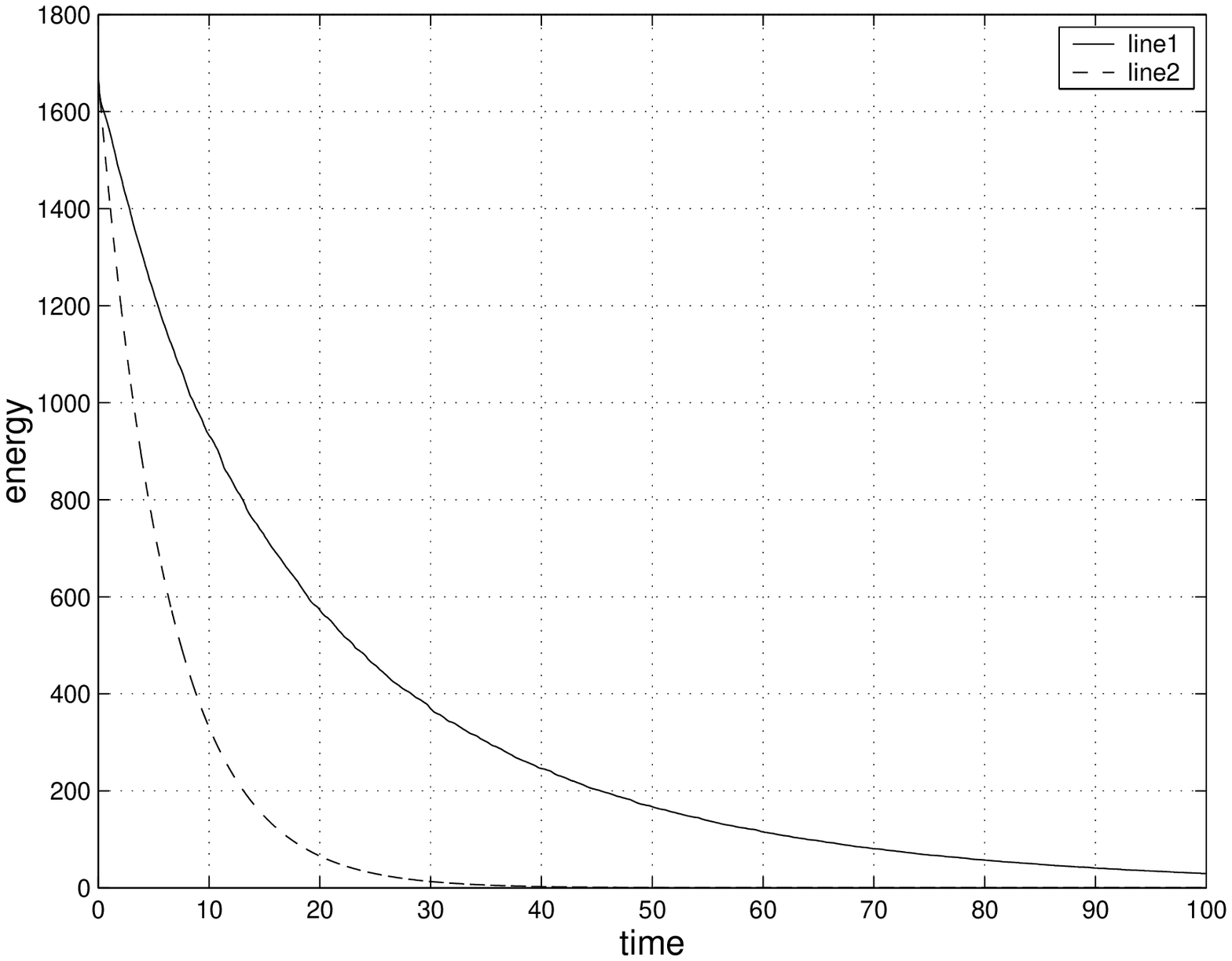}
\end{center}
\begin{center}
figure 31 : energy of the damped plate, $a(x)=1,\omega  = (0, 0.6) \times  (0,1)$, $n=m=2$ (left),  and $n=m=10$ (right).
\end{center}

\section{Optimization of the position of the damped region $\omega$ }

Starting from a fixed damping domain $\omega$ we try to see the influence of the position on decay energy. In the first example we plot the energy for five different damping position :  $\omega=(0.,0.4)\times (0,0.4)$, $\omega=(0.3,0.7)\times (0.3,0.7)$ (dashed line), $\omega=(0.5,0.9)\times (0.4,0.8)$ , $\omega=(0.1,0.5)\times (0.5,0.9)$, $\omega=(0.3,0.7)\times (0,0.4)$. The best position in this case is the middle of the plate. 
\begin{center}
\includegraphics[angle=0,width=6.5cm]{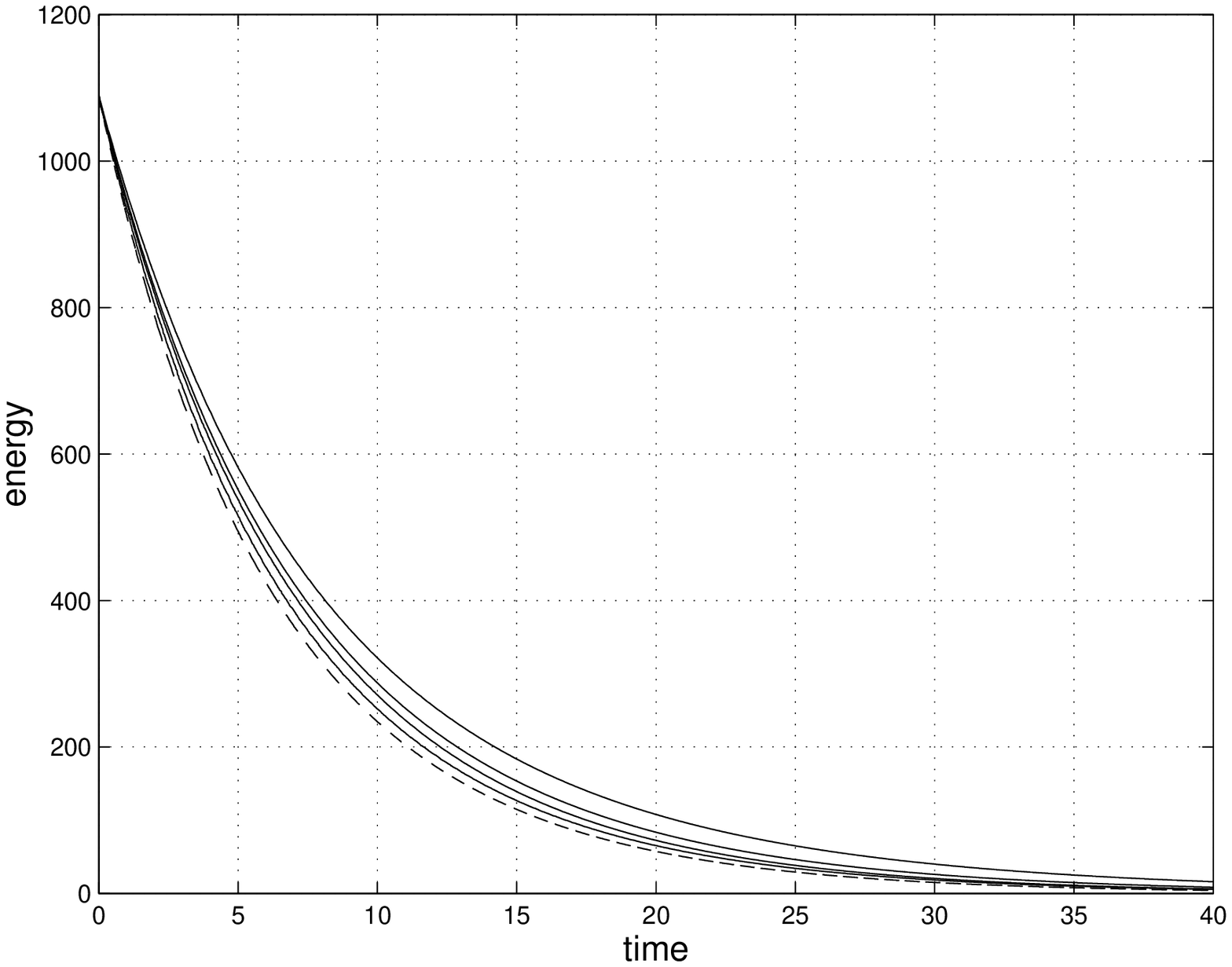}
\end{center}
\begin{center}
figure 32 : energy of the damped plate, $a(x)=1$, $n=m=8$.
\end{center}
In the second example we plot the energy for five damping region $\omega=(0.,0.3)\times (0,0.3)$ (green), $\omega=(0.35,0.65)\times (0.35,0.65)$ (red), $\omega=(0.1,0.4)\times (0.6,0.9)$ (blue) , $\omega=(0.7,1)\times (0.5,0.8)$ (black), $\omega=(0.7,1)\times (0.1,0.4)$, yellow. The best position in this case is the corner : $\omega=(0.,0.3)\times (0,0.3)$.
\begin{center}
\includegraphics[angle=0,width=6.5cm]{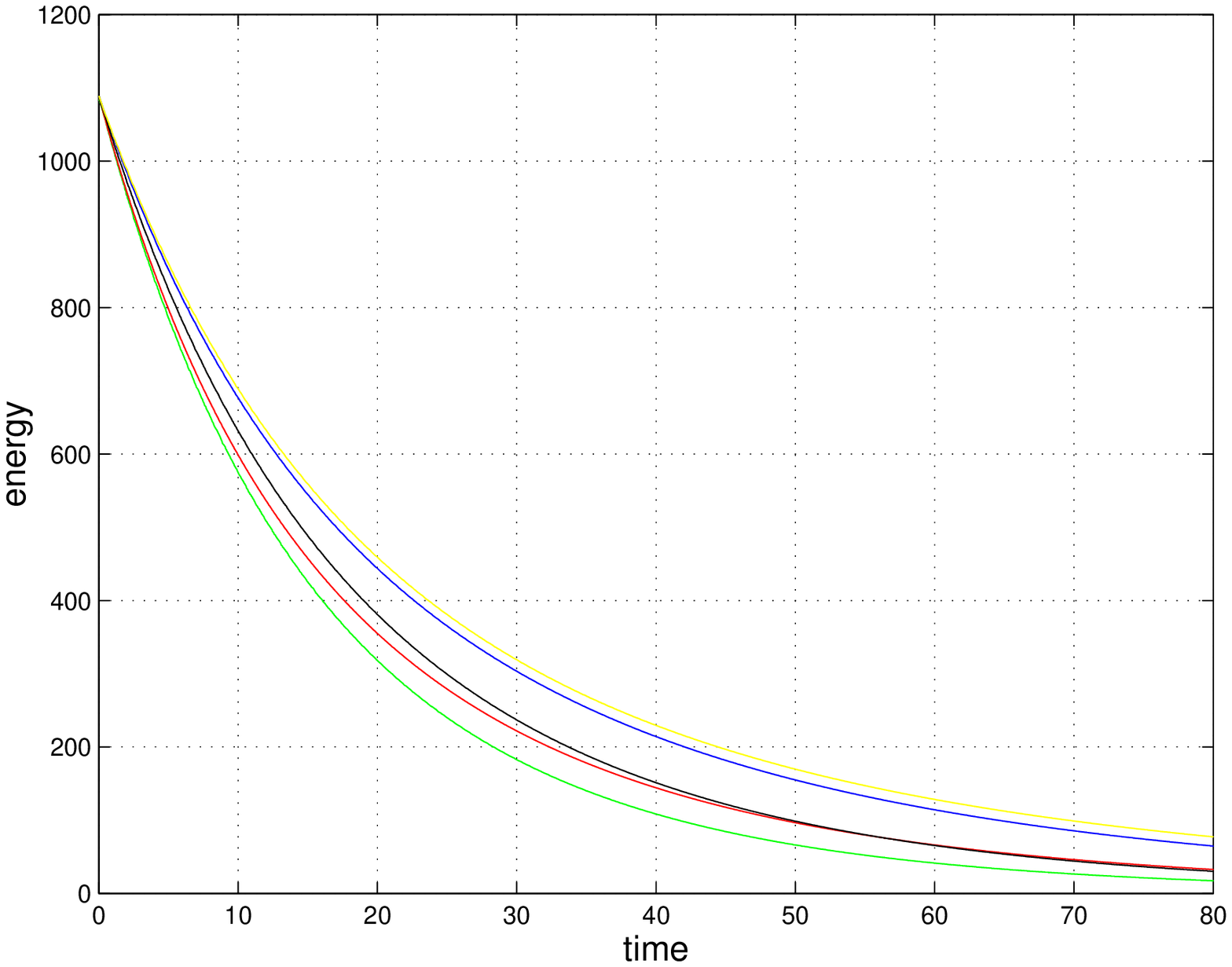}
\end{center}
\begin{center}
figure 33 : energy of the damped plate, $a(x)=1$, $n=m=8$.
\end{center}
\begin{center}
\includegraphics[angle=0,width=6.5cm]{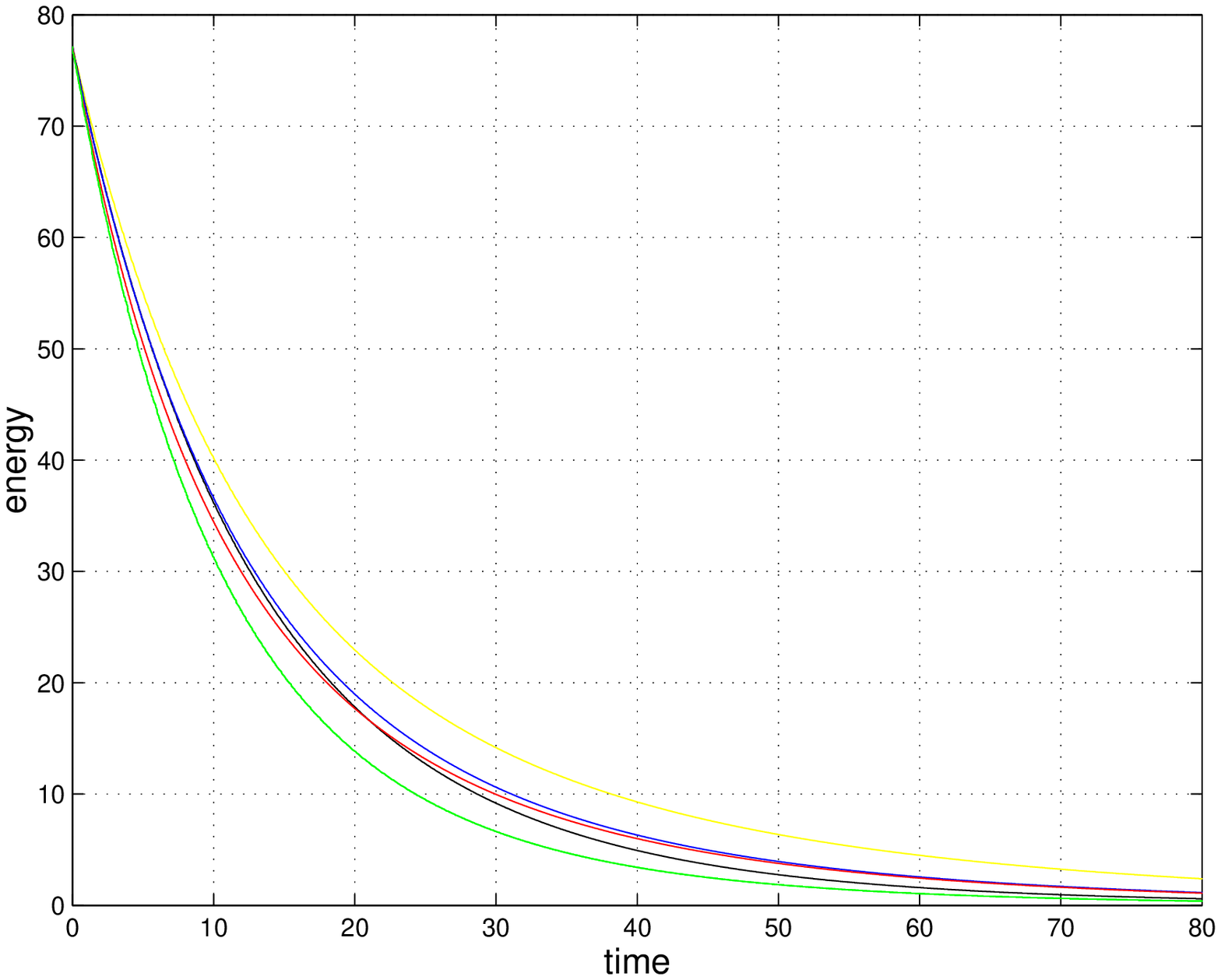}
\end{center}
\begin{center}
figure 34 : energy of the damped plate, $a(x)=1$, $n=m=5$.
\end{center}

In the third example we plot the energy for five damping region $\omega=(0.,0.2)\times (0,0.2)$ (black), $\omega=(0.1,0.3)\times (0.1,0.3)$ (green), $\omega=(0.2,0.4)\times (0.2,0.4)$ (blue) , $\omega=(0.3,0.5)\times (0.3,0.5)$ (yellow), $\omega=(0.4,0.6)\times (0.4,0.6)$, (red). for $n=m=5$,$n=m=6$ and $n=m=7$,
\begin{center}
\includegraphics[angle=0,width=6.5cm]{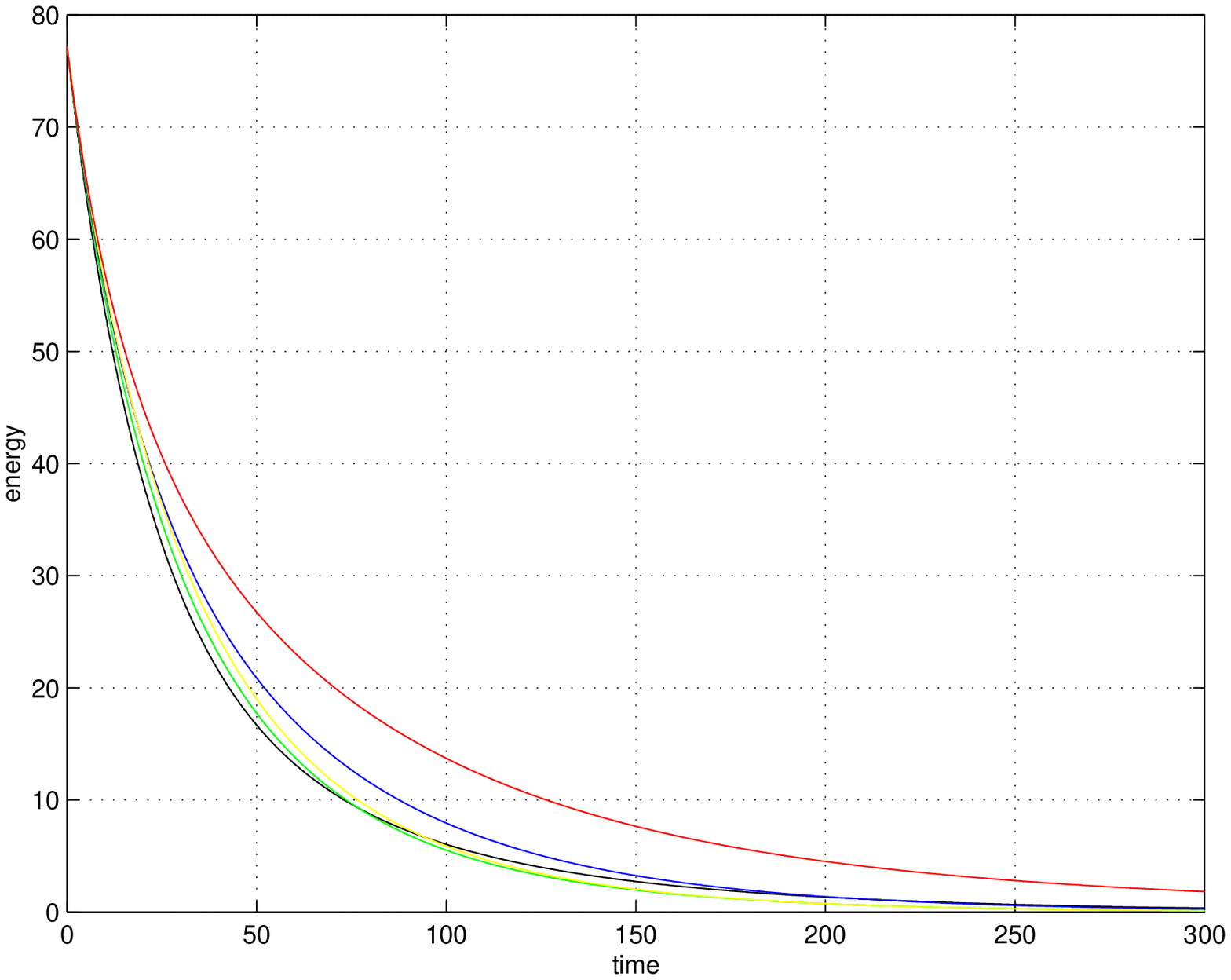}
\includegraphics[angle=0,width=6.5cm]{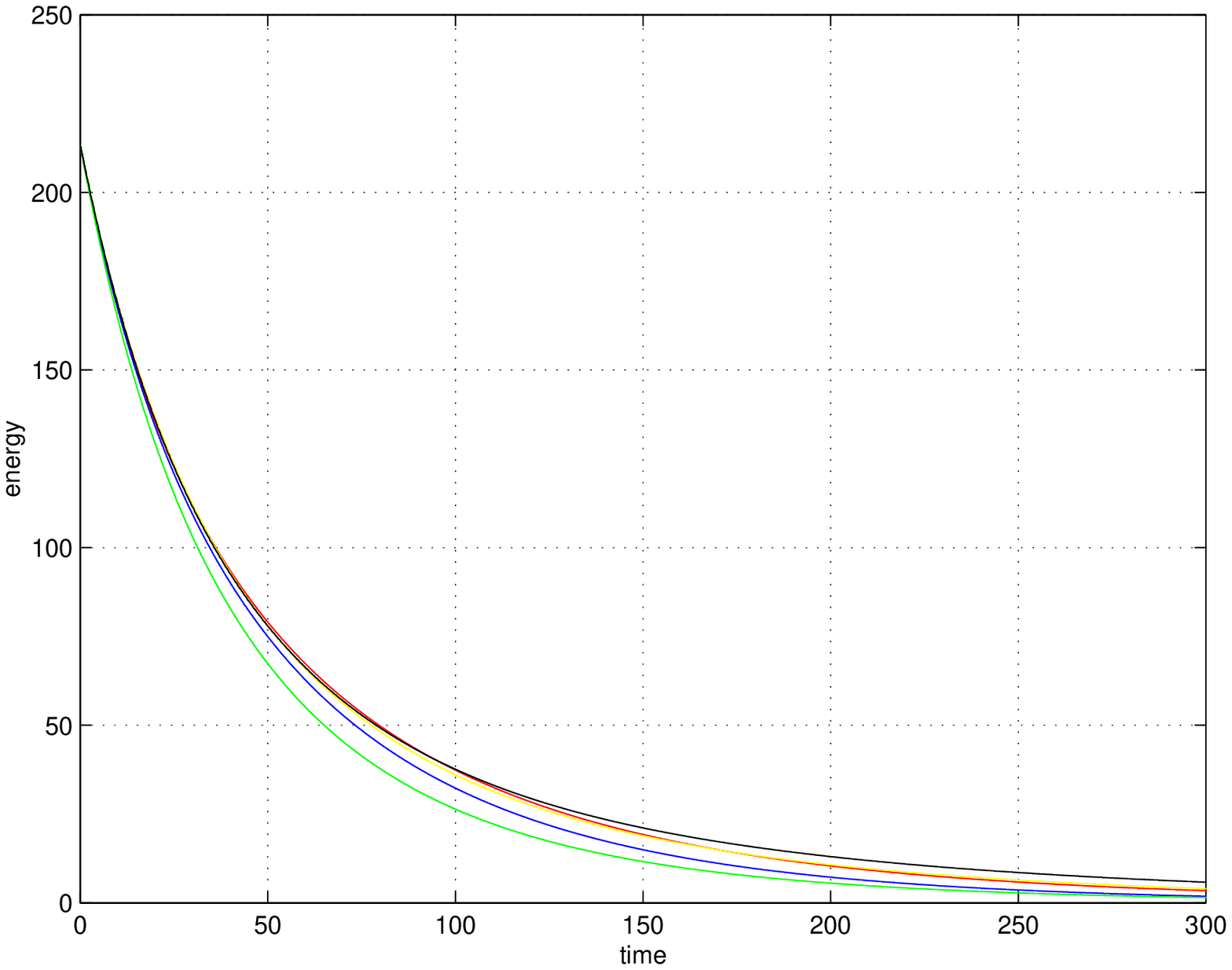}
\end{center}
\begin{center}
figure 35 : energy of the damped plate, $a(x)=1$, $n=m=5$ and $n=m=7$.
\end{center}
\begin{center}
\includegraphics[angle=0,width=6.5cm]{optim03a.ps}
\end{center}
\begin{center}
figure 36 : energy of the damped plate, $a(x)=1$, $n=m=7$.
\end{center}

\end{document}